\documentclass[final,1p,times,number]{elsarticle}

\usepackage{amsmath}
\usepackage{amsthm}
\usepackage{amssymb}

\usepackage{amsfonts}
\usepackage{graphicx}
\usepackage{epstopdf}
\usepackage{color}
\usepackage{bm}
\usepackage{multirow}
\usepackage{natbib}
\usepackage{hyperref}
\usepackage{cleveref}

\usepackage{color}
\usepackage{float}
\usepackage{caption}
\usepackage{booktabs}

\usepackage{algorithm}
\usepackage{algpseudocode}

\usepackage{geometry}

\usepackage[caption=false]{subfig}

\usepackage{ifpdf}
\ifpdf%
\usepackage{pdflscape}
\else
\usepackage{lscape}
\fi

\hypersetup{
	colorlinks,
	linkcolor={red!50!black},
	citecolor={blue!50!black},
	urlcolor={blue!80!black}
}

\def\f{\mbox{\boldmath $f$}}

\def\m{\mbox{\boldmath $m$}}

\def\0{\mbox{\boldmath $0$}}

\biboptions{sort&compress}

\begin{document}

\begin{frontmatter}
% \title{Comparative Study of Iterative Linear Solvers for Discrete Laudau-Lifshitz Systems}

% \title{Efficient Multigrid Method for Semi-implicit Schemes of The Landau-Lifshitz Model with Applications to Micromagnetics Simulations}

\title{Efficient Multigrid Methods for Semi‑implicit Landau‑Lifshitz Schemes in Micromagnetic Simulations}

\author[XJTLU]{Changjian Xie\corref{cor2}}
\cortext[cor2]{Corresponding author.} 
\ead{Changjian.Xie@xjtlu.edu.cn}

\address[XJTLU]{School of Mathematics and Physics, Xi'an-Jiaotong-Liverpool University, Re'ai Rd. 111, Suzhou, 215123, Jiangsu, China.}

\begin{abstract}

% This work presents an efficient aggregation multigrid method to solve the linear discrete system stemming from the first-order backward differentiation formula discretization of the Landau–Lifshitz equation, which describes the magnetization dynamics in micromagnetic simulations. Standard iterative and conventional multigrid solvers exhibit mesh-dependent convergence degradation and poor efficiency for such a linear system, as existing algorithms fail to accommodate its unique sparse spectral and asymmetric matrix structure. To resolve such issues, we employ a multigrid framework with refined smoothing and coarse-grid correction operators specifically adapted to the linear system. The proposed method confirms its unconditional stability and mesh-independent convergence. Numerical benchmarks verify that the developed multigrid solver prominently reduces iteration counts and computational cost compared with classic Krylov solvers and standard multigrid approaches, while maintaining accurate magnetization dynamic responses. 

An efficient aggregation‑based multigrid approach is developed in this work to solve linear algebraic systems generated by discretizing the Landau‑Lifshitz equation in micromagnetics. The discretization couples the first‑order backward differentiation formula with first‑order extrapolation, and the underlying equation describes magnetization dynamics within micromagnetic simulations. For these linear systems, standard iterative methods and conventional multigrid solvers intended for symmetric problems suffer from deteriorated convergence under mesh refinement and poor computational performance. Existing algorithms struggle to handle the system's intrinsic features, namely its sparse spectral properties and non‑symmetric matrix structure. To overcome these limitations, a multigrid framework is constructed, where smoothing and coarse‑grid correction operators are customized for the target linear system. Numerical tests confirm the robustness and mesh‑independent convergence of the resulting method. Compared with well‑established Krylov‑subspace solvers and conventional multigrid techniques, the aggregation‑based multigrid solver cuts down iteration counts and overall computational cost considerably, whilst producing faithful representations of magnetization dynamics.
\end{abstract}


\end{frontmatter}

\section{Introduction}

The Landau–Lifshitz (LL) equation is the canonical nonlinear partial differential equation that characterizes magnetization precession dynamics in ferromagnetic materials, serving as the core mathematical model \cite{Landau1935On} for micromagnetic analysis, spintronic device optimization \cite{vzutic2004spintronics}, and ultrafast magnetic dynamic simulation \cite{koopmans2005unifying}. Owing to its strong vector nonlinearity, nonconvex unit-length magnetization constraint, and coupled spatial-temporal gradients, analytical solutions of the nonlinear LL equation are generally unavailable, rendering numerical discretization and iterative solving the only feasible research approach. In semi-implicit (or implicit) time-marching simulations, the overall computational efficiency and numerical fidelity are predominantly determined by the performance of linear solvers for the resulting discrete linear systems. Interested readers are referred to \cite{cimrak2007survey,xie2026semi,cheng2023length,huang2026new} for further details on numerous semi‑implicit schemes for LL equation, including finite‑difference and finite‑element methods.

From the perspective of stability analysis, explicit schemes impose stringent constraints on the time‑step size. To obtain stable numerical solutions, an extremely small time step must be adopted, which gives rise to substantial overhead and various difficulties for long‑time dynamical evolution. For this reason, semi‑implicit and fully implicit formulations are favoured owing to their favourable stability properties. Nevertheless, fully implicit discretisations yield nonlinear algebraic systems. A practical remedy is to perform linearisation at each stage, which generates large non‑symmetric linear systems that require repeated solves within every iterative loop, leading to prohibitive computational costs. For instance, the Crank‑Nicolson scheme employed in \cite{jeong2010crank} produces nonlinear discrete systems, where multigrid techniques are utilised to accelerate computation. Semi‑implicit schemes for solving the Landau‑Lifshitz equation have attracted considerable attention due to their promising theoretical and practical characteristics, including robustness, stability and convergence behaviour. Different from linear systems with constant coefficients, the resulting linear systems from such semi‑implicit discretisations are large‑scale and feature non‑symmetric coefficient matrices; consequently, fast Fourier transform‑based solvers cannot be applied. Meanwhile,
\cite{jeong2014accurate} propose a numerical method which achieves second‑order accuracy in both space and time. A nonlinear multigrid approach is adopted to resolve nonlinear terms at each time step. Furthermore, the second‑order convergence in space‑time and the inherent energy conservation or dissipation property are validated.

In micromagnetic simulations, extremely fine spatial meshes are typically required, which gives rise to large‑scale sparse non‑symmetric linear systems. For our proposed semi‑implicit method, numerical experiments using GMRES incur substantial computational overhead and demand a large number of solver iterations, rendering the approach practically infeasible in certain scenarios. Classical numerical solving strategies adopted in micromagnetic simulations suffer from notable deficiencies when applied to the BDF1-discretized LL linear system. Traditional Krylov subspace iterative solvers for non-symmetric linear system of equations, such as BiCGstab or GMRES, TFQMR, typically exhibit mesh-dependent convergence degradation, where iteration counts grow sharply with grid refinement, leading to prohibitive computational costs for large-scale fine-grid simulations and a preconditioner need to be used \cite{pearson2020preconditioners}. For linear systems arising from discretized PDEs, efficient iterative solvers are often schemes of the multigrid or the multilevel type \cite{hackbusch1985multigrid,trottenberg2001multigrid}. These combine the effects of a smoother and of a coarse grid correction. The smoother is generally based on a simple iterative scheme such as the Gauss-Seidel method. The coarse grid correction consists in computing an approximate solution to the residual equation on a coarser grid which has fewer unknowns. The approach is applied recursively until the coarse system is small enough to make negligible the cost of an exact solution. Standard geometric and algebraic multigrid methods, though capable of accelerating long-wavelength error elimination for general linear systems, but not to be considered to the semi-implicit schemes for the LL dynamics. To address the aforementioned research gaps and numerical defects, this work proposes and employs an aggregation algebraic multigrid method (AGMG) designed for the linear discrete system arising from BDF1 temporal discretization of the Landau–Lifshitz equation and applies it to the micromagnetics simulations. 

The rest of this paper is structured as follows. \Cref{sec: model} presents the mathematical formulation of the LL equation and \Cref{sec:method} elaborates the detailed BDF1 temporal discretization procedure, deriving the final linear discrete system to be solved in \Cref{sec:system}. \Cref{sec:mgmethod} illustrates the fundamental multigrid principle and details the customized algorithm design, including the coarsening of aggregation-based algebraic multigrid method and the solution algorithm for the asymmetric linear system of equations. \Cref{sec:experiments} conducts comprehensive numerical experiments and comparative analysis to quantitatively demonstrate the efficiency and accuracy advantages of the proposed method. Finally, \Cref{sec:conclusions} concludes the entire work and prospects future research directions on the magnetic dynamic systems.

\section{The physical governing equation and main linear system}\label{sec: model}

In this paper, we first kick off studying the Landau-Lifshitz equation without damping as below,
\begin{equation}\label{eq:LL_undamped}
\partial_t \boldsymbol{m} = -\boldsymbol{m} \times \Delta \boldsymbol{m},
\end{equation}
which can be derived from microscopic modeling via continuumization under the semiclassical limit
The model with damping term is given as below,
\begin{equation}\label{eq:LL_damped}
\partial_t \boldsymbol{m} = -\boldsymbol{m} \times \Delta \boldsymbol{m}-\alpha\boldsymbol{m} \times (\boldsymbol{m} \times \Delta \boldsymbol{m}),
\end{equation}
where the damping term is added by the physical phenomenological experiments.
Let $\Omega\subset\mathbb{R}^d$ ($d=2,3$) be a bounded domain with smooth boundary $\partial\Omega$. The unknown $\boldsymbol{m}(\boldsymbol{x},t): \Omega\times(0,T)\to\mathbb{R}^3$ denotes the magnetization field subject to the unit-length constraint
\begin{equation}
|\boldsymbol{m}(\boldsymbol{x},t)| = 1,\quad \forall \boldsymbol{x}\in\Omega,\; t\in[0,T].
\end{equation}
Here $\partial_t$ represents the temporal partial derivative, $\Delta$ is the Laplace operator acting componentwisely on $\boldsymbol{m}$, and $\times$ stands for the vector cross product in $\mathbb{R}^3$. This equation describes purely conservative magnetization precession without Gilbert damping, and the total exchange energy
\[
E(\boldsymbol{m})=\frac12\int_\Omega |\nabla\boldsymbol{m}|^2\,\mathrm{d}\boldsymbol{x}
\]
is conserved along smooth solutions, i.e., $\frac{\mathrm{d}}{\mathrm{d}t}E(\boldsymbol{m})=0$. We impose suitable boundary conditions and initial condition
\[
\boldsymbol{m}(\boldsymbol{x},0)=\boldsymbol{m}_0(\boldsymbol{x}),\quad |\boldsymbol{m}_0(\boldsymbol{x})|=1,\quad \boldsymbol{x}\in\Omega.
\]
For the real micromagnetics simulations, we take the model in the nondimensionalized form as below,
\begin{equation}
\boldsymbol{m}_t = -\boldsymbol{m} \times (\epsilon\Delta\boldsymbol{m}+\boldsymbol{f}) -\alpha \boldsymbol{m} \times \boldsymbol{m} \times (\epsilon\Delta\boldsymbol{m}+\boldsymbol{f}),
\label{eq:llg_original}
\end{equation}
where the following source term is defined
\begin{equation}
\boldsymbol{f} = -q\big(m_2\boldsymbol{e}_2 + m_3\boldsymbol{e}_3\big) + \boldsymbol{h}_s + \boldsymbol{h}_e,
\label{eq:source_f}
\end{equation}
where $\boldsymbol{h}_e$ is the external field and $\boldsymbol{h}_s$ is the stray field with the following formula
\begin{equation}
\boldsymbol{h}_s = \frac{1}{4\pi}\nabla \int_{\Omega} \nabla\left(\frac{1}{|\boldsymbol{x}-\boldsymbol{y}|}\right)\cdot \boldsymbol{m}(\boldsymbol{y})\,d\boldsymbol{y}.
\label{eq:stray_field}
\end{equation}

Here, the dimensionless parameters become $\epsilon = C_{ex}/(\mu_0 M_s^2 L^2)$ and $q = K_u/(\mu_0 M_s^2)$ with $L$ the diameter of the ferromagnetic body and $\mu_0$ the permeability of vacuum. The unit vectors are given by $\boldsymbol{e}_2=(0,1,0)$, $\boldsymbol{e}_3=(0,0,1)$, and $\Delta$ denotes the standard Laplacian operator. For the Permalloy, an alloy of Nickel (80\%) and Iron (20\%), typical values of the physical parameters are given by: the exchange constant $C_{ex}=1.3\times 10^{-11}\,\mathrm{J/m}$, the anisotropy constant $K_u=100\,\mathrm{J/m^3}$, the saturation magnetization constant $M_s=8.0\times 10^{5}\,\mathrm{A/m}$. If $\Omega$ is a rectangular domain, the evaluation of \eqref{eq:stray_field} can be efficiently done by the Fast Fourier Transform (FFT).

\subsection{BDF1 temporal discretization}\label{sec:method}
We adopt the first-order backward differentiation formula (BDF1, backward Euler) for time discretization. Denote the uniform time step size by $\Delta t>0$ and write $t^n=n\Delta t$. Let $\boldsymbol{m}^n$ approximate $\boldsymbol{m}(\cdot,t^n)$. The BDF1 discretization of \eqref{eq:LL_undamped} reads
\begin{equation}\label{eq:BDF1_nonlinear}
\frac{\boldsymbol{m}^{n+1}-\boldsymbol{m}^n}{\Delta t} = -\boldsymbol{m}^{n+1}\times\Delta \boldsymbol{m}^{n+1}.
\end{equation}
Equation \eqref{eq:BDF1_nonlinear} is nonlinear due to the cross-product term. Following a standard semi-implicit linearization strategy, we freeze the magnetization vector on the right-hand side at the known previous time level $\boldsymbol{m}^n$, which yields the linearized time-discrete equation
\begin{equation}\label{eq:BDF1_linearized}
\boldsymbol{m}^{n+1} + \Delta t\,\boldsymbol{m}^n\times \Delta \boldsymbol{m}^{n+1} = \boldsymbol{m}^n.
\end{equation}

In this section, we study the first order fully discrete semi-implicit schemes as below,
\begin{align*}
    \frac{\m_h^{n+1}-\m_h^n}{\Delta t}=-\m_h^n \times \Delta_h \m_h^{n+1}.
\end{align*}
We have the linear system of equations as below,
\begin{align*}
    (I+\Delta t\m_h^n \times \Delta_h)\m_h^{n+1}=\m_h^n,
\end{align*}
which is
\begin{align*}
    m_1^{n+1}+\Delta t (m_2^n\cdot \Delta_h m_3^{n+1}-m_3^n \cdot\Delta_hm_2^{n+1})&=m_1^n\\
     m_2^{n+1}+\Delta t (m_3^n\cdot \Delta_h m_1^{n+1}-m_1^n \cdot\Delta_hm_3^{n+1})&=m_2^n\\
      m_3^{n+1}+\Delta t (m_1^n\cdot \Delta_h m_2^{n+1}-m_2^n \cdot\Delta_hm_1^{n+1})&=m_3^n.
\end{align*}
For the component-wise, we have
\begin{align*}
    \begin{pmatrix}
        I&-\Delta t m_3^n\cdot \Delta_h &\Delta t m_2^n\cdot \Delta_h\\
        \Delta t m_3^n\cdot \Delta_h&I&-\Delta t m_1^n\cdot \Delta_h\\
        -\Delta t m_2^n\cdot \Delta_h&\Delta t m_1^n\cdot \Delta_h&I
    \end{pmatrix}\begin{pmatrix}
        m_1^{n+1}\\
        m_2^{n+1}\\
        m_3^{n+1}
    \end{pmatrix}=\begin{pmatrix}
        m_1^{n}\\
        m_2^{n}\\
        m_3^{n}
    \end{pmatrix}.
\end{align*}
% If we take the form below,
% \begin{align*}
%     A=\begin{pmatrix}
%         I&-\Delta t m_3^n\cdot \Delta_h &\Delta t m_2^n\cdot \Delta_h\\
%         \Delta t m_3^n\cdot \Delta_h&I&-\Delta t m_1^n\cdot \Delta_h\\
%         -\Delta t m_2^n\cdot \Delta_h&\Delta t m_1^n\cdot \Delta_h&I
%     \end{pmatrix},
% \end{align*}
% we then have
% \begin{align*}
%     A+A^T=2I.
% \end{align*}

% We note that 
% \begin{align*}
%     \begin{pmatrix}
%         I&-\Delta t m_3^n\cdot \Delta_h &\Delta t m_2^n\cdot \Delta_h\\
%         \Delta t m_3^n\cdot \Delta_h&I&-\Delta t m_1^n\cdot \Delta_h\\
%         -\Delta t m_2^n\cdot \Delta_h&\Delta t m_1^n\cdot \Delta_h&I
%     \end{pmatrix}=\begin{pmatrix}
%         I&\Delta t m_3^n\cdot \Delta_h &\Delta t m_2^n\cdot \Delta_h\\
%         \Delta t m_3^n\cdot \Delta_h&I&\Delta t m_1^n\cdot \Delta_h\\
%         \Delta t m_2^n\cdot \Delta_h&\Delta t m_1^n\cdot \Delta_h&I
%     \end{pmatrix}+\begin{pmatrix}
%         0&-2\Delta t m_3^n\cdot \Delta_h &0\\
%         0&0&-2\Delta t m_1^n\cdot \Delta_h\\
%         -2\Delta t m_2^n\cdot \Delta_h&0&0
%     \end{pmatrix}
% \end{align*}
The semi-implicit BDF1 scheme for the Landau-Lifshitz equation is given by
\begin{align*}
    \frac{\m_h^{n+1}-\m_h^n}{\Delta t}=-\m_h^n \times \Delta_h \m_h^{n+1}-\alpha \m_h^n\times (\m_h^n \times \Delta_h \m_h^{n+1}).
\end{align*}
Thus, the linear system leads to 
\begin{align*}
    (I+\Delta t\m_h^n \times \Delta_h+\alpha \Delta t \m_h^n \times (\m_h^n \times \Delta_h))\m_h^{n+1}=\m_h^n,
\end{align*}
which is
\begin{align*}
   m_1^n&= m_1^{n+1}+\Delta t (m_2^n\cdot \Delta_h m_3^{n+1}-m_3^n \cdot\Delta_hm_2^{n+1})+\alpha\Delta t[-(m_2^2+m_3^2)\cdot\Delta_h m_1+m_1m_2 \cdot\Delta_h m_2+m_1m_3 \cdot\Delta_h m_3]\\
     m_2^n&=m_2^{n+1}+\Delta t (m_3^n\cdot \Delta_h m_1^{n+1}-m_1^n \cdot\Delta_hm_3^{n+1})+\alpha\Delta t[m_1m_2 \cdot\Delta_h m_1-(m_1^2+m_3^2)\cdot\Delta_h m_2 +m_2m_3\cdot \Delta_h m_3]\\
      m_3^n&=m_3^{n+1}+\Delta t (m_1^n\cdot \Delta_h m_2^{n+1}-m_2^n \cdot\Delta_hm_1^{n+1})+\alpha\Delta t[m_1m_3 \cdot\Delta_h m_1+m_2m_3 \cdot\Delta_h m_2-(m_1^2+m_2^2)\cdot\Delta_h m_3].
\end{align*}

For the full LLG model \eqref{eq:llg_original}, we have the scheme as below,
\begin{align*}
    \frac{\m_h^{n+1}-\m_h^n}{\Delta t}=-\m_h^n \times (\epsilon \Delta_h \m_h^{n+1}+\f_h^n)-\alpha \m_h^n\times (\m_h^n \times (\epsilon \Delta_h \m_h^{n+1}+\f_h^n)).
\end{align*}

\subsection{The resulting linear algebraic system}\label{sec:system}
After spatial discretization (finite element or finite difference) on a given computational mesh, the continuous linear partial differential equation \eqref{eq:BDF1_linearized} leads to a sparse linear system at each time step of the form
\begin{equation}\label{eq:linear_system}
\mathbf{A}_h^n\mathbf{u}_h^n=\mathbf{b}_h^n.
\end{equation}
Here $\mathbf{u}_h^n$ is the vector collecting discrete unknowns of $\boldsymbol{m}_h^{n+1}$, the coefficient matrix $\mathbf{A}_h^n$ arises from the combination of BDF1 temporal discretization and spatial discretization, and the right-hand side $\mathbf{b}_h^n$ is explicitly constructed from the known state $\boldsymbol{m}_h^n$. Here the notation $(\cdot)_h^n$ denotes the quantity associated with time discretization and spatial discretization.

\section{Multigrid method}\label{sec:mgmethod}

We consider the iterative solution of large sparse $n\times n$ linear systems \eqref{eq:linear_system}
% \[
% A\mathbf{u} = \mathbf{b}
% \]
arising from the discretization of the semi-implicit scheme by BDF1 for the LL equation. In this context, multigrid methods \cite{trottenberg2001multigrid} are among the most efficient solution techniques. Whereas geometric multigrid methods require a predetermined hierarchy of grids and discretizations, algebraic multigrid (AMG) methods are set up using only the information present in the system matrix \cite{brandt1984amg}.

These algorithms combine the effect of a \textit{smoother} and a \textit{coarse grid correction}. In AMG schemes, the smoother is fixed and generally based on a simple iterative method such as the (symmetric) Gauss-Seidel method. The coarse grid correction consists of computing an approximate solution to the residual equation on a coarser grid, that is, solving a linear system of smaller size. This solution is then transferred back to the actual grid by means of an appropriate prolongation. In AMG methods, this coarse grid correction is entirely defined once the prolongation is known, that is, once an appropriate prolongation matrix has been set up by applying a so-called \textit{coarsening algorithm} to the system matrix.

The improvement of AMG schemes is a hot research topic. The current trend (e.g., \cite{brezina2000algebraic,brezina2005adaptive,chartier2004spectral,henson2002elementfree,jones2001amge}) leads to more involved algorithms with denser prolongation matrices, which increases setup costs and memory requirements. In this paper, we consider coarsening by aggregation of the unknowns, which leads to prolongation matrices with at most one nonzero entry per row, which are much sparser than the ones obtained by the classical AMG approach, as developed in \cite{brandt1984algebraic,ruge1987algebraic,stuben1983algebraic}.

\subsection{Aggregation‑based algebraic multigrid method}

We use an aggregation‑based algebraic multigrid method introduced by \cite{notay2010aggAMG}. An algebraic multigrid method is presented to solve large systems of linear equations. The coarsening is obtained by aggregation of the unknowns. The aggregation scheme uses two passes of a pairwise matching algorithm applied to the matrix graph, resulting in most cases in a decrease of the number of variables by a factor slightly less than four. The matching algorithm favors the strongest negative coupling(s), inducing a problem dependent coarsening. This aggregation is combined with piecewise constant (unsmoothed) prolongation, ensuring low setup cost and memory requirements. Compared with previous aggregation‑based multigrid methods, the scalability is enhanced by using a so‑called K‑cycle multigrid scheme, providing Krylov subspace acceleration at each level.

\subsubsection{Coarsening}
An algebraic coarsening algorithm sets up a prolongation matrix $P$ using only the information available in $A$.
The prolongation is an $n\times n_c$ matrix, where $n_c < n$ is the number of coarse variables.
It allows to transfer on the fine grid a vector defined on the coarse variable set $[1,n_c]$.
Further, it entirely determines the coarse grid correction. Formally, the latter also depends on a restriction matrix and on a coarse grid matrix, but, as it is usual with AMG methods, one takes the restriction equal to the transpose of the prolongation, and the coarse grid matrix is computed from the Galerkin formula
\begin{equation}
A_c = P^T AP.
\label{eq:galerkin}
\end{equation}

In the classical AMG coarsening \cite{stueben2001introAMG}, one first selects a subset of fine grid variables as coarse variable, by inspecting the graph of $A$. Next, the matrix entries are used to build interpolation rules, that define $P$.

Coarsening by aggregation works differently. One needs to define aggregates $G_i$, which are disjoint subsets of the variable set. The number of coarse variables $n_c$ is then the number of such subsets, and $P$ is given by
\begin{equation}
P_{ij}=
\begin{cases}
1, & \text{if } i\in G_j,\\
0, & \text{otherwise},
\end{cases}
\quad (1\le i\le n,\, 1\le j\le n_c).
\label{eq:Pagg}
\end{equation}

If $\bigcup_i G_i = [1,n]$ (i.e., if the aggregates form a partitioning of $[1,n]$), $P$ is a Boolean matrix with exactly one nonzero entry per row. As seen below, it is however sometimes advantageous to leave some variables outside the set of aggregates, in which case the rows of $P$ corresponding to these variables are zero. Note that there is no need to explicitly form $P$, and the coarse grid matrix \eqref{eq:galerkin} is in practice computed by
\begin{equation}
(A_c)_{ij}= \sum_{k\in G_i}\sum_{\ell\in G_j} a_{k\ell},
\quad (1\le i,j\le n_c).
\label{eq:Ac_direct}
\end{equation}

Many aggregation algorithms proposed in the literature (e.g., \cite{braess1985towards,kim2003graphmatchingmg,notay2006aggregationamg}) starts by forming pairs, or ``matchings'', in the matrix graph. Here we reuse the algorithm from \cite{notay2006aggregationamg}, because it discriminates between different neighbors of a node, giving preference to the strongest negative coupling(s). For two-dimensional anisotropic model problems, this produces aggregates aligned with the strong coupling direction, as desired according to the analysis in \cite{muresan2008analysisaggregationmg}. Some connection may also be made with the classical AMG coarsening, which is also based on strong negative couplings.

\subsubsection{The solution algorithm}

In nonsymmetric cases, we use a preconditioned variant of GCR \cite{eisenstat1983variational}, referred to as GMRESR in \cite{vandervorst1994gmresr}. This method provides the minimal residual norm solution and allows for variable preconditioning. We use an improved implementation, given in \Cref{alg:precon_gcr_econ} for the sake of completeness. In the standard implementation, at step 3(b) one applies to $\mathbf z_j$ a recursion similar to the one applied to $\mathbf c_j^{(\cdot)}$. This allows us to obtain the approximate solution at each step with a simple recursion, but doubles the cost of step 3(b), which is the most expensive part of the algorithm. In \Cref{alg:precon_gcr_econ}, instead, the computation of the approximate solution is performed only upon completion of the main loop. Note that one has effectively
\[
\mathbf r_m = \mathbf b - A\mathbf u_m,
\]
because $\mathbf r_m = \mathbf r_0 - (\mathbf c_1 \ \cdots \ \mathbf c_m)\mathbf a$ whereas $A(\mathbf z_1 \ \cdots \ \mathbf z_m) = (\mathbf c_1 \ \cdots \ \mathbf c_m)\Gamma$, entailing $\mathbf r_m = \mathbf r_0 - A(\mathbf z_1 \ \cdots \ \mathbf z_m)(\Gamma^{-1}\mathbf a)$. This further shows that \Cref{alg:precon_gcr_econ} is mathematically equivalent to the original GCR/GMRESR algorithm, since the residual is computed as in the latter. This improved implementation is discussed in detail in \cite{jiranek2008how} for the unpreconditioned case, where it is shown to also have superior stability properties.

\begin{algorithm}[htbp]
\caption{Preconditioned GCR}
\label{alg:precon_gcr_econ}
\begin{algorithmic}[1]
\Require Matrix $A$; right-hand-side $\mathbf{b}$; initial approximation $\mathbf{u}_0$;
Maximal number of iterations $m$; tolerance $\varepsilon$.
\Ensure Approximate solution $\mathbf{u}_m$; residual $\mathbf{r}_m = \mathbf{b} - A\mathbf{u}_m$.
\State Initialization: $\mathbf{r}_0 = \mathbf{b} - A\mathbf{u}_0$.
\For{$j = 1,\dots, m$}
    \State 1. Apply preconditioner: $\mathbf{z}_j = \operatorname{Prec}(\mathbf{r}_{j-1})$.
    \State 2. $\mathbf{c}_j^{(1)} = A\mathbf{z}_j$.
    \For{$i = 1,\dots, j-1$}
        \State (a) $\gamma_{ij} = \mathbf{c}_i^T \mathbf{c}_j^{(i)}$,
        \State (b) $\mathbf{c}_j^{(i+1)} = \mathbf{c}_j^{(i)} - \gamma_{ij}\mathbf{c}_i$,
    \EndFor
    \State 4. $\gamma_{jj} = \|\mathbf{c}_j^{(j)}\|$;\quad $\mathbf{c}_j = \gamma_{jj}^{-1}\mathbf{c}_j^{(j)}$.
    \State 5. $\alpha_j = \mathbf{c}_j^T \mathbf{r}_{j-1}$;\quad $\mathbf{r}_j = \mathbf{r}_{j-1} - \alpha_j \mathbf{c}_j$.
    \State 6. \textbf{If} $\|\mathbf{r}_j\| < \varepsilon\|\mathbf{b}\|$, \textbf{exit do loop and reset} $m = j$.
\EndFor
\State $\mathbf{u}_m = \begin{bmatrix}\mathbf{z}_1 & \cdots & \mathbf{z}_m\end{bmatrix}\bigl(\Gamma^{-1}\mathbf{a}\bigr)$
\State where
\[
\mathbf{a}=
\begin{bmatrix}
\alpha_1\\
\vdots\\
\alpha_m
\end{bmatrix},
\qquad
\Gamma_{ij}=
\begin{cases}
\gamma_{ij}, & \text{if } i\le j,\\
0, & \text{otherwise}.
\end{cases}
\]
\end{algorithmic}
\end{algorithm}

% \section{Theoretical analysis}\label{sec:theory}

\section{Numerical experiments}
\label{sec:experiments}

\subsection{Spectral property and condition number}

\begin{figure}[htbp]
    \centering
    \subfloat[$Nx=64, N_t=10$]{\includegraphics[width=0.25\linewidth]{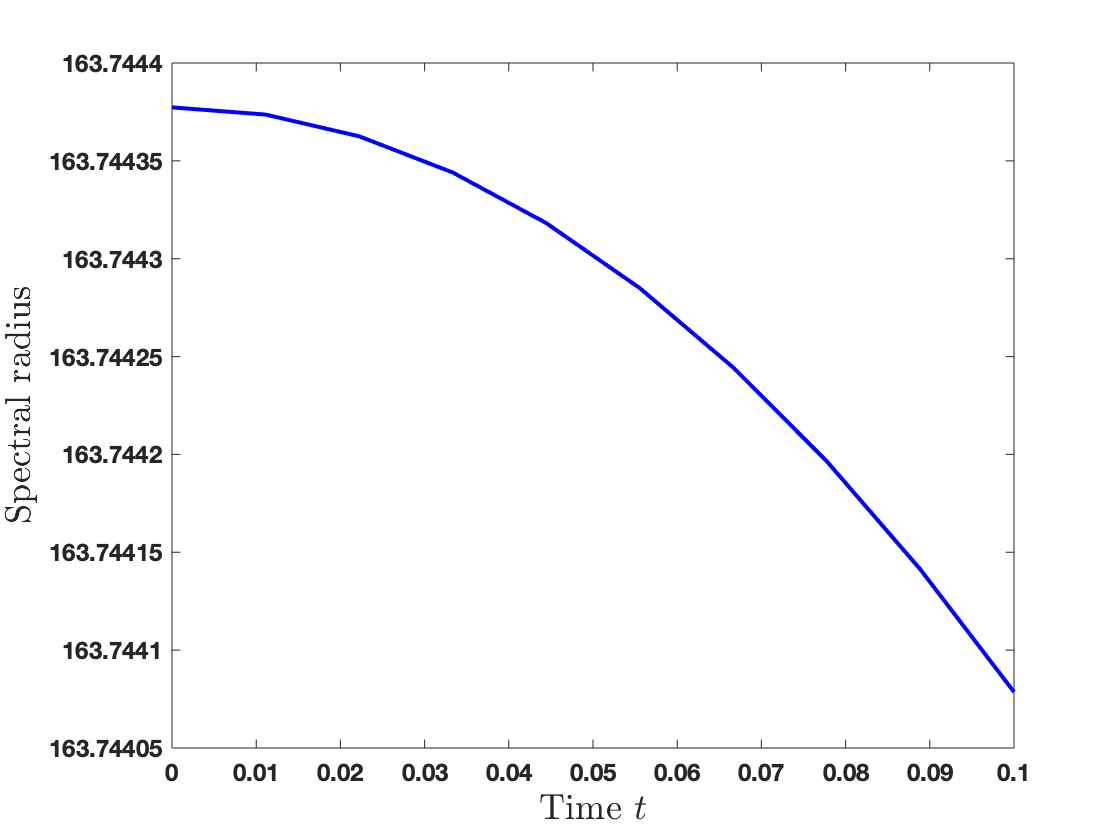}}
    \subfloat[$Nx=64, N_t=100$]{\includegraphics[width=0.25\linewidth]{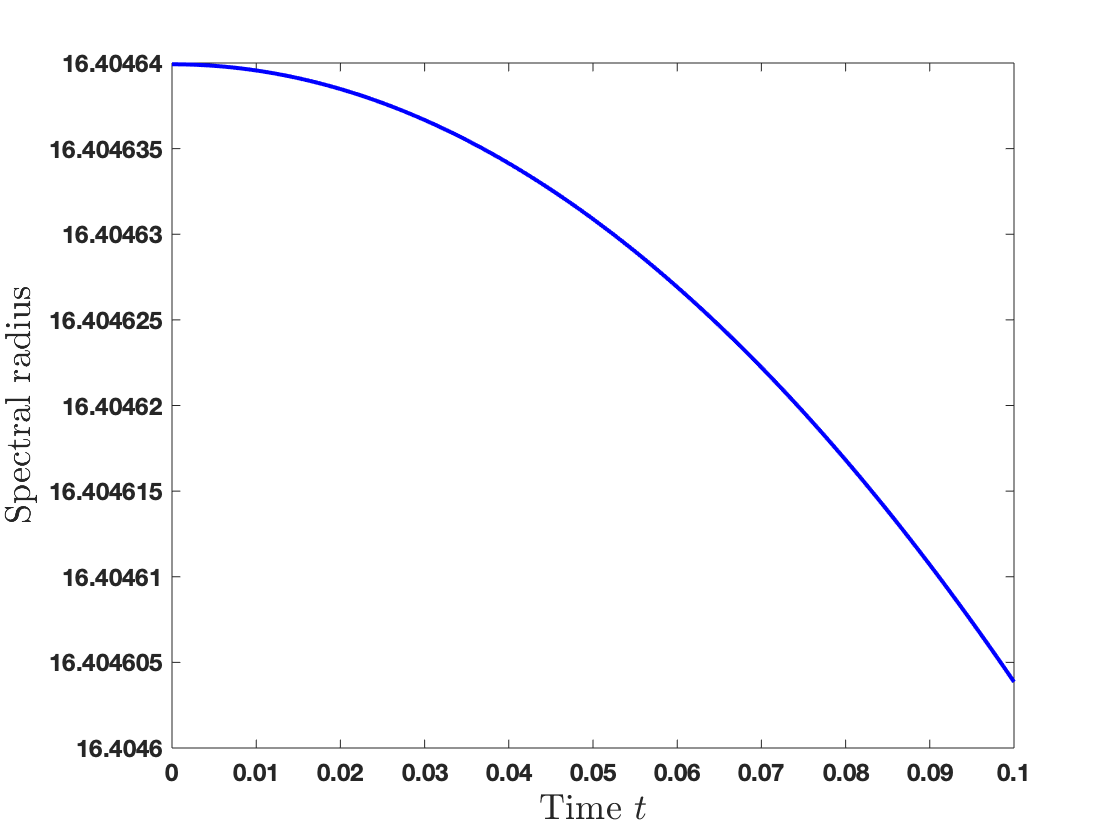}}
     \subfloat[$Nx=64, N_t=1000$]{\includegraphics[width=0.25\linewidth]{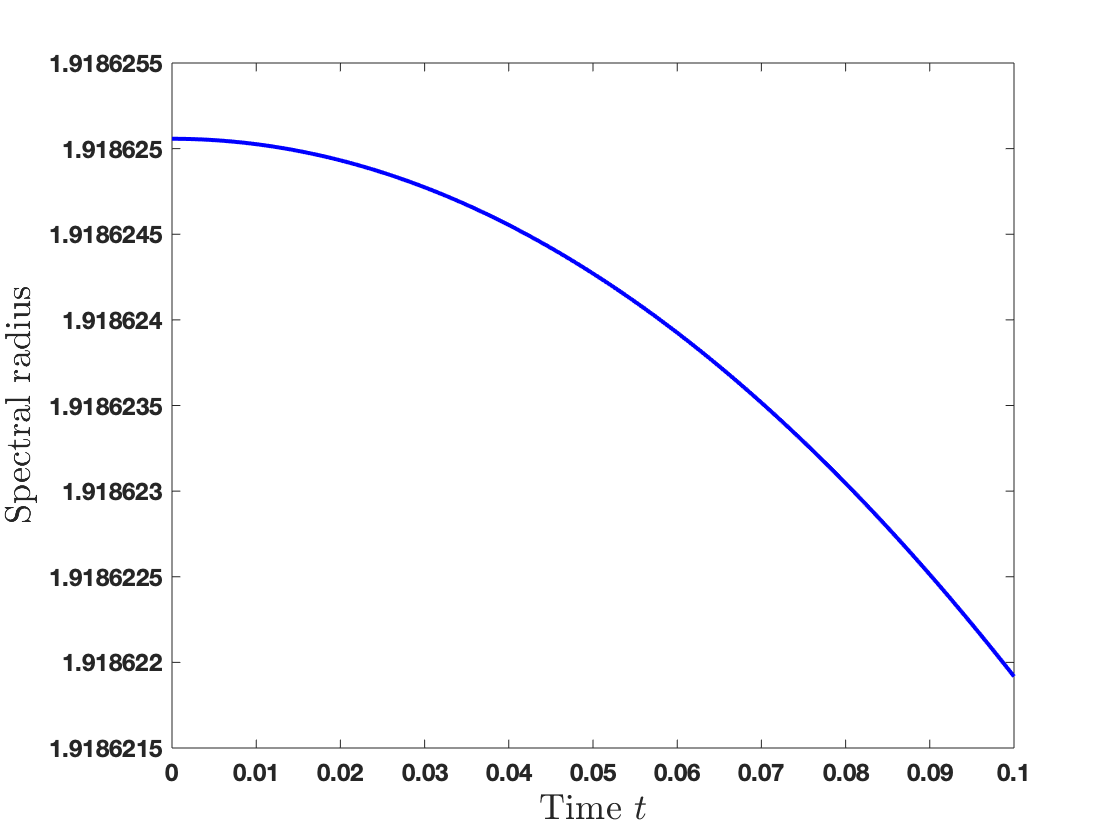}}
     \hspace{0.1in}
      \subfloat[$Nx=200, N_t=10$]{\includegraphics[width=0.25\linewidth]{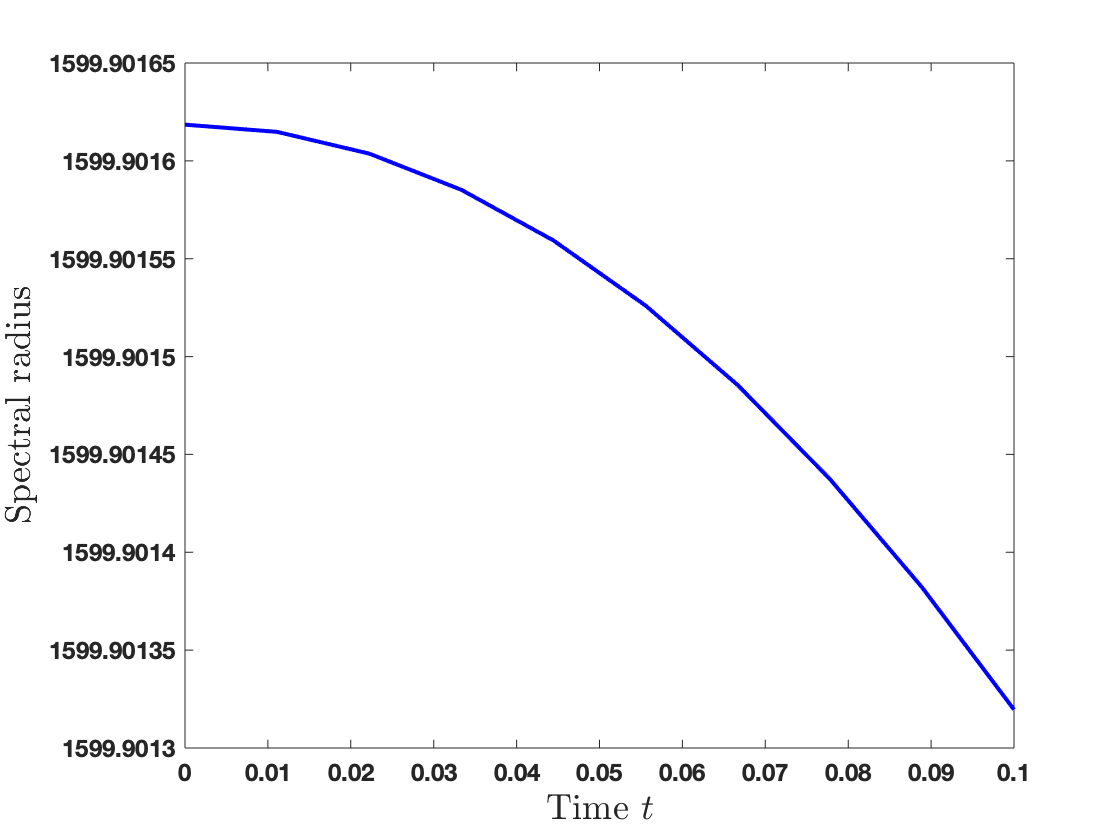}}
    \subfloat[$Nx=200, N_t=100$]{\includegraphics[width=0.25\linewidth]{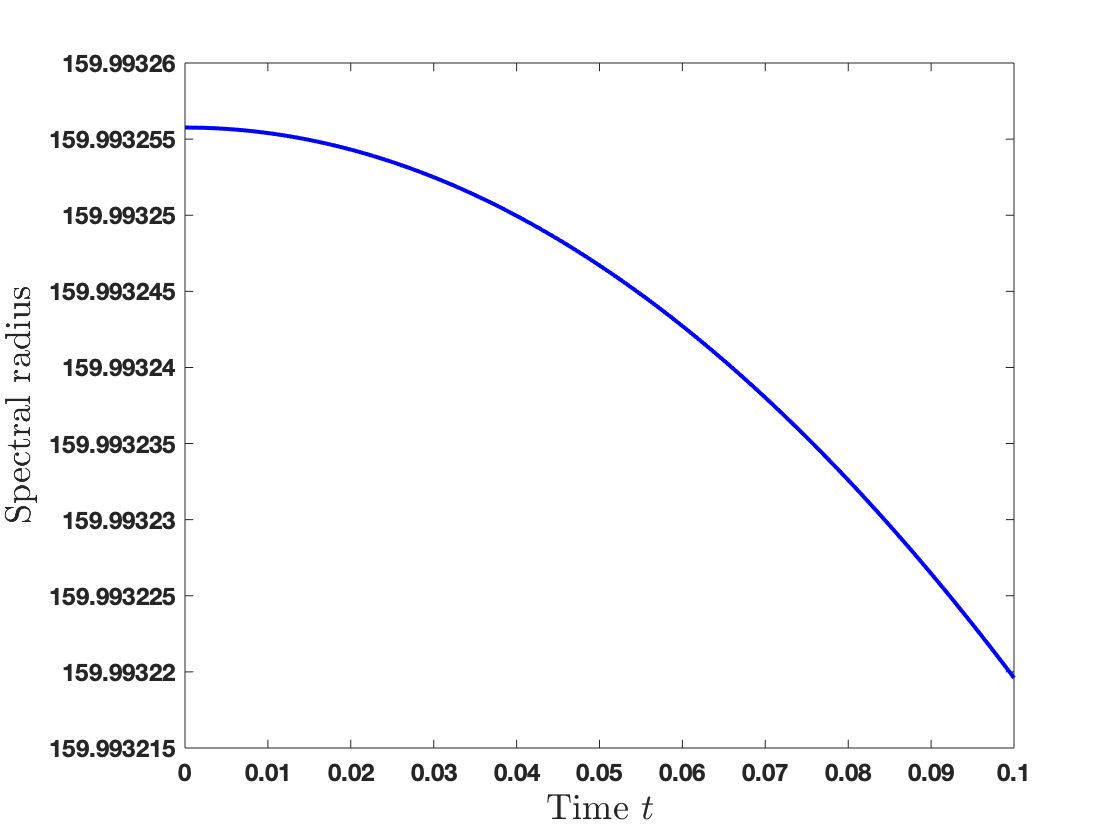}}
     \subfloat[$Nx=200, N_t=1000$]{\includegraphics[width=0.25\linewidth]{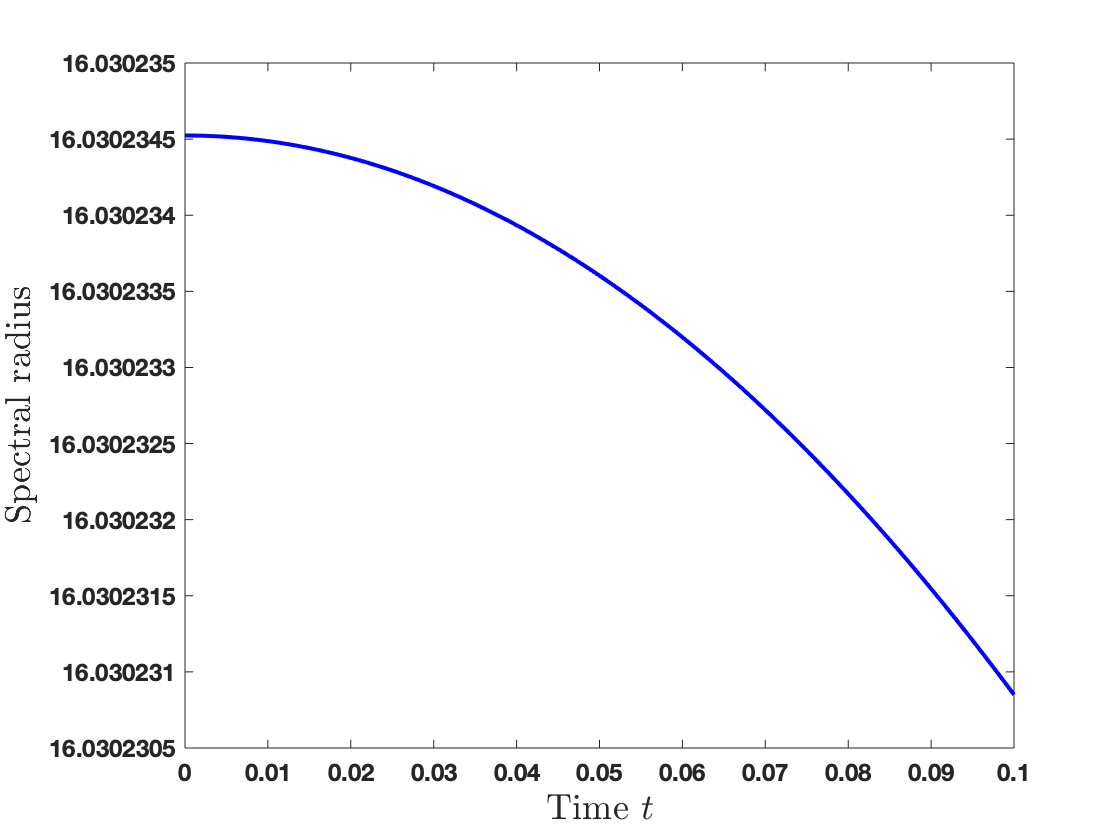}}
 \caption{The spectral radius of the matrix $A$ up to the final time $T=0.1$ for the 1D case. BDF1 without projection.}
    \label{fig:spectral_rad-1}
\end{figure}

\begin{figure}[htbp]
    \centering
    \subfloat[$Nx=64, N_t=10$]{\includegraphics[width=0.25\linewidth]{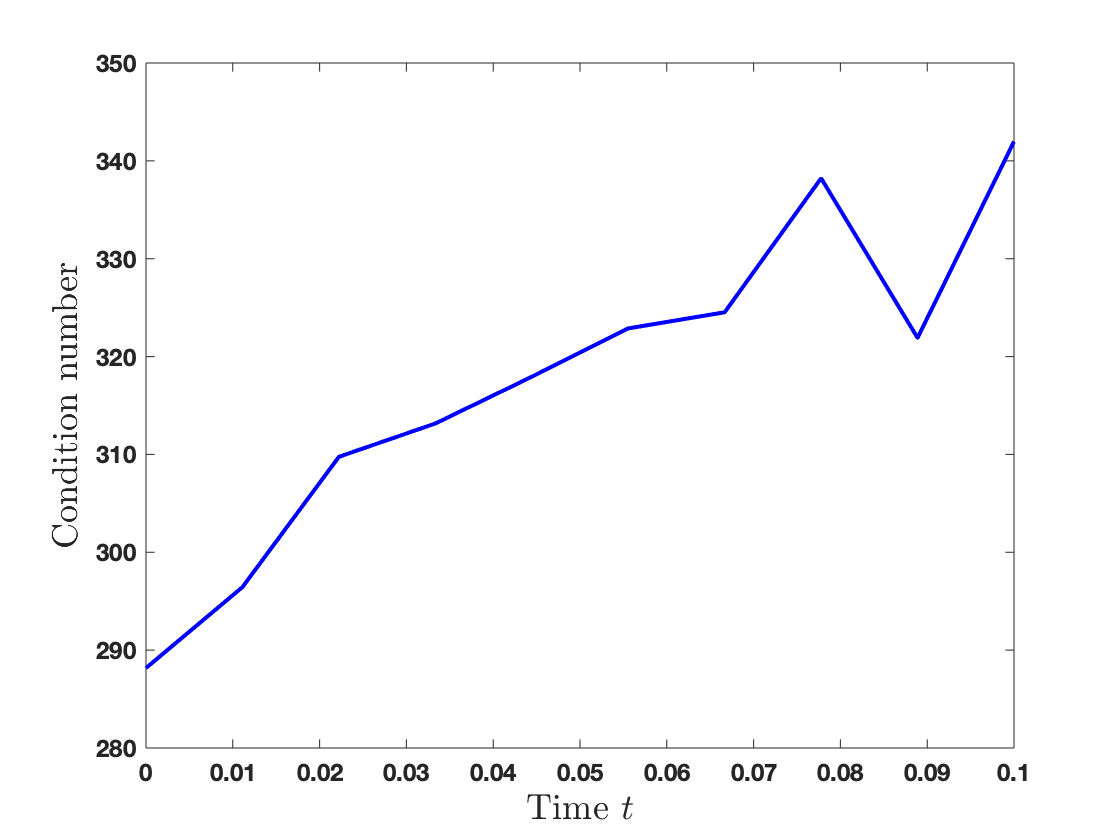}}
    \subfloat[$Nx=64, N_t=100$]{\includegraphics[width=0.25\linewidth]{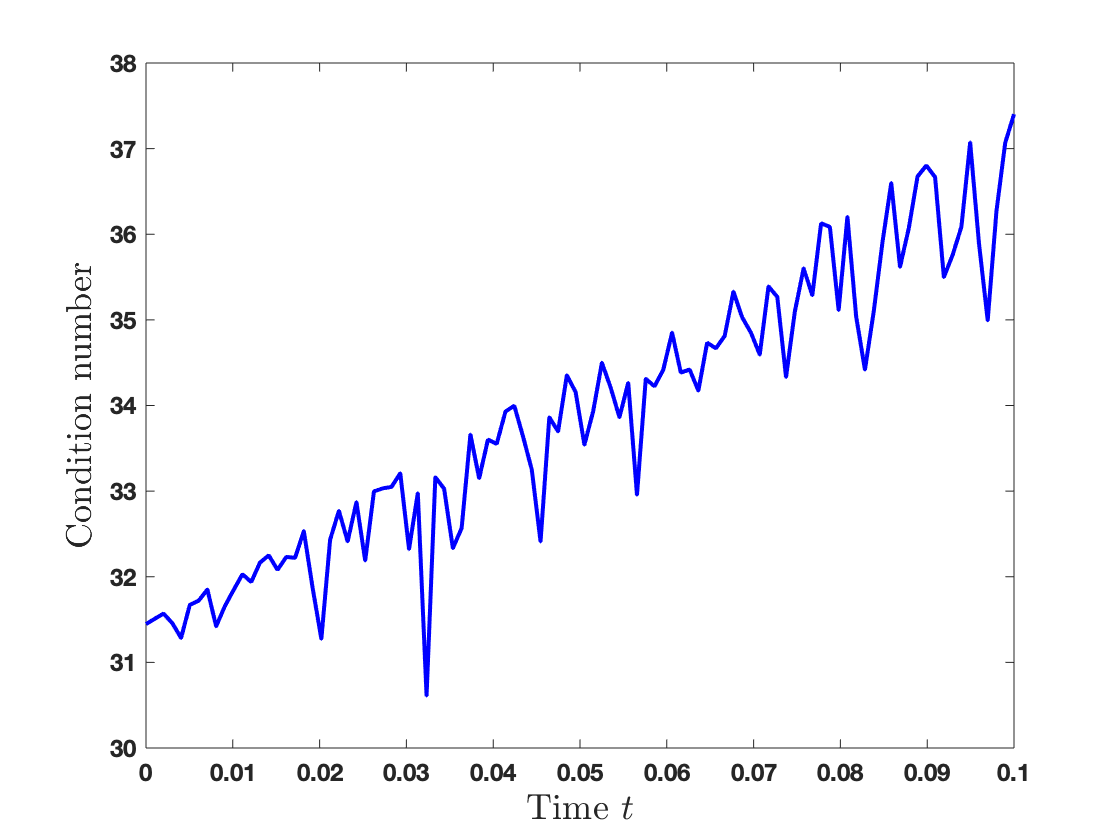}}
     \subfloat[$Nx=64, N_t=1000$]{\includegraphics[width=0.25\linewidth]{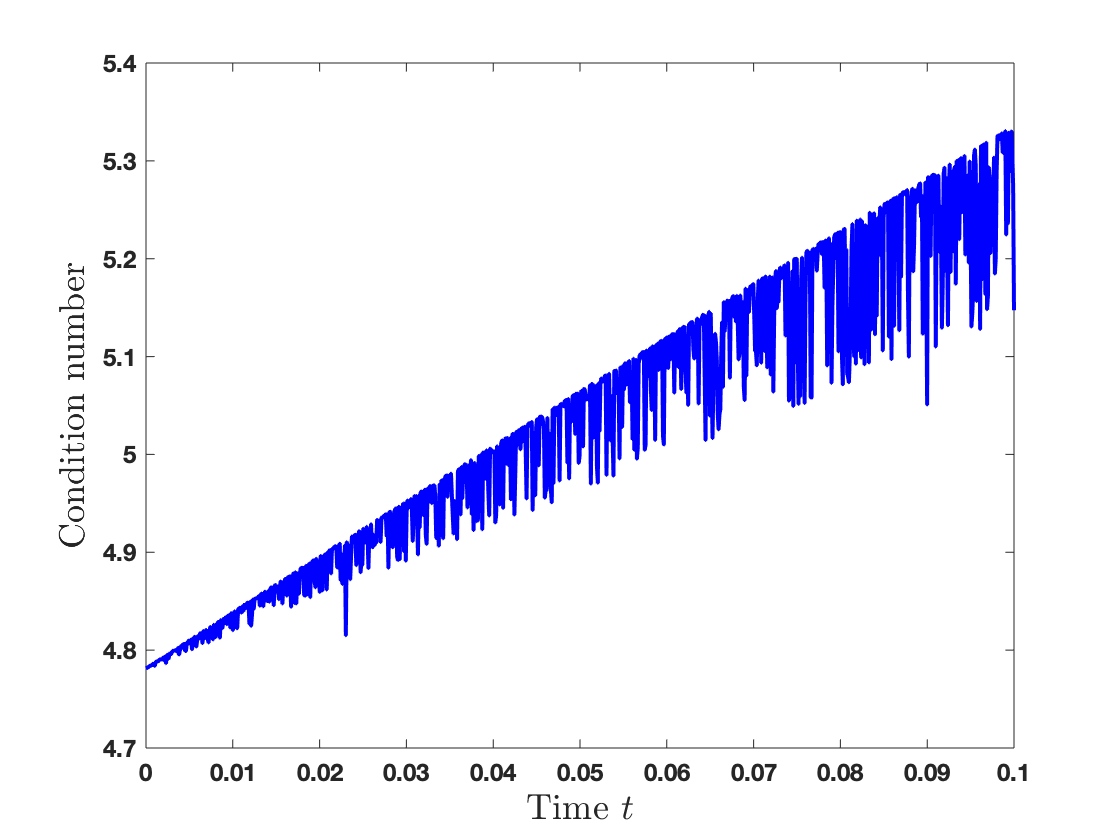}}
     \hspace{0.1in}
      \subfloat[$Nx=200, N_t=10$]{\includegraphics[width=0.25\linewidth]{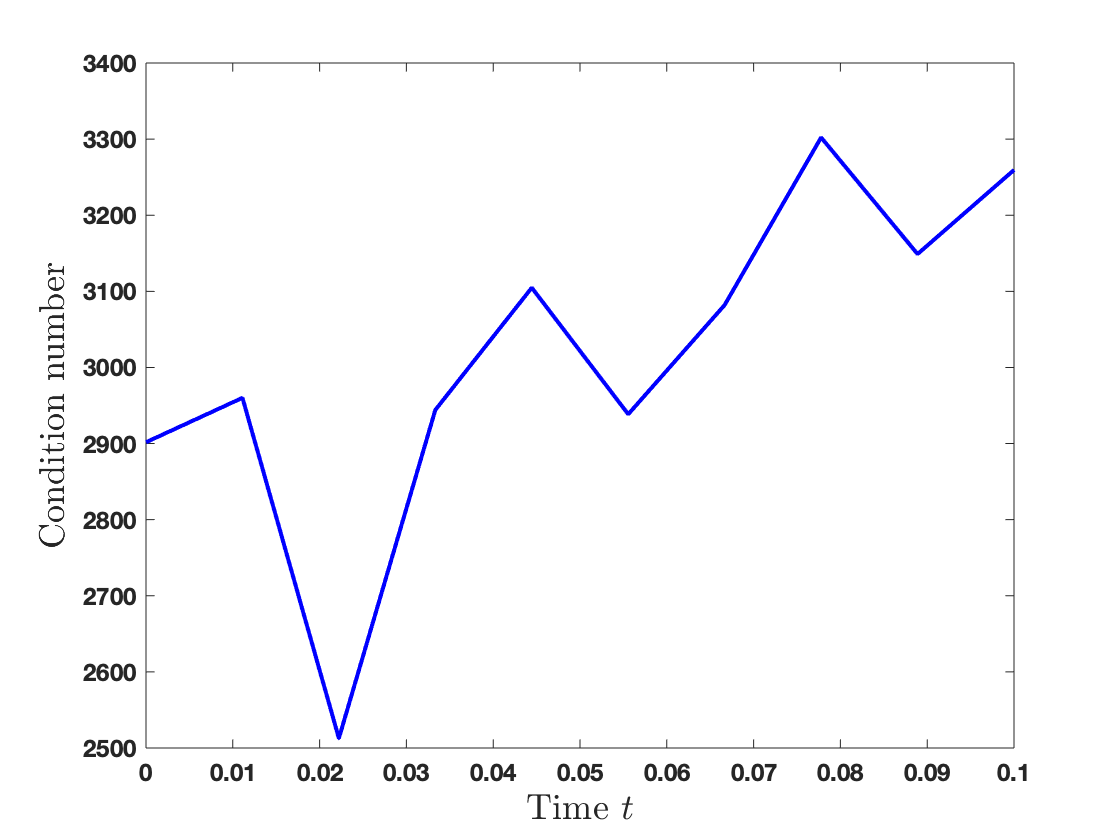}}
    \subfloat[$Nx=200, N_t=100$]{\includegraphics[width=0.25\linewidth]{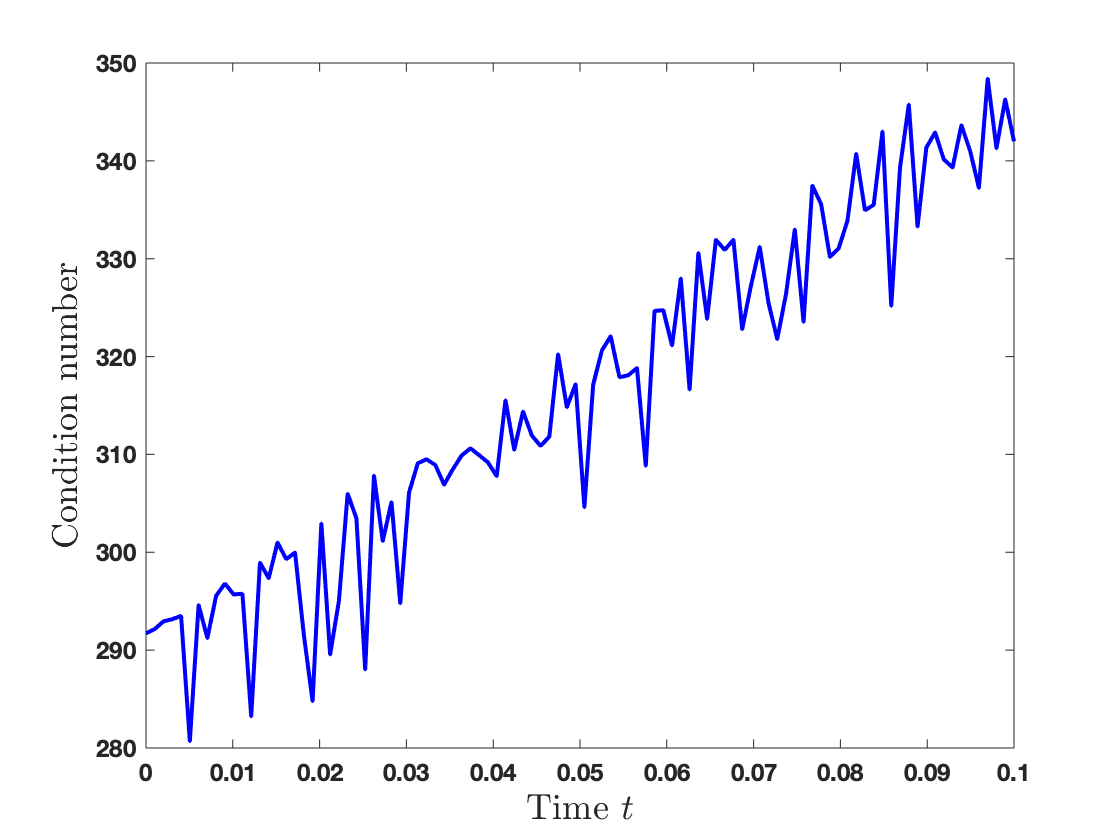}}
     \subfloat[$Nx=200, N_t=1000$]{\includegraphics[width=0.25\linewidth]{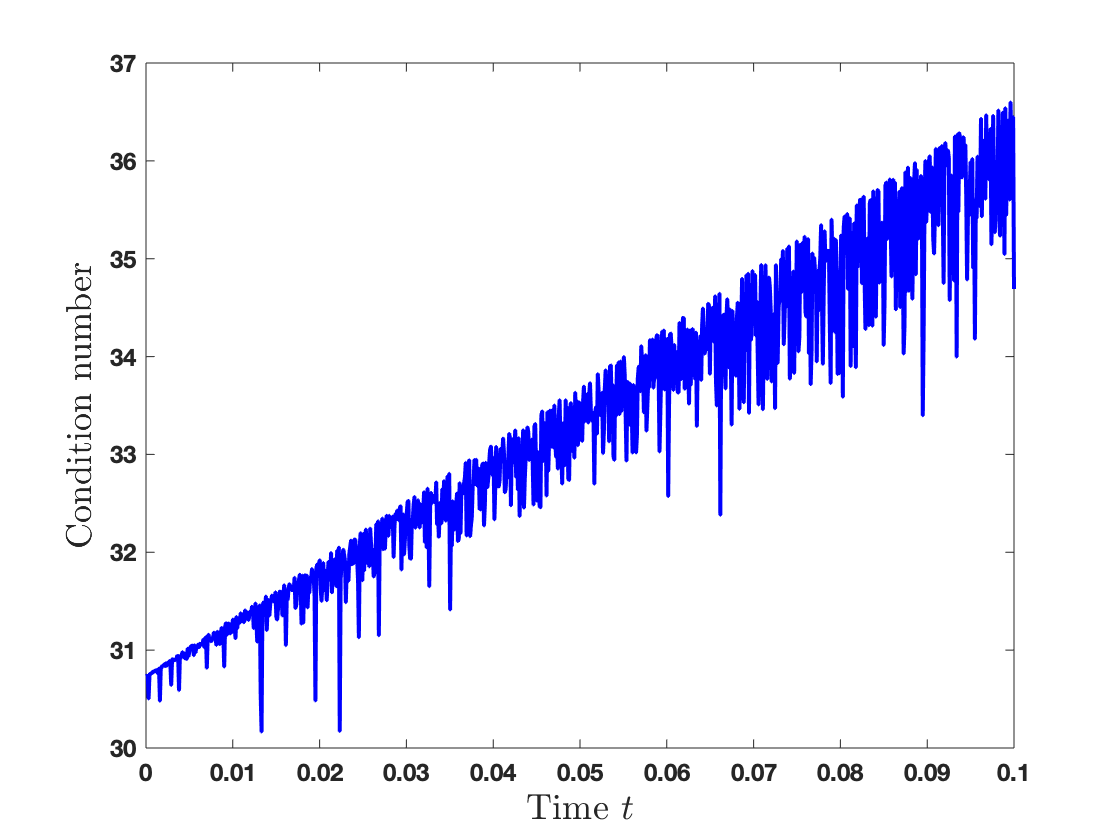}}
 \caption{The condition number of the matrix $A$ up to the final time $T=0.1$ for 1D case. BDF1 without projection.}
    \label{fig:cond-1}
\end{figure}

The spectral radius of the matrix $A$ is presented in \Cref{fig:spectral_rad-1}, which presents the time‑dependent spectral radius of matrix A for the 1D problem solved by the projection‑free BDF1 scheme up to \(T=0.1\). Results are shown for two spatial resolutions \(N_x=64, 200\) and three time‑step numbers \(N_t=10,100,1000\). In all cases, the spectral radius decreases monotonically over time. Higher spatial resolution significantly increases the spectral‑radius magnitude, while refining time steps moderately reduces its value. Regardless of discretization parameters, the decaying profile remains consistent, illustrating the joint influence of spatial and temporal discretization on the discrete operator’s spectral behaviour. The condition number of the matrix $A$ is presented in \Cref{fig:cond-1}. Simulations adopt two spatial resolutions \(N_x=64,200\) and three temporal discretizations \(N_t=10,100,1000\). The condition number generally grows over time for all configurations. Coarse time grids produce large‑amplitude oscillations, whereas finer time steps suppress violent fluctuations and yield smoother increasing trends. Higher spatial resolution substantially amplifies the condition number, showing spatial discretization strongly affects the ill‑posedness of the discrete operator.

\begin{figure}[htbp]
    \centering
    \subfloat[$Nx=16, N_t=10$]{\includegraphics[width=0.3\linewidth]{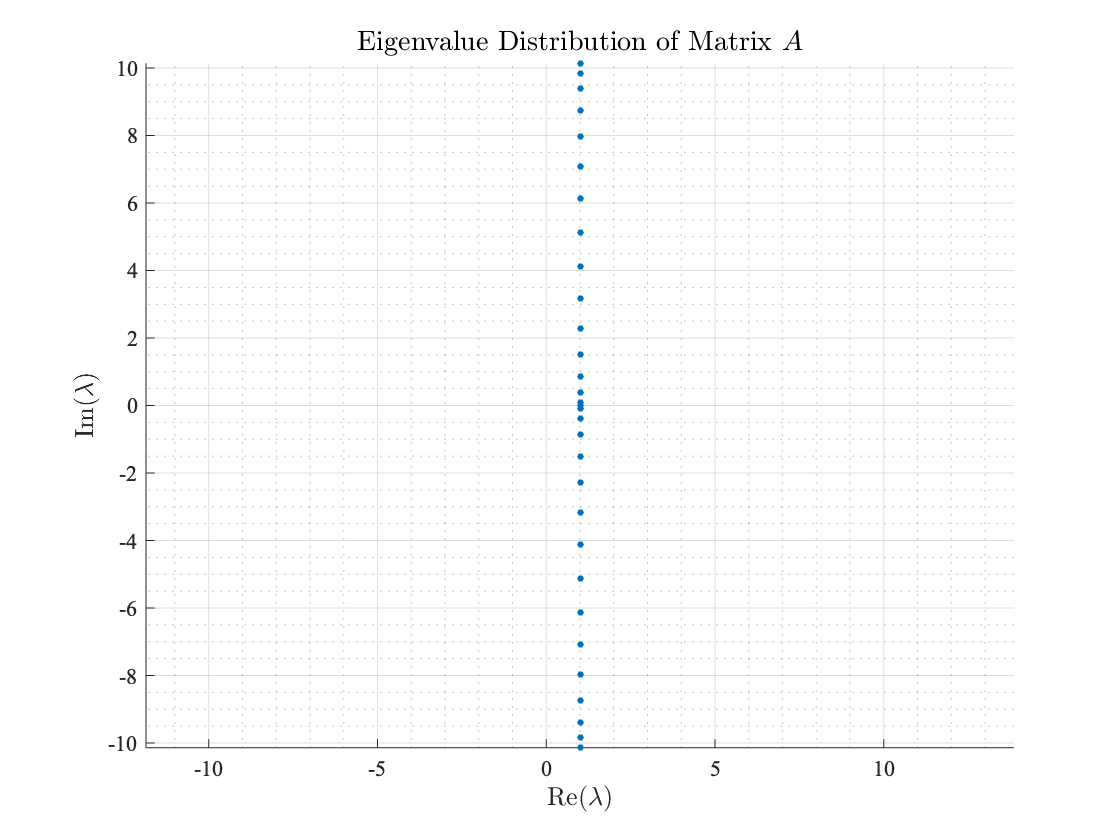}}
    \subfloat[$Nx=16, N_t=100$]{\includegraphics[width=0.3\linewidth]{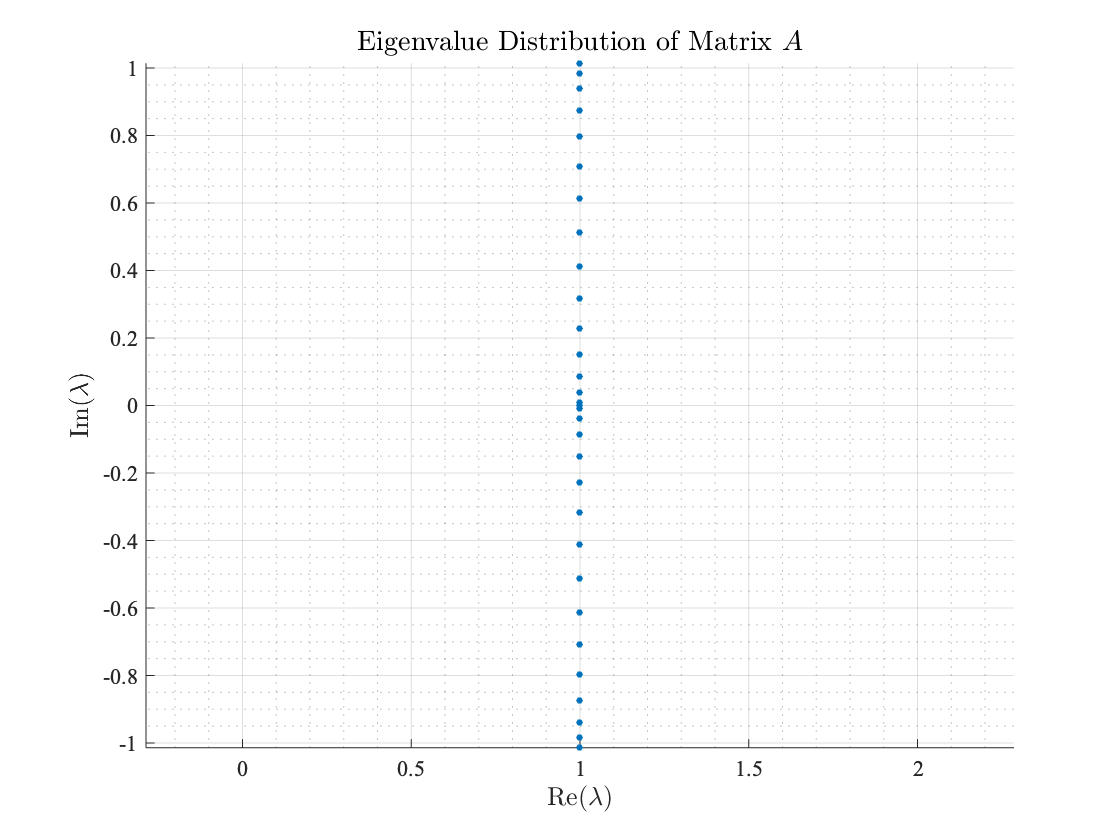}}
     \subfloat[$Nx=16, N_t=1000$]{\includegraphics[width=0.3\linewidth]{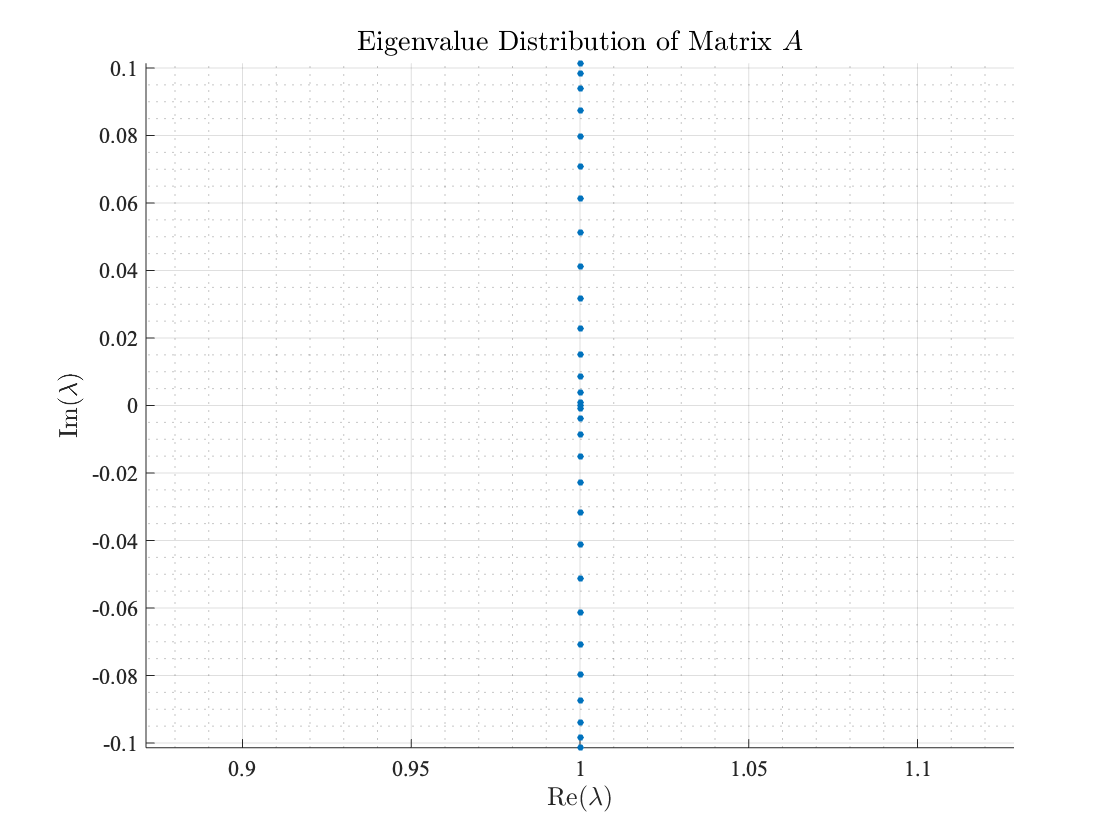}}
     \hspace{0.1in}
      \subfloat[$Nx=100, N_t=10$]{\includegraphics[width=0.3\linewidth]{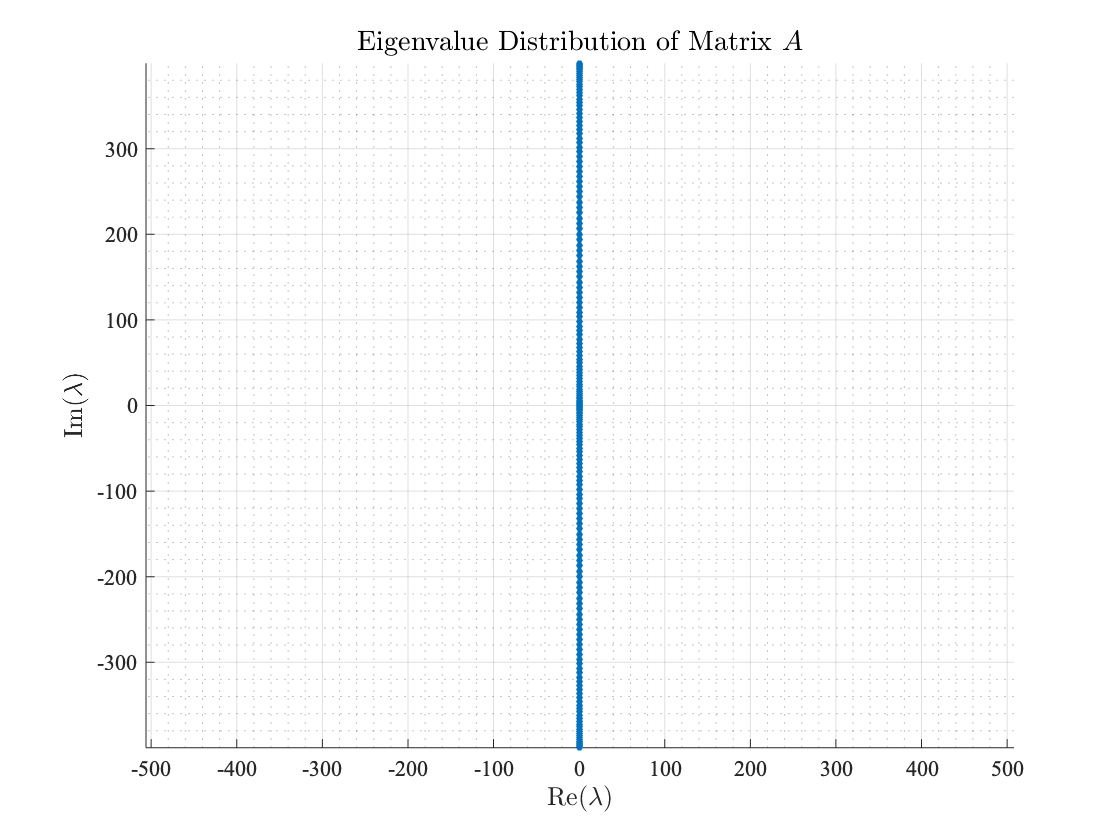}}
    \subfloat[$Nx=100, N_t=100$]{\includegraphics[width=0.3\linewidth]{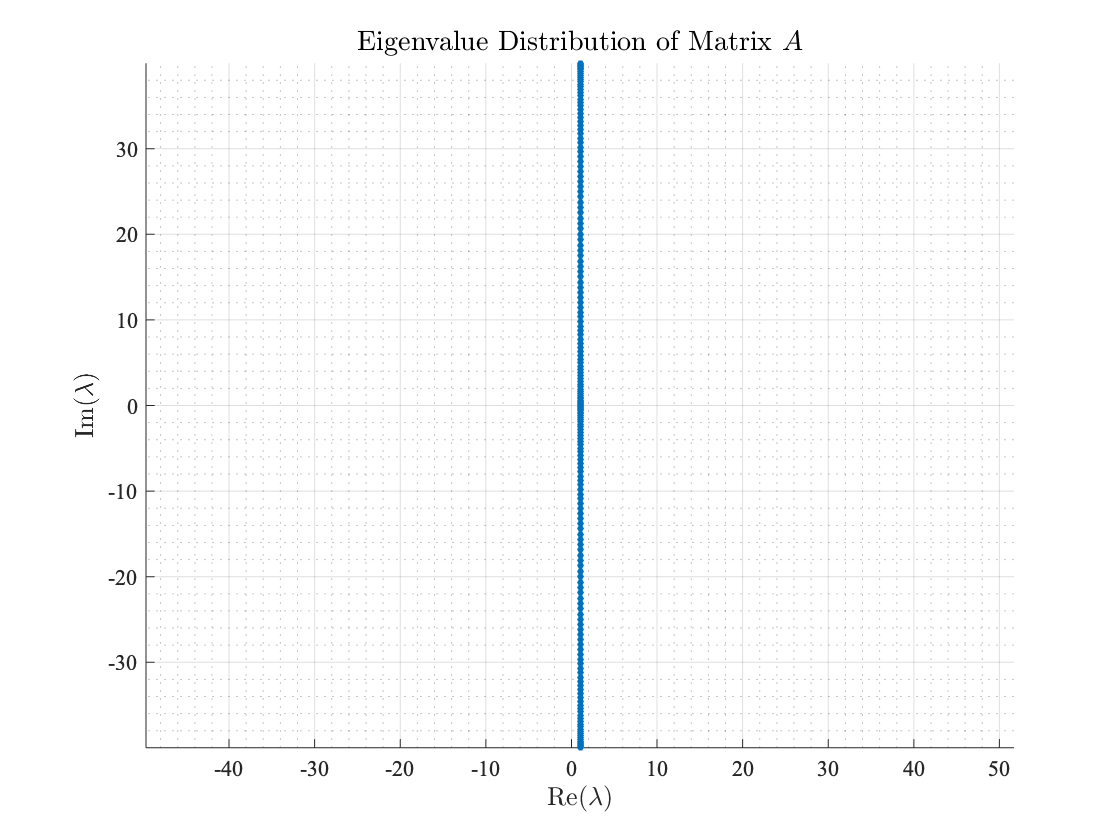}}
     \subfloat[$Nx=100, N_t=1000$]{\includegraphics[width=0.3\linewidth]{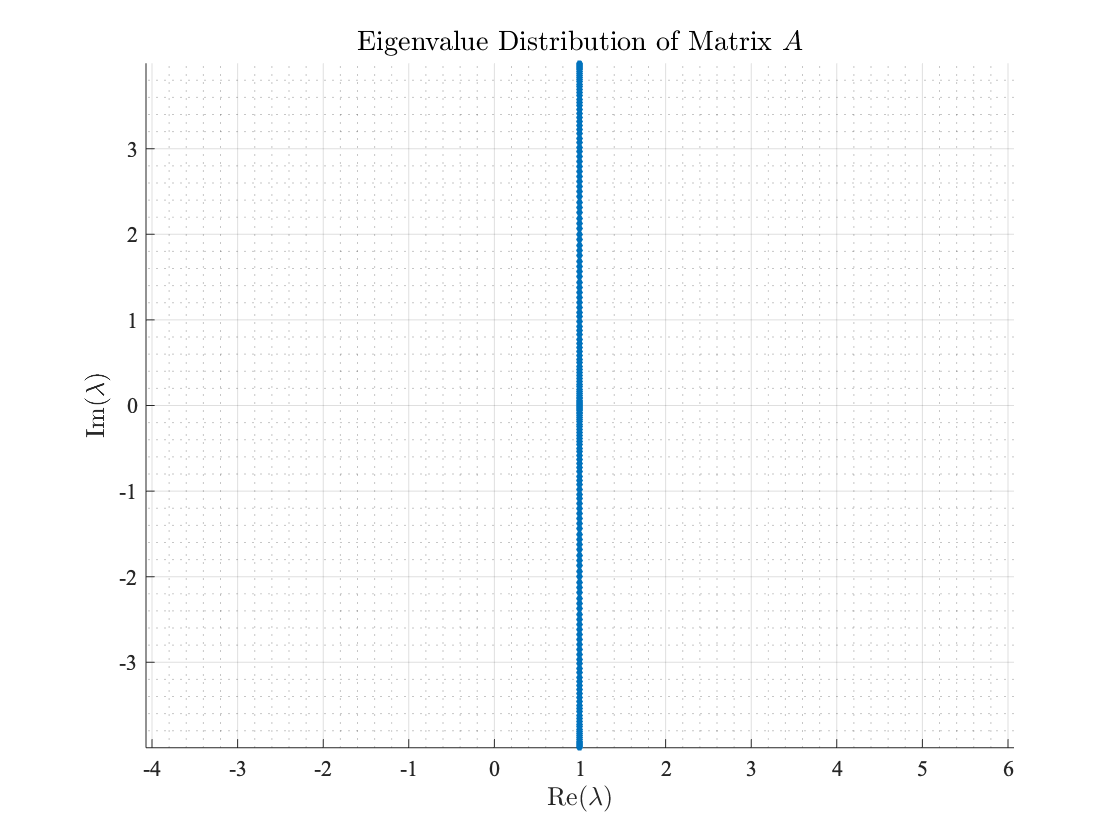}}
     \hspace{0.1in}
      \subfloat[$Nx=200, N_t=10$]{\includegraphics[width=0.3\linewidth]{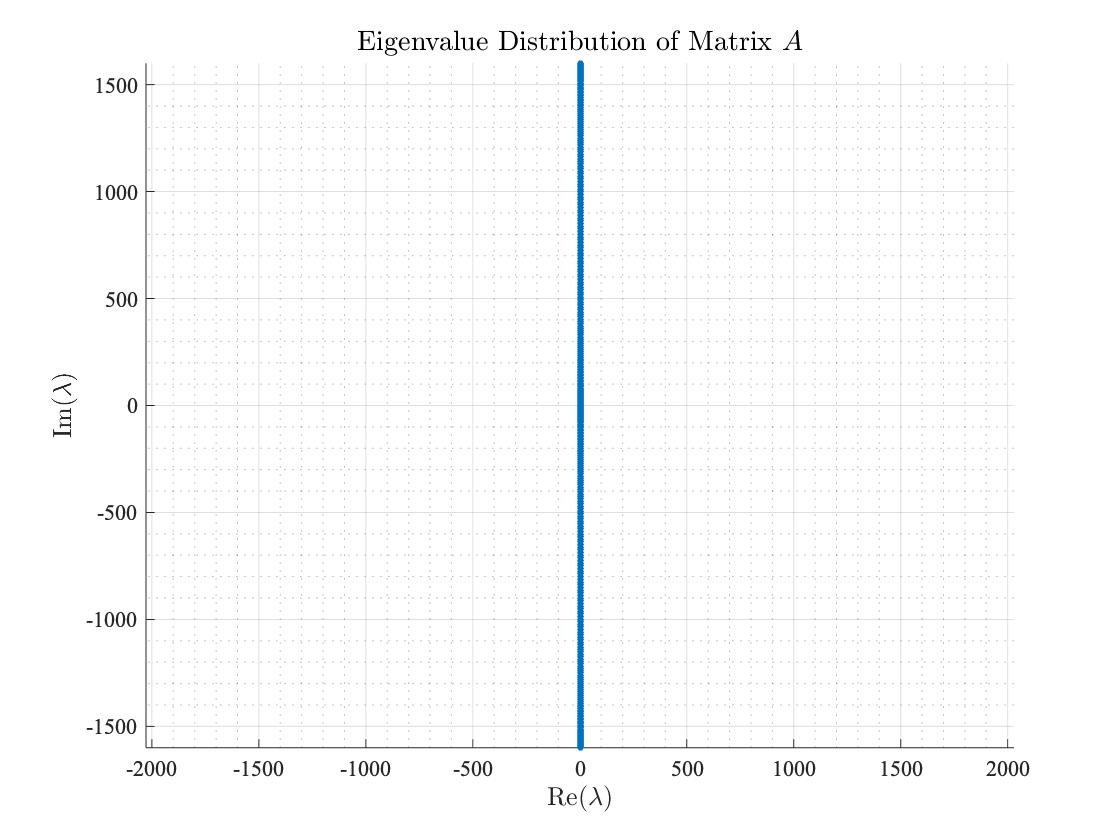}}
    \subfloat[$Nx=200, N_t=100$]{\includegraphics[width=0.3\linewidth]{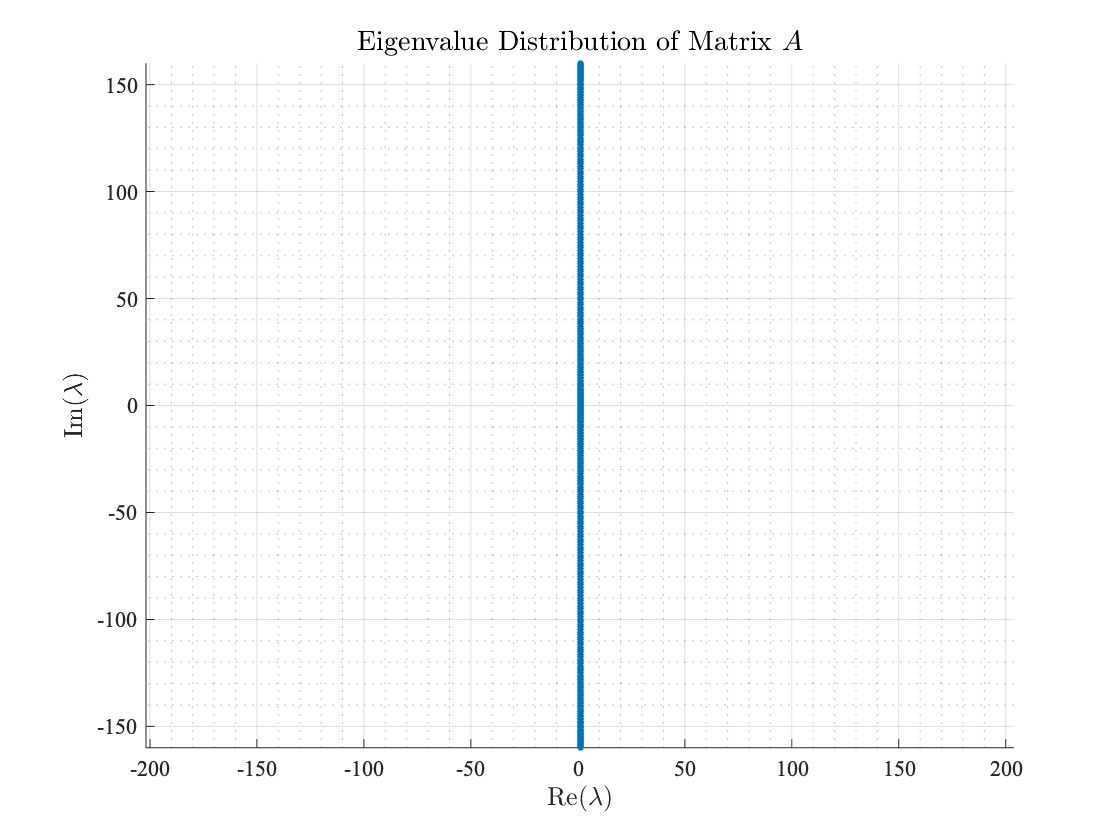}}
     \subfloat[$Nx=200, N_t=1000$]{\includegraphics[width=0.3\linewidth]{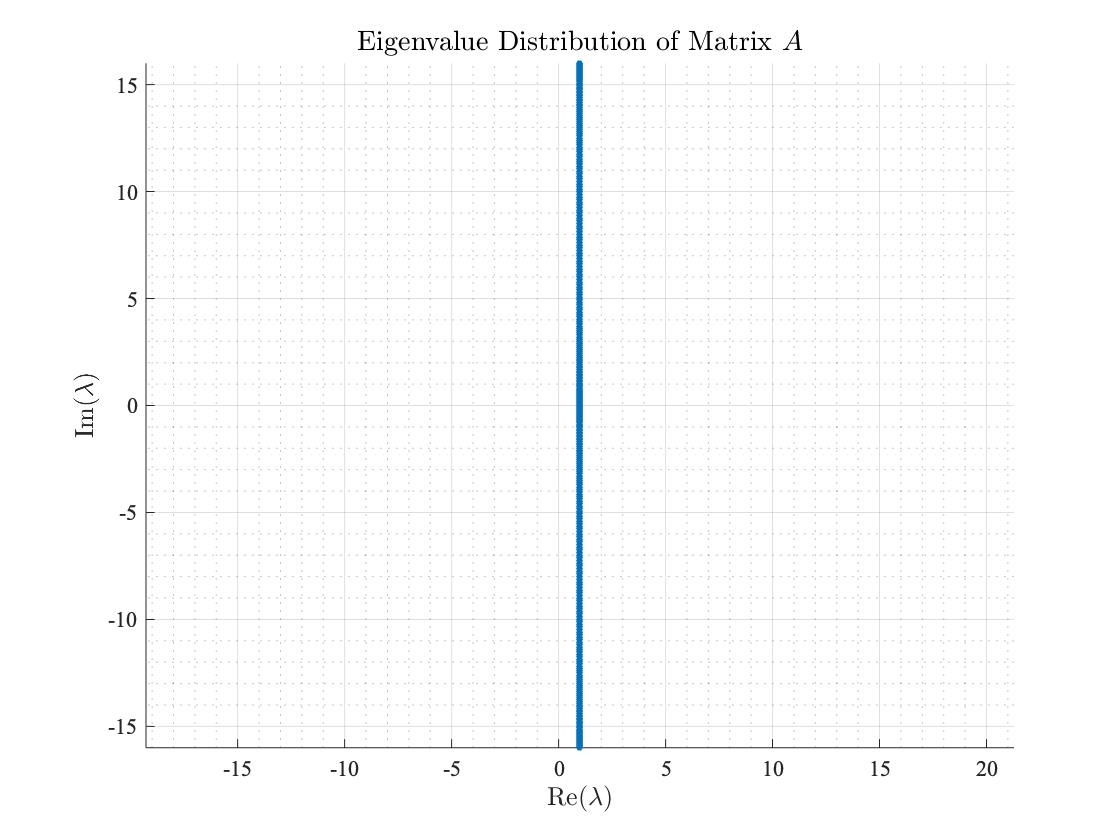}}
 \caption{The distribution of eigenvalues of the matrix $A$ at the final time $T=0.1$ for 1D case. BDF1 without projection.}
    \label{fig:dist-1}
\end{figure}

The distribution of eigenvalues of the matrix $A$ is presented in \Cref{fig:dist-1}. Numerical tests cover three spatial resolutions \(N_x=16,100,200\) and three time‑step numbers \(N_t=10,100,1000\). All eigenvalues lie along the real axis with zero imaginary components. Increasing spatial resolution expands the range of real eigenvalues significantly. By contrast, refining temporal discretization greatly compresses the eigenvalue spectrum, shrinking the spread of real‑valued eigenvalues for fixed spatial grids, which demonstrates the combined discretization effects on spectral properties. The matrix A arising from BDF1 discretization possesses eigenvalues spread over a large range along the negative real axis, leading to severe ill‑conditioning as spatial resolution \(N_x\) increases. Multigrid (AMG) methods are well‑suited here: local relaxation efficiently damps high‑frequency large‑magnitude error components on fine grids, while coarse‑grid corrections eliminate slowly‑converging low‑frequency modes associated with small‑magnitude eigenvalues near the origin. When used as a preconditioner for GMRES, AMG clusters the preconditioned‑matrix spectrum near unity, yields nearly grid‑size‑independent convergence, linear computational complexity \(O(N)\), and avoids prohibitive memory cost from direct LU factorizations for large \(N_x\).

\subsection{Traditional iterative solver and preconditioning}

\begin{figure}[htbp]
    \centering
    \subfloat[GMRES iter]{\includegraphics[width=0.25\linewidth]{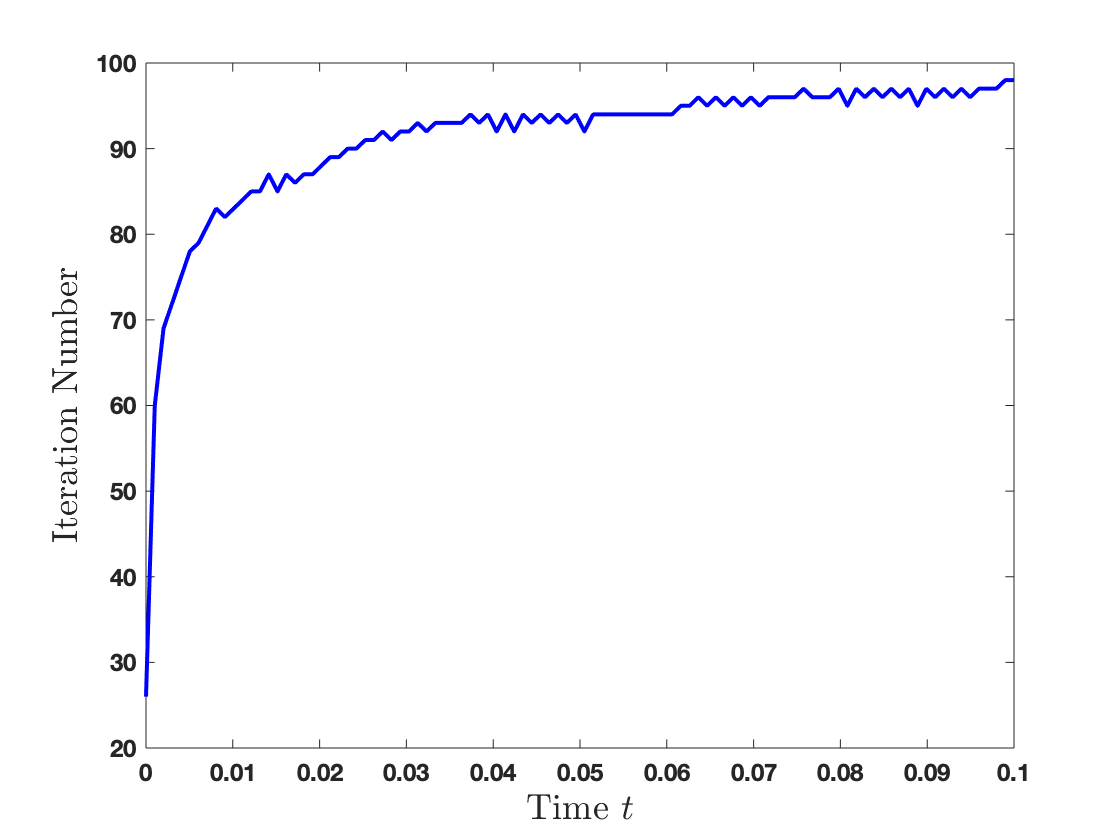}}
     \subfloat[BiCGstab iter]{\includegraphics[width=0.25\linewidth]{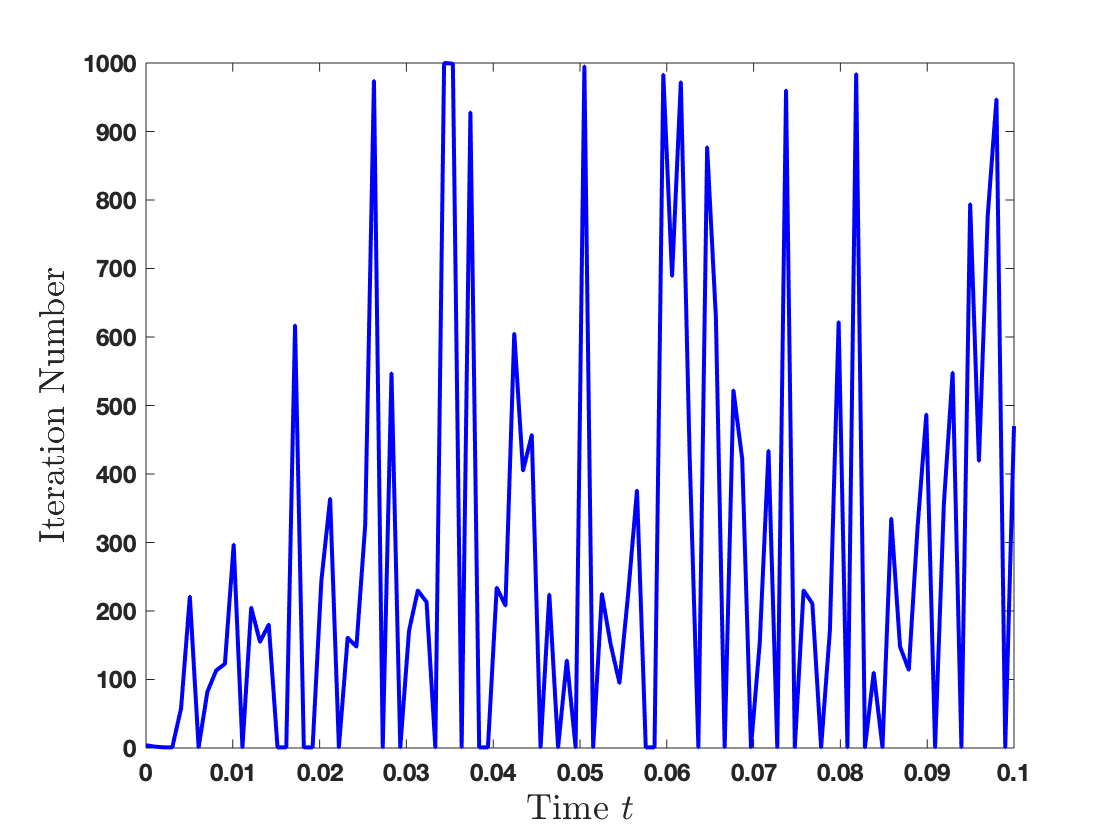}}
     \subfloat[TFQMR iter]{\includegraphics[width=0.25\linewidth]{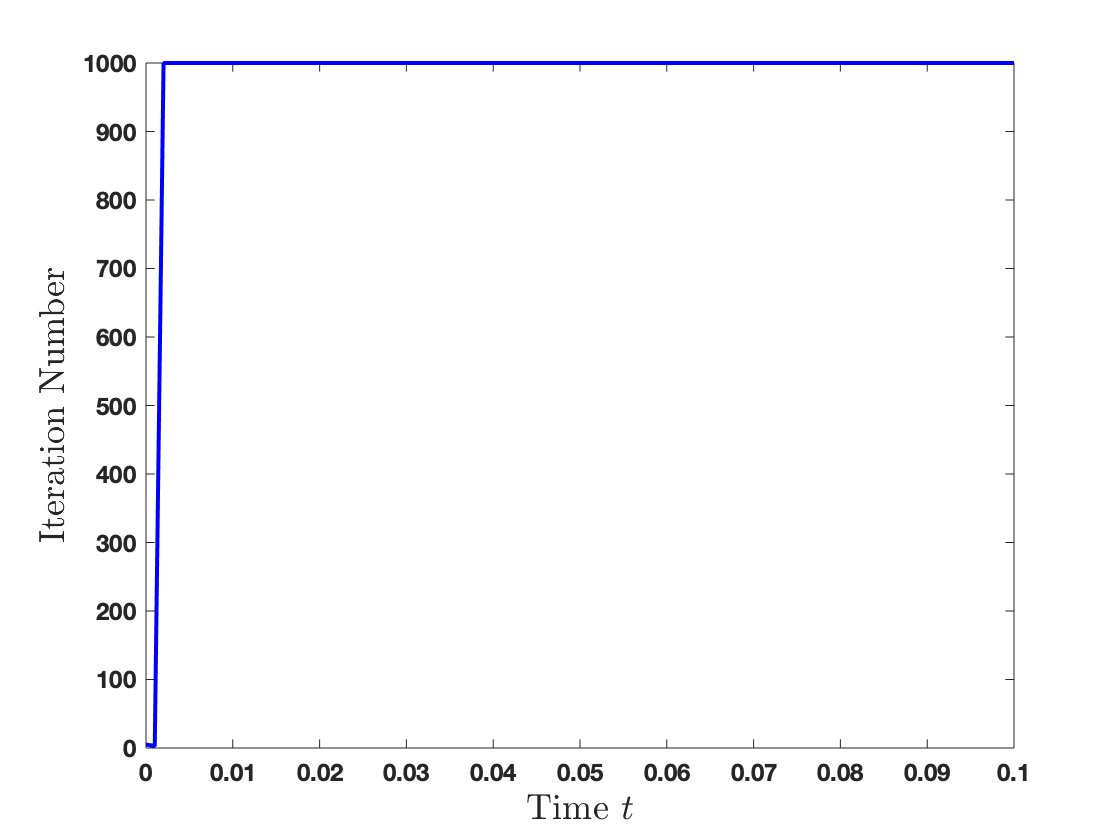}}
      \hspace{0.1in}
      \subfloat[GMRES relres]{\includegraphics[width=0.25\linewidth]{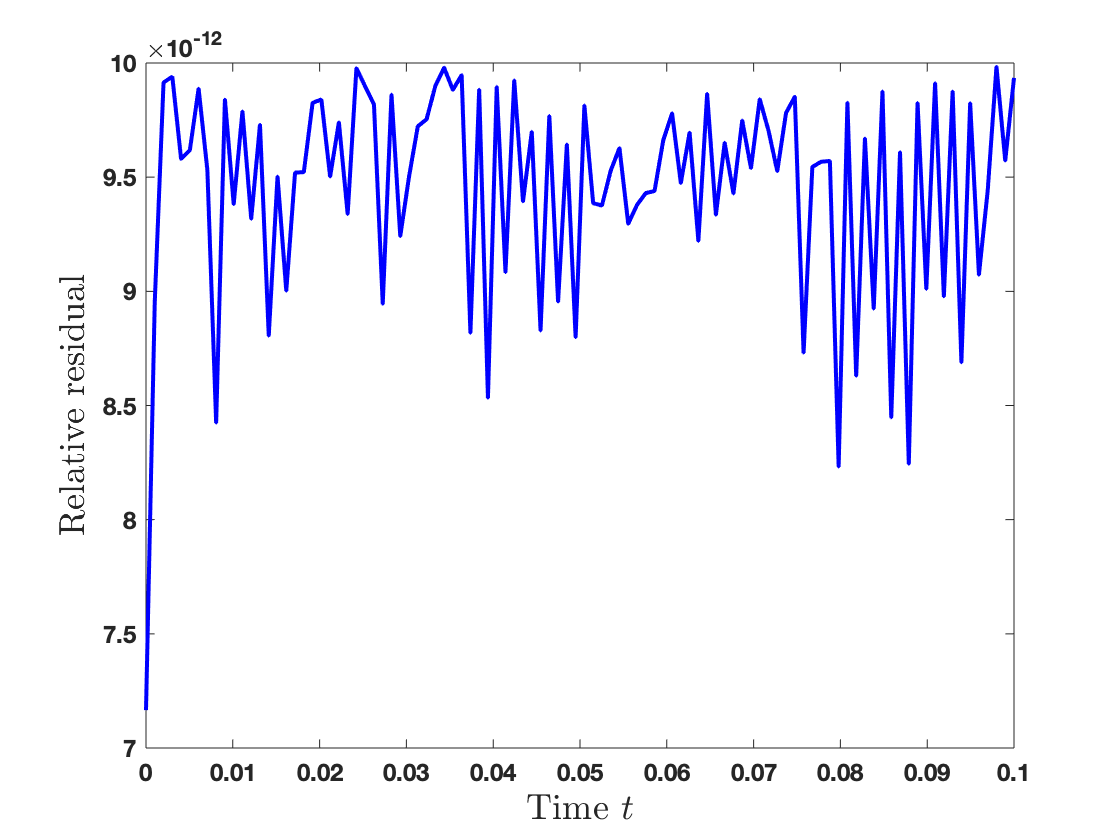}}
     \subfloat[BiCGstab relres]{\includegraphics[width=0.25\linewidth]{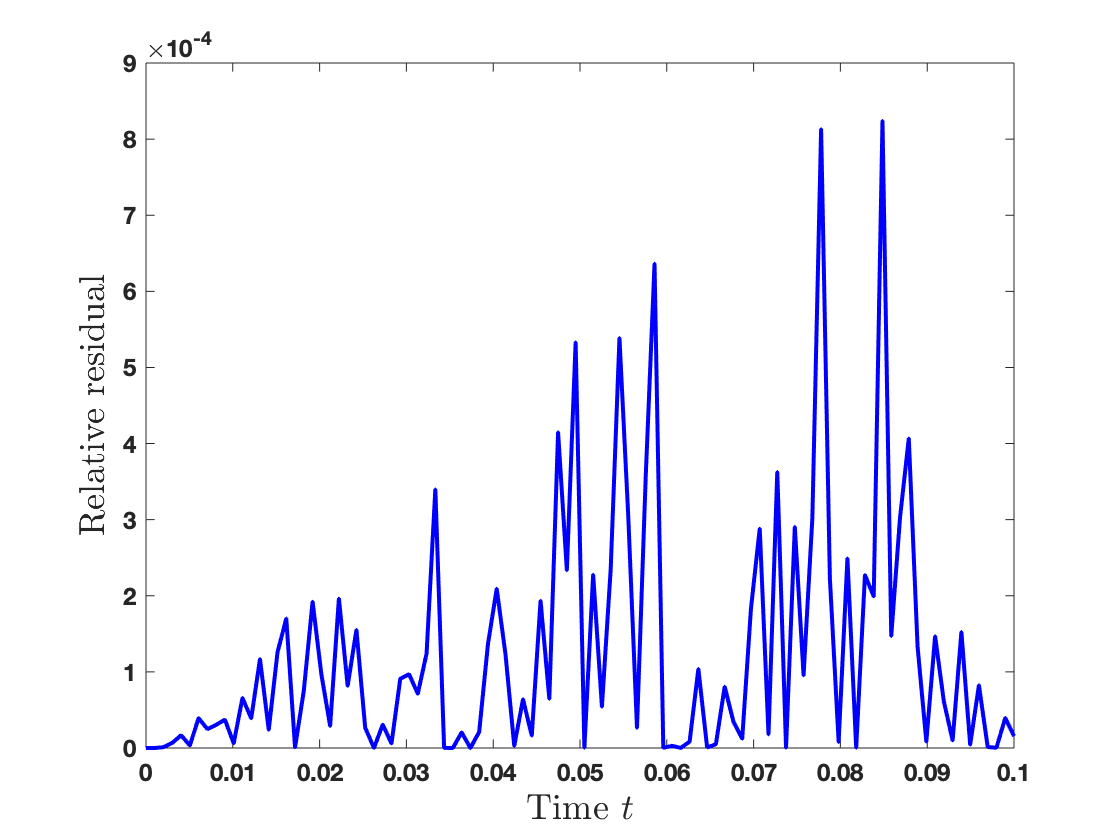}}
     \subfloat[TFQMR relres]{\includegraphics[width=0.25\linewidth]{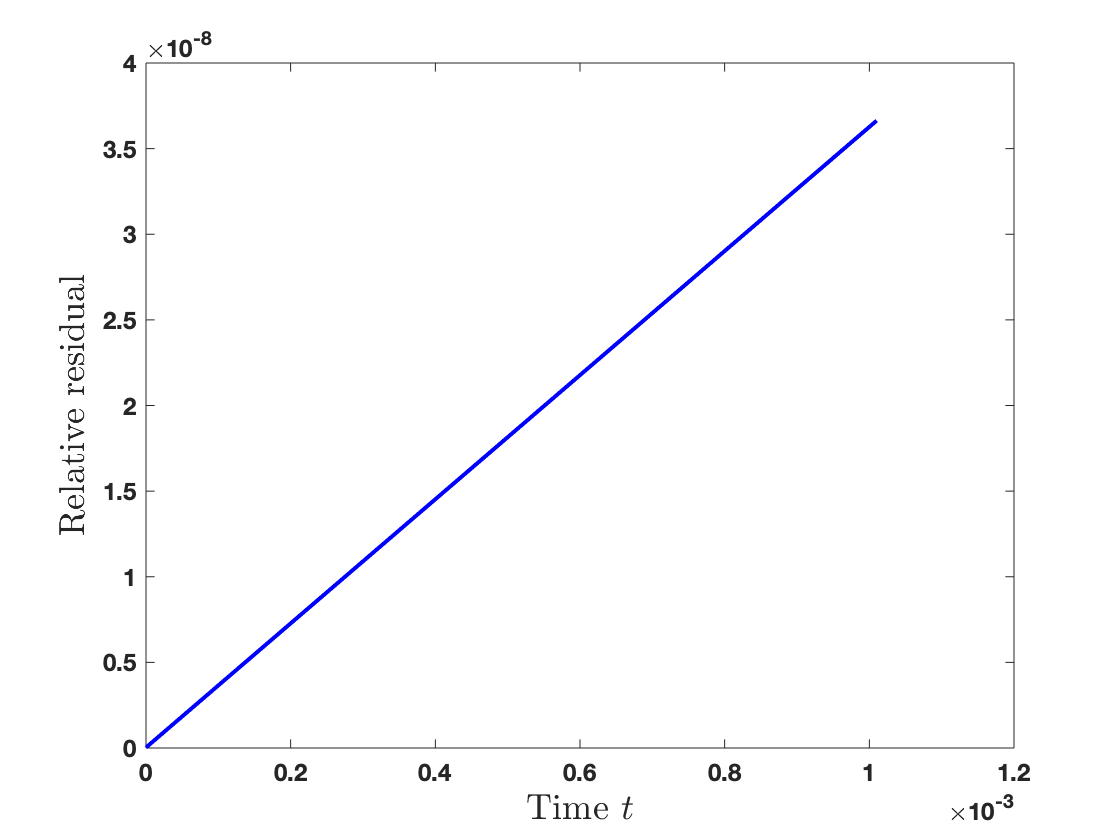}}
      \hspace{0.1in}
    \subfloat[GMRES LU iter]{\includegraphics[width=0.25\linewidth]{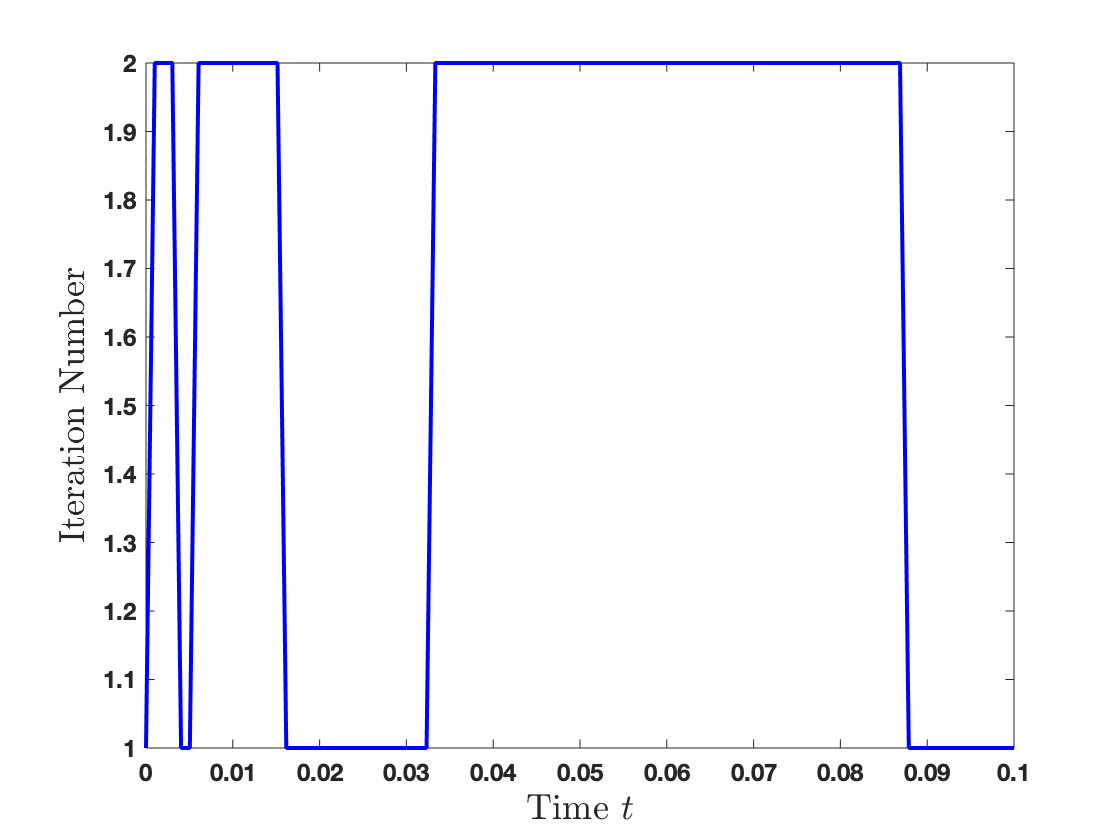}}
    % \hspace{0.1in}
    \subfloat[BiCGstab LU iter]{\includegraphics[width=0.25\linewidth]{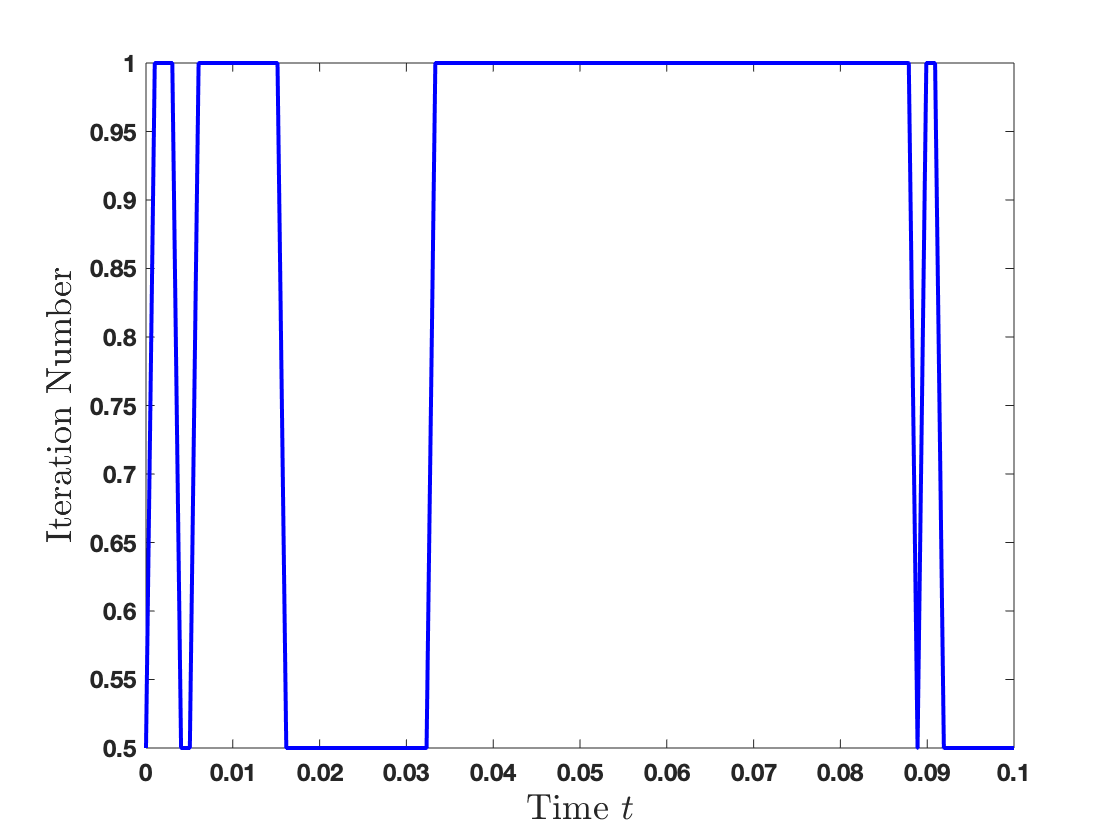}}
    % \hspace{0.1in}
    \subfloat[TFQMR LU iter]{\includegraphics[width=0.25\linewidth]{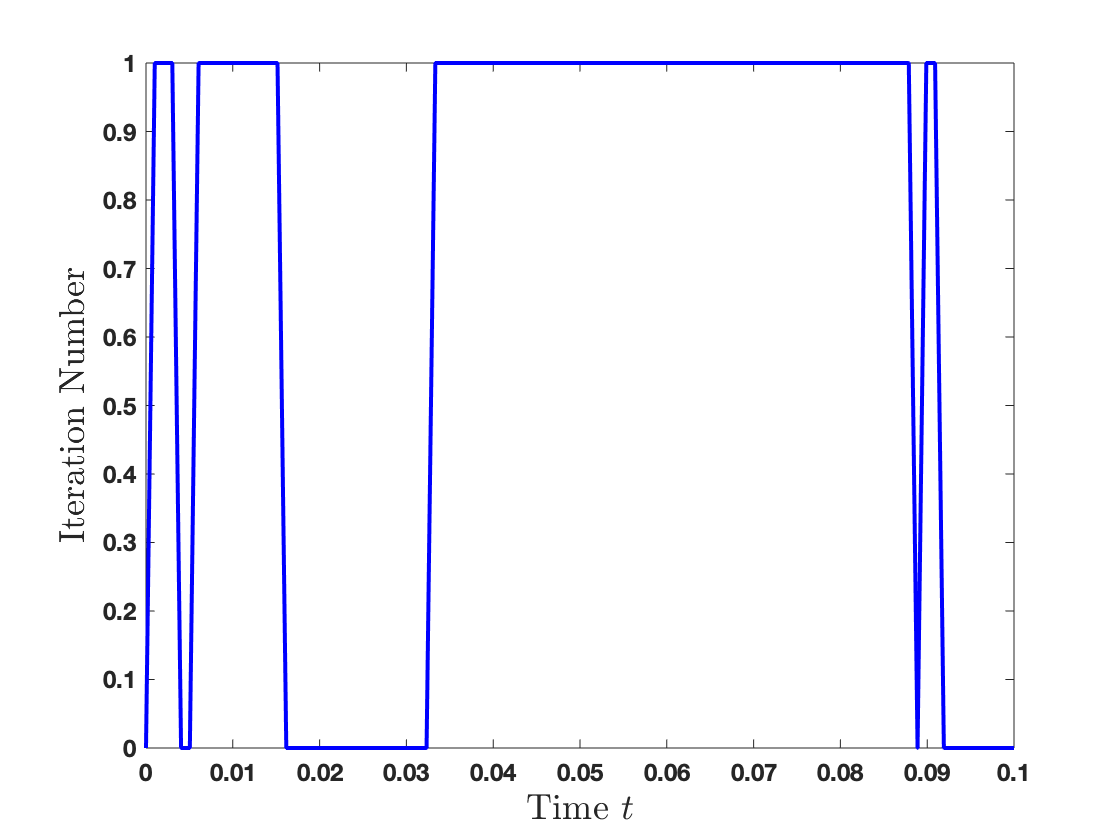}}
    \hspace{0.1in}
    \subfloat[GMRES LU relres]{\includegraphics[width=0.25\linewidth]{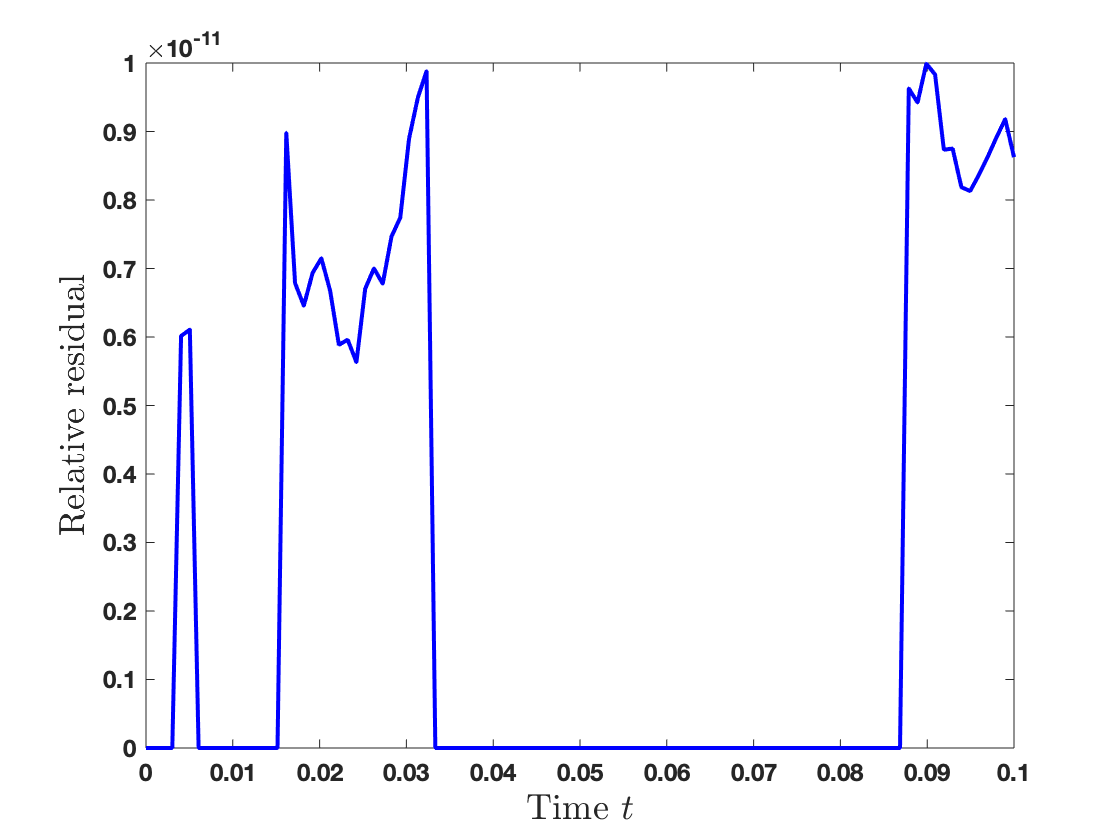}}
    % \hspace{0.1in}
    \subfloat[BiCGstab LU relres]{\includegraphics[width=0.25\linewidth]{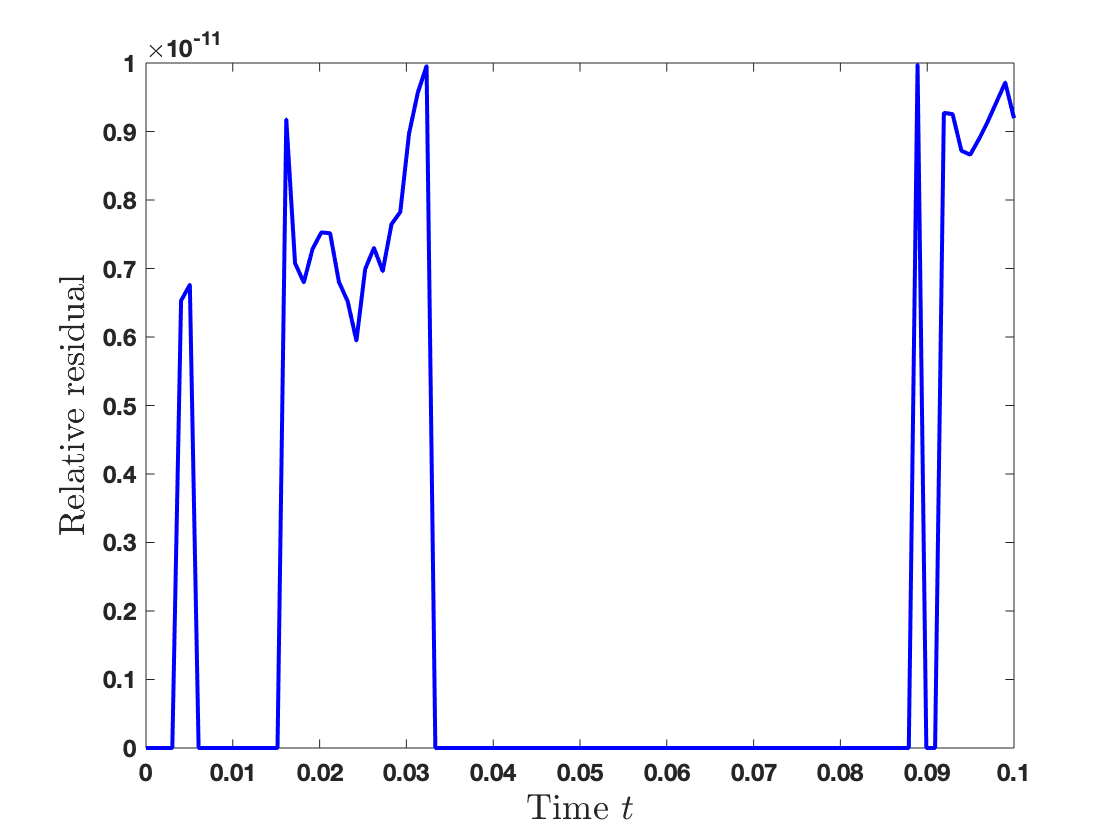}}
    % \hspace{0.1in}
    \subfloat[TFQMR LU relres]{\includegraphics[width=0.25\linewidth]{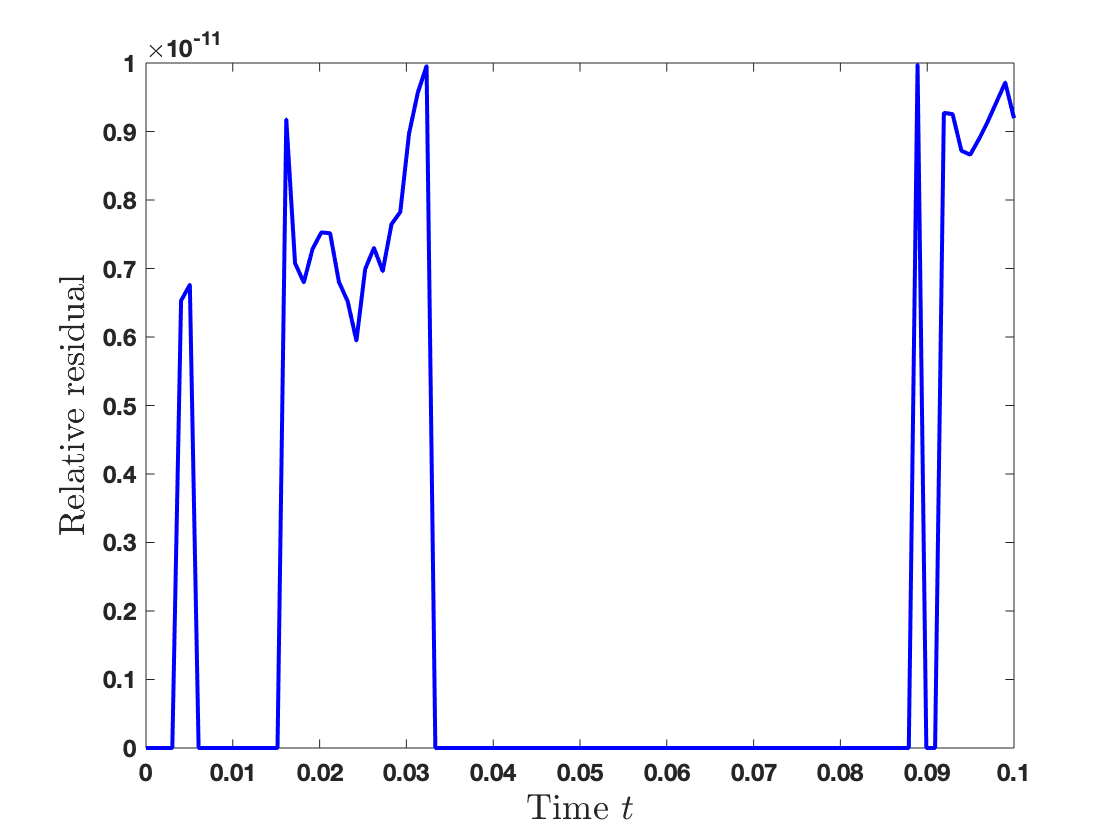}}
    \caption{Different linear solvers with $\alpha=0$, $N_x=64$, $N_t=100$ up to the final time $T=0.1$ for 1D case. BDF1 without projection.}
    \label{fig:BDF1-1}
\end{figure}

\begin{figure}[htbp]
    \centering
    \subfloat[GMRES iter]{\includegraphics[width=0.25\linewidth]{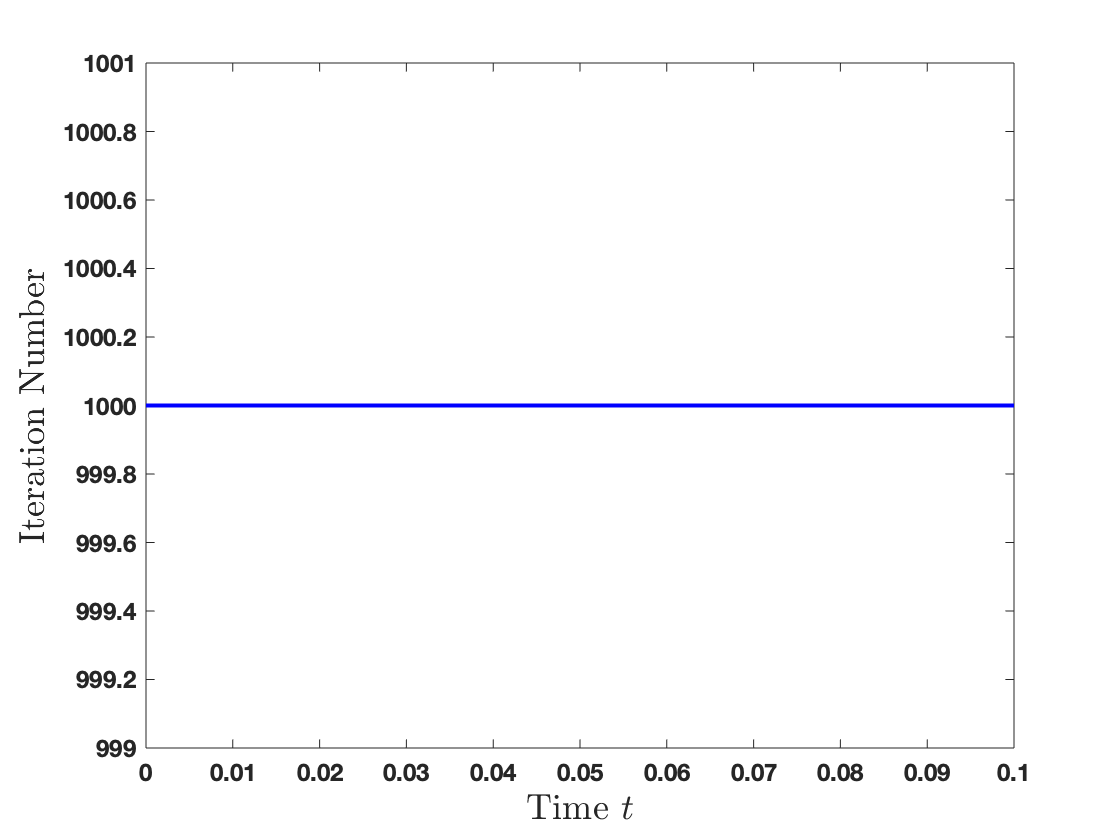}}
      \subfloat[BiCGstab iter]{\includegraphics[width=0.25\linewidth]{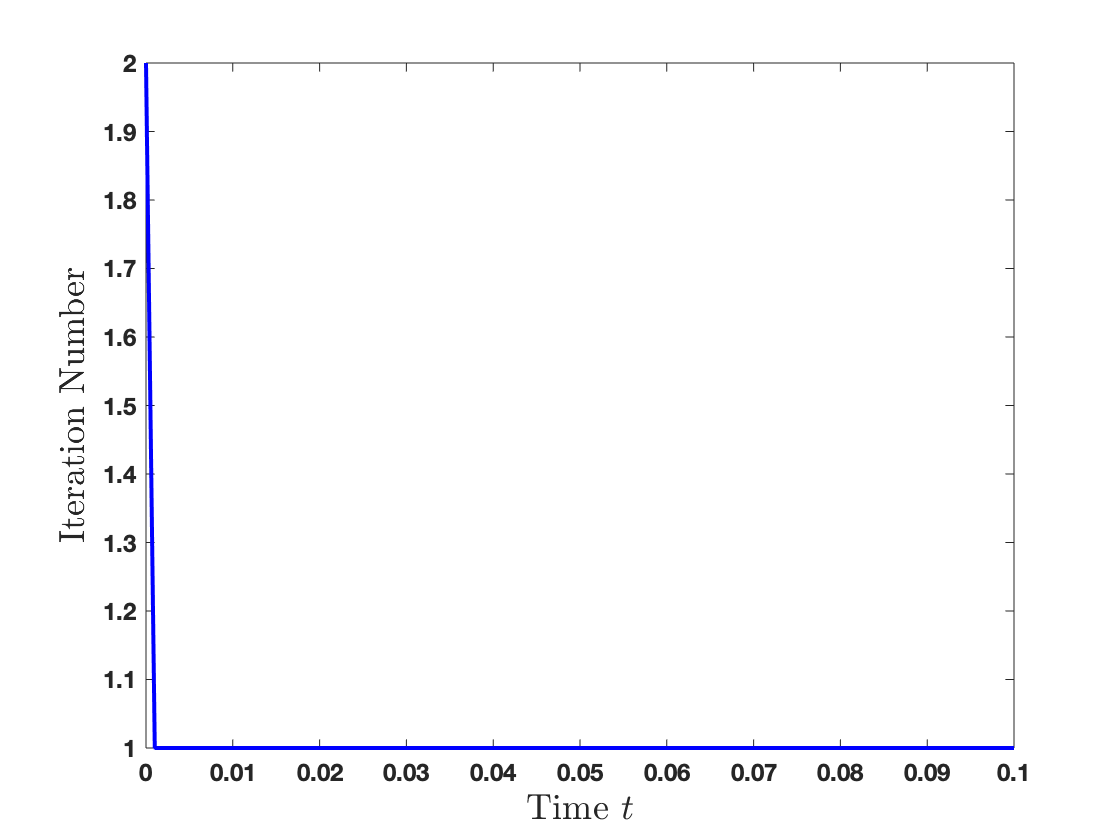}}
      \subfloat[TFQMR iter]{\includegraphics[width=0.25\linewidth]{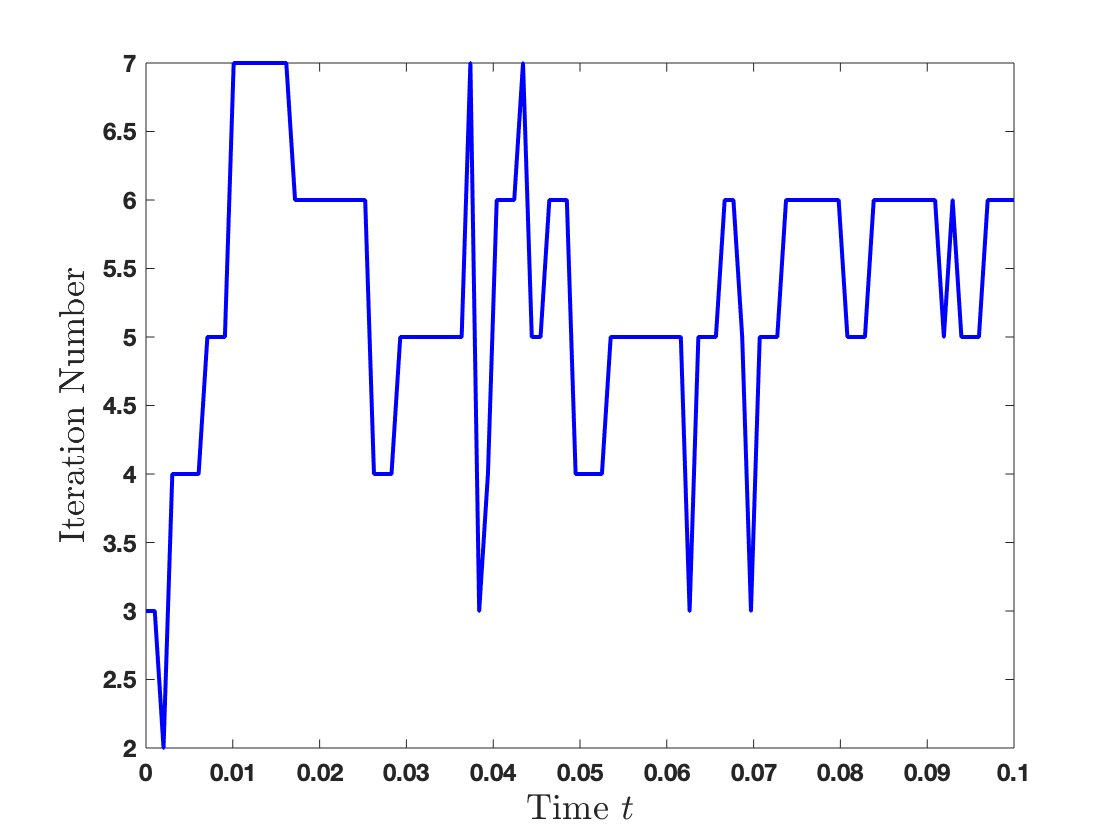}}
      \hspace{0.1in}
      \subfloat[GMRES relres]{\includegraphics[width=0.25\linewidth]{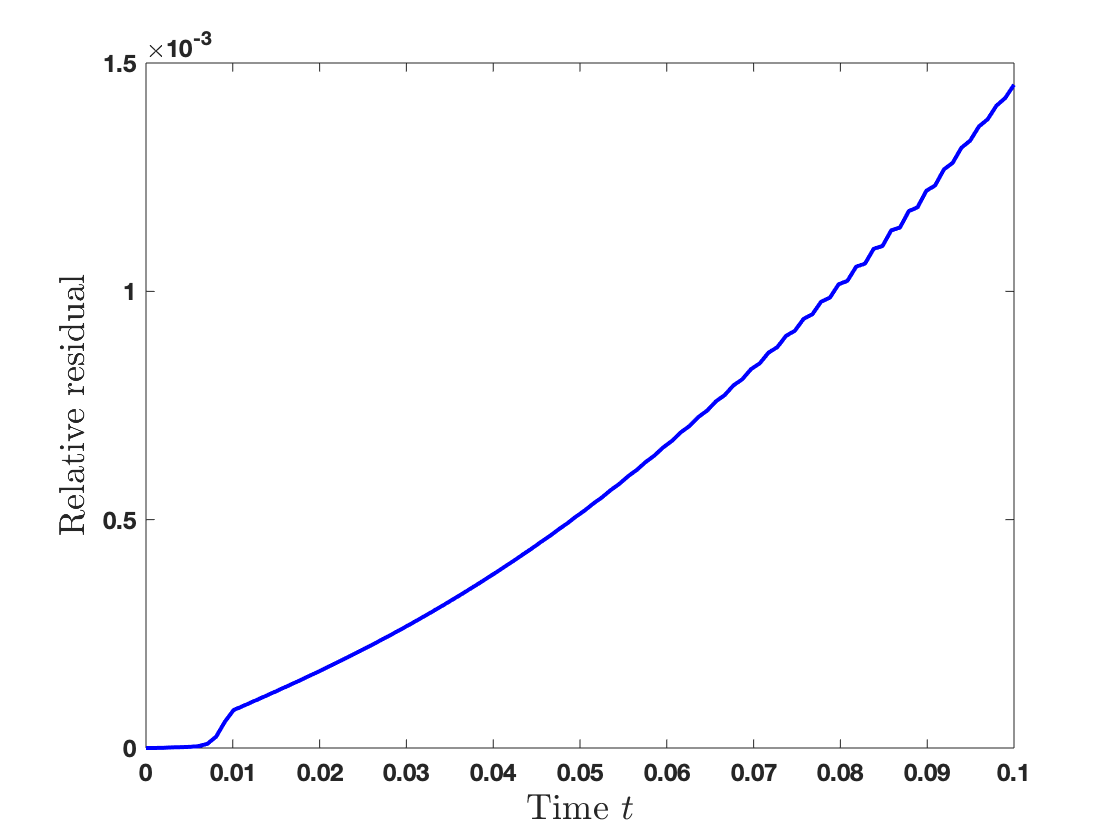}}
      \subfloat[BiCGstab relres]{\includegraphics[width=0.25\linewidth]{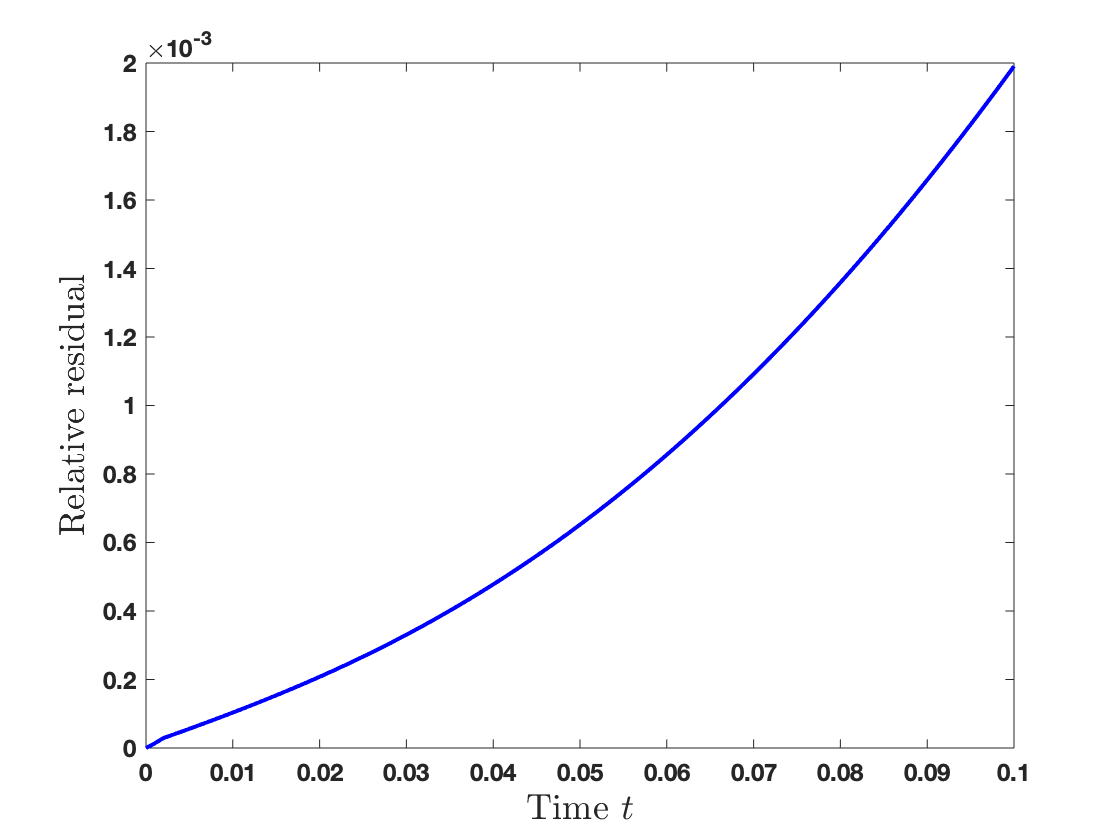}}
      \subfloat[TFQMR relres]{\includegraphics[width=0.25\linewidth]{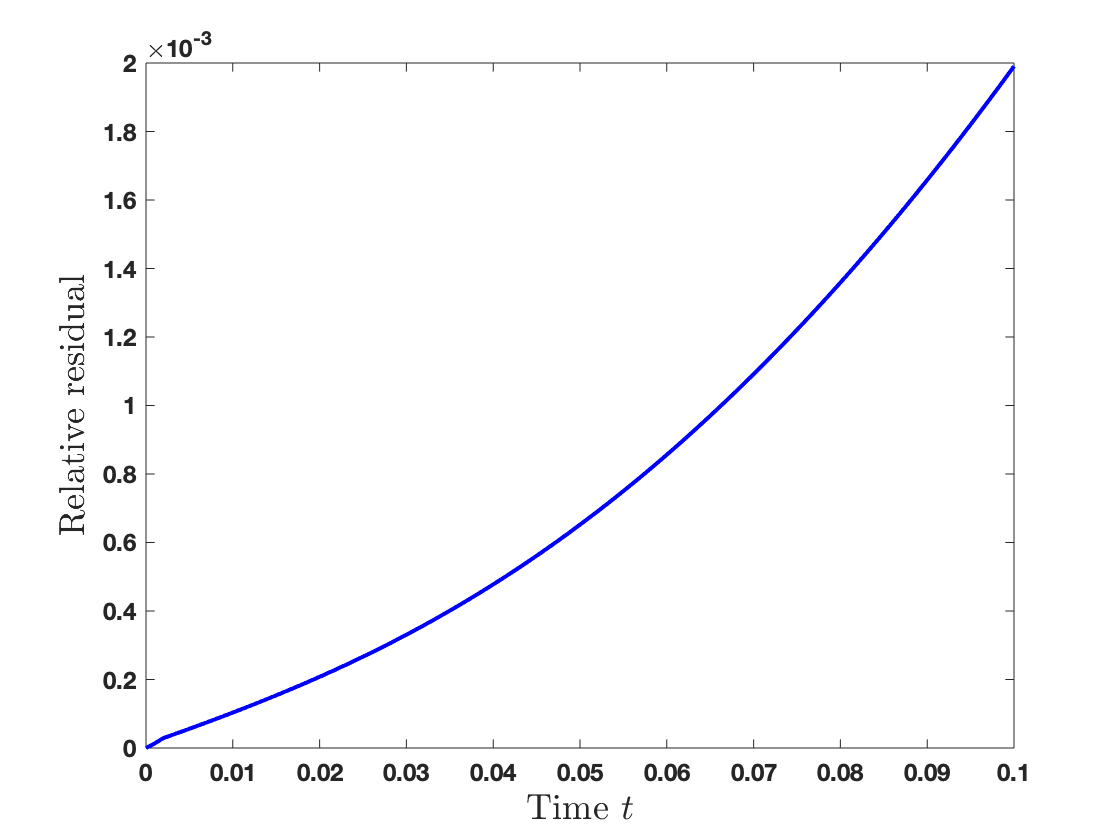}}
      \hspace{0.1in}
    \subfloat[GMRES LU iter]{\includegraphics[width=0.25\linewidth]{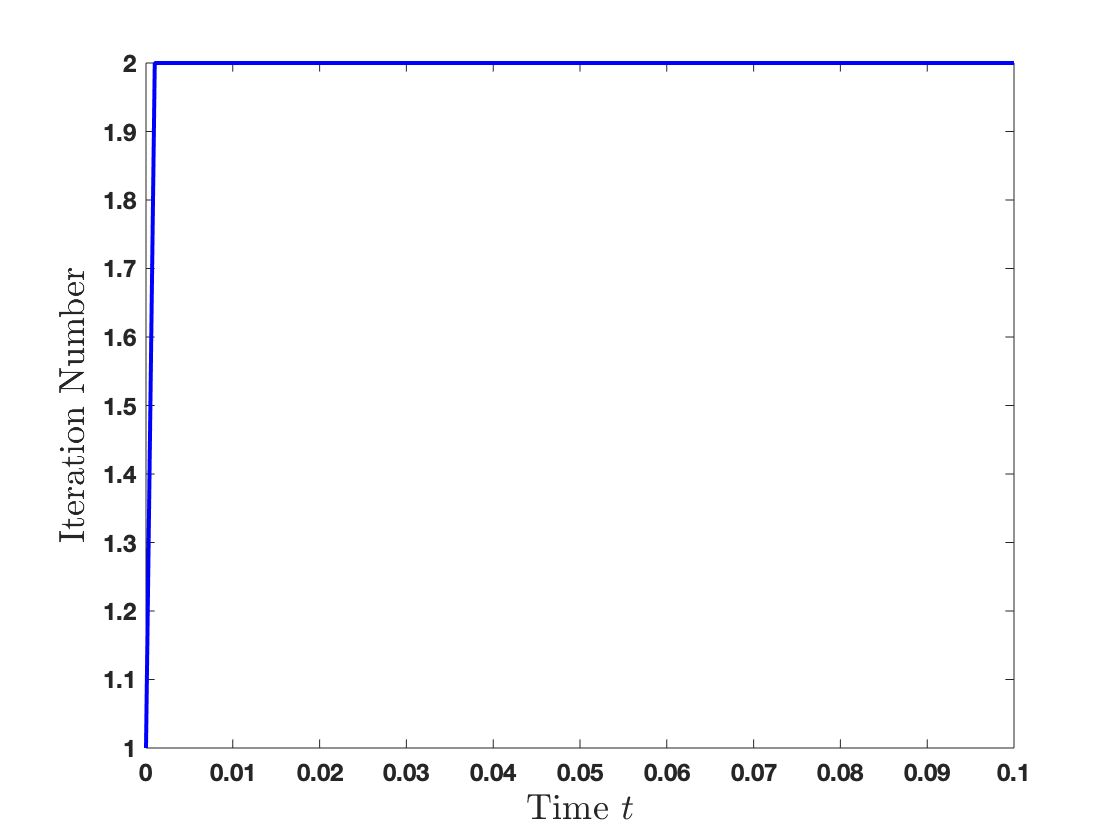}}
    \subfloat[BiCGstab LU iter]{\includegraphics[width=0.25\linewidth]{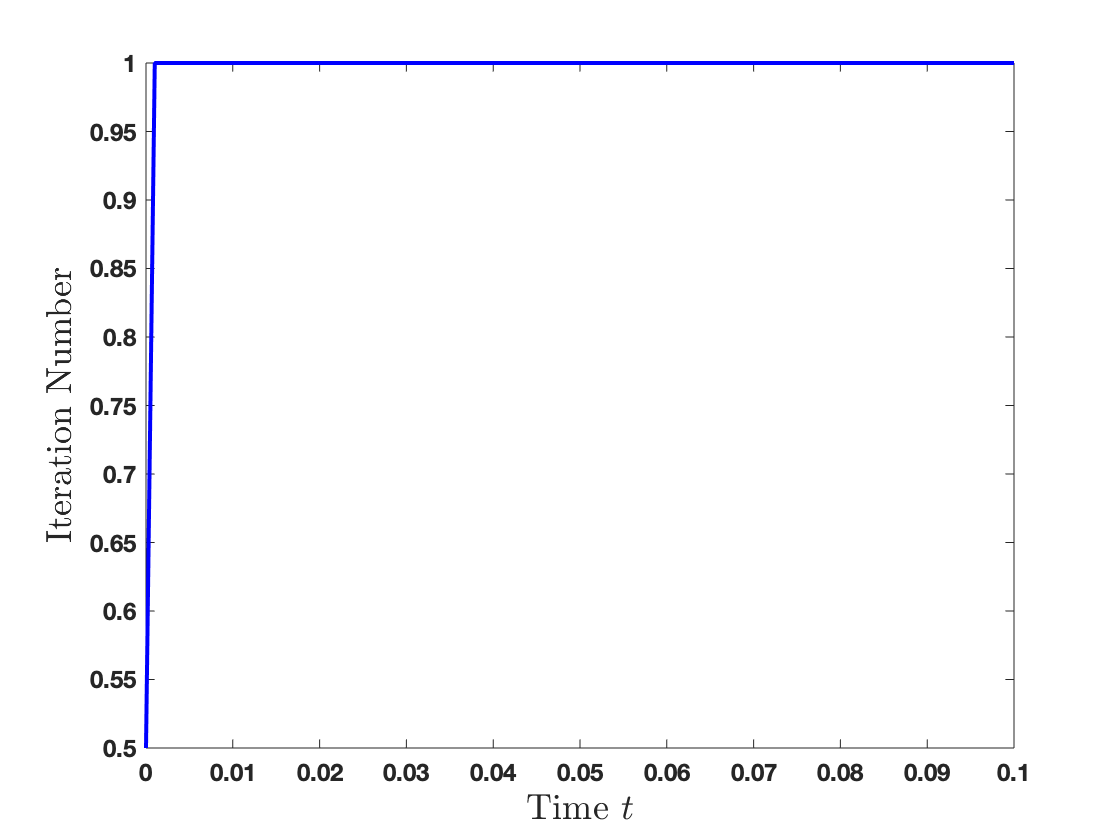}}
    % \hspace{0.1in}
    % \hspace{0.1in}
    \subfloat[TFQMR LU iter]{\includegraphics[width=0.25\linewidth]{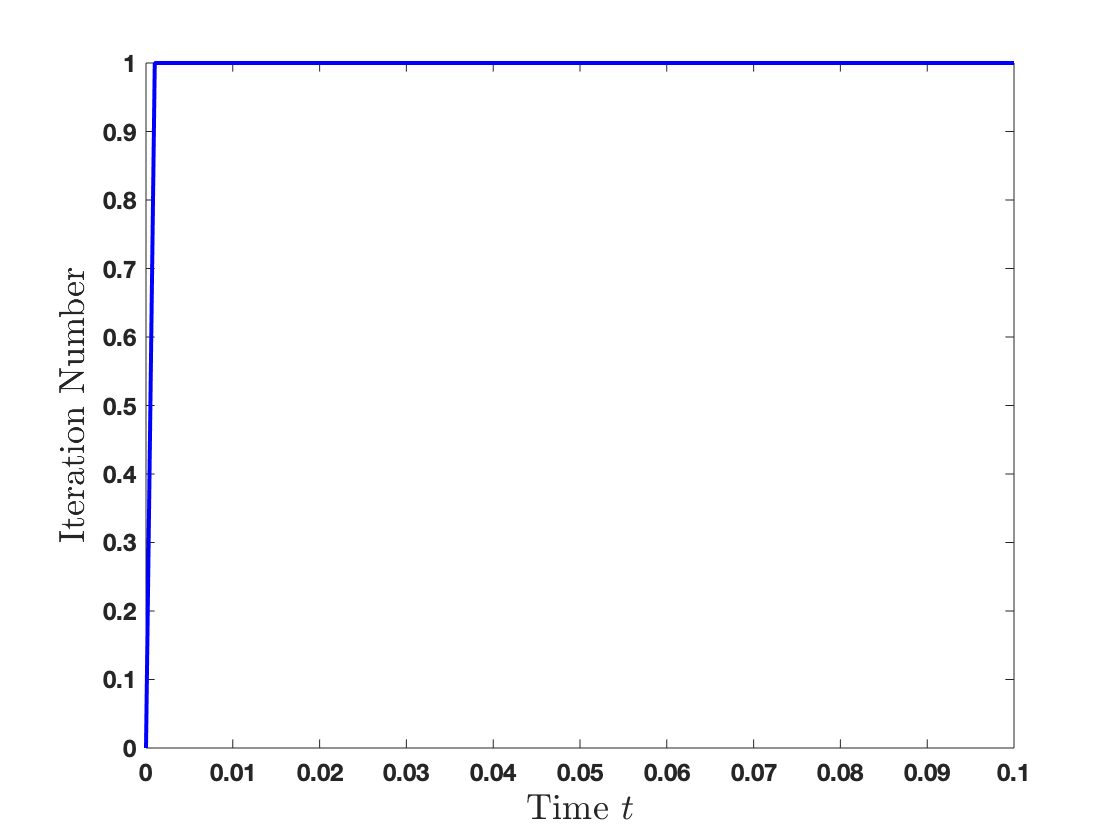}}
     \hspace{0.1in}
    \subfloat[GMRES LU relres]{\includegraphics[width=0.25\linewidth]{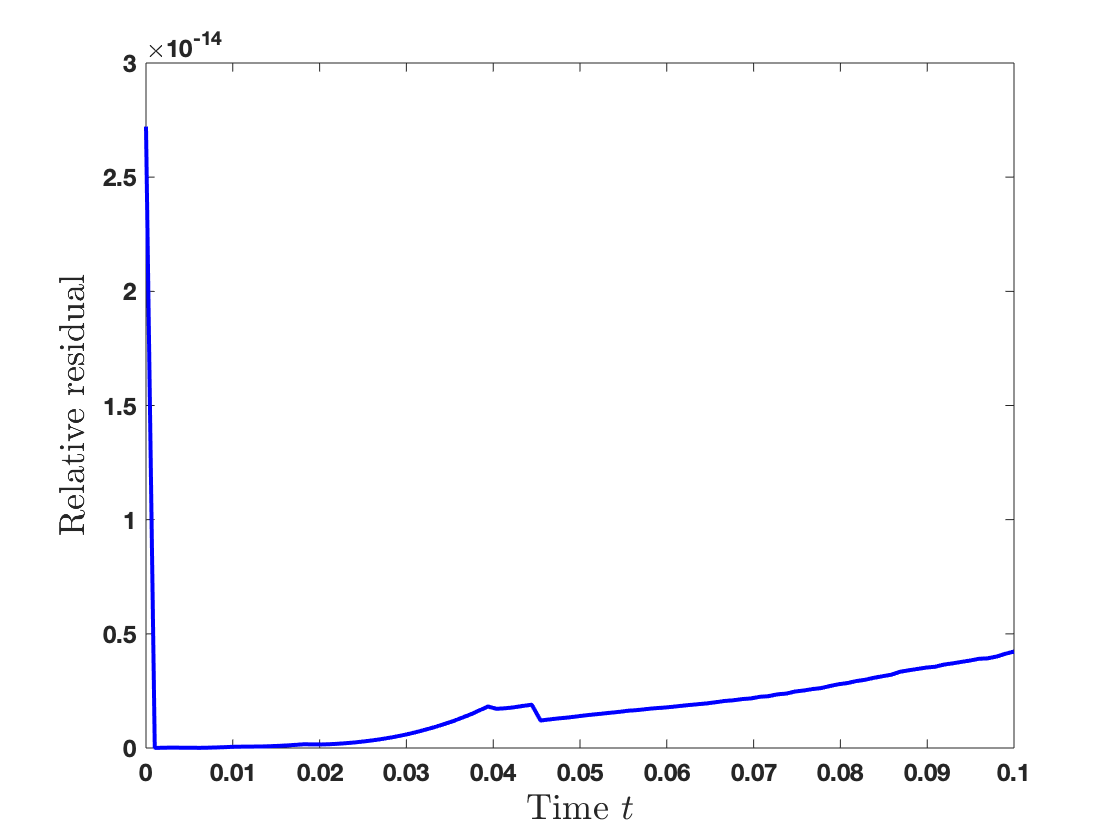}}
    \subfloat[BiCGstab LU relres]{\includegraphics[width=0.25\linewidth]{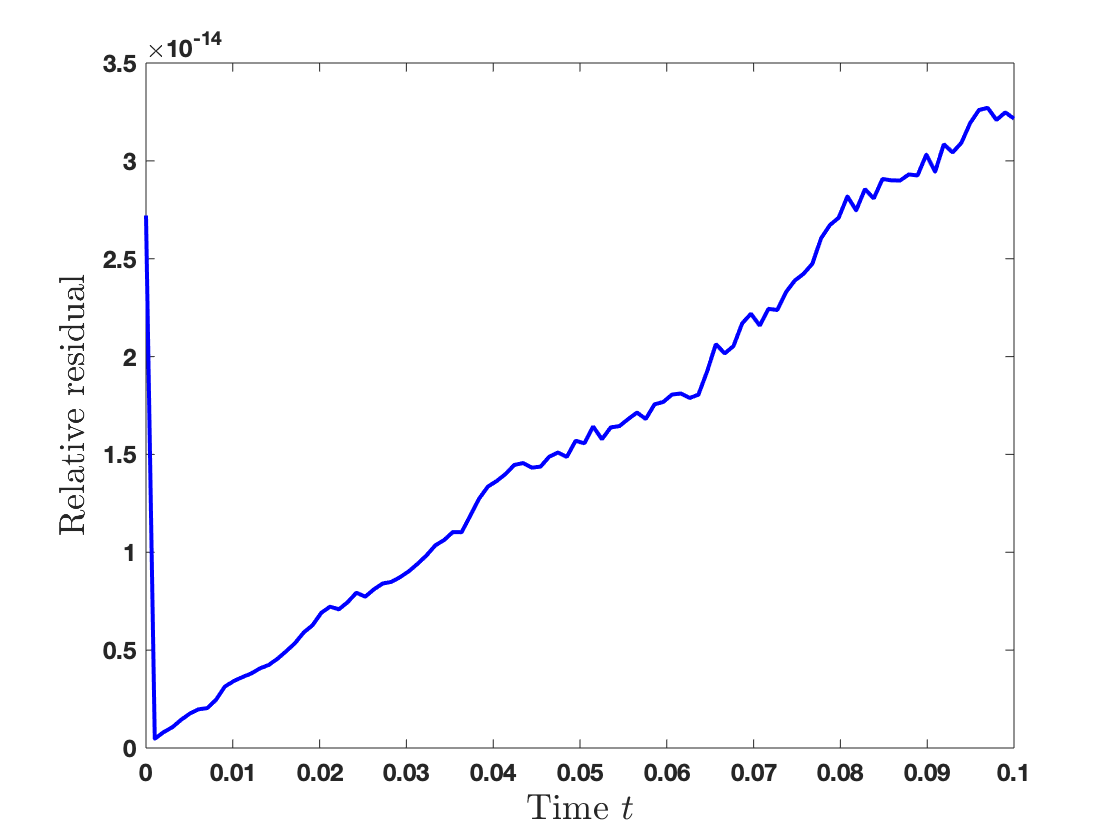}}
    % \hspace{0.1in}
    % \hspace{0.1in}
    \subfloat[TFQMR LU relres]{\includegraphics[width=0.25\linewidth]{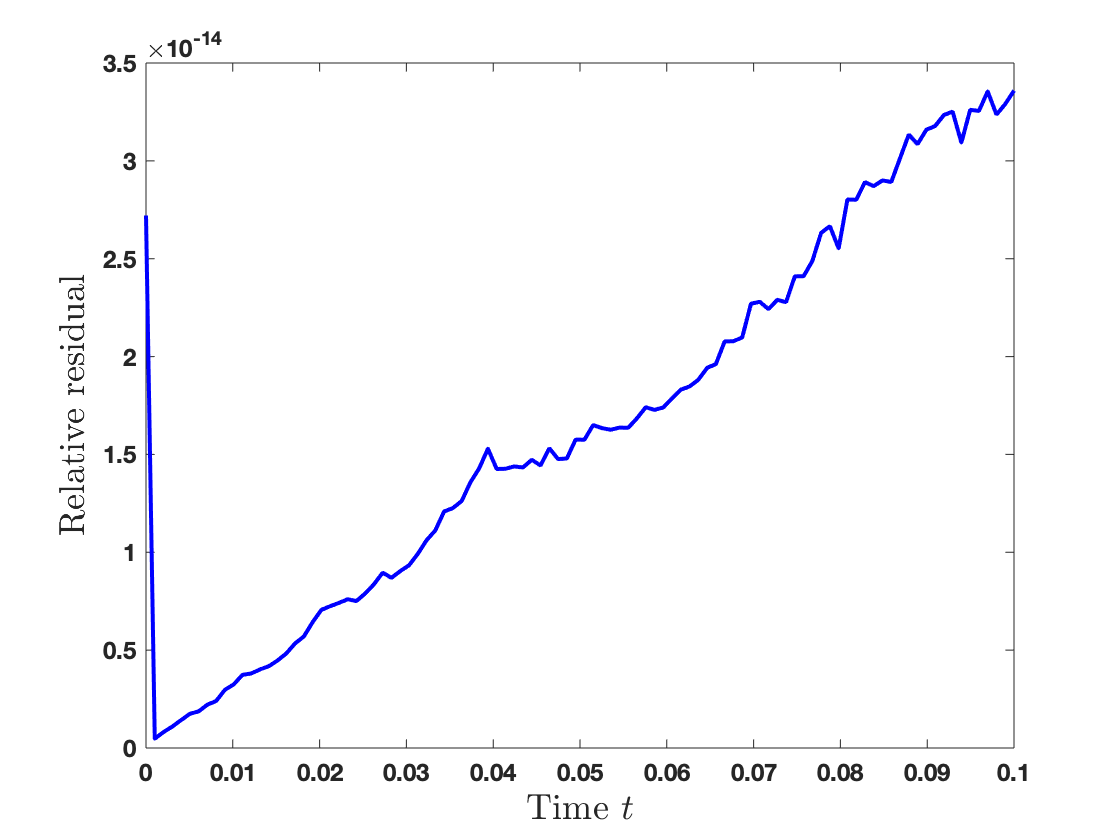}}
    \caption{Different linear solvers with $\alpha=0$, $N_x=1000$, $N_t=100$ up to the final time $T=0.1$ for 1D case. BDF1 without projection.}
    \label{fig:BDF1-2}
\end{figure}

\begin{figure}[htbp]
    \centering
    \subfloat[GMRES iter]{\includegraphics[width=0.25\linewidth]{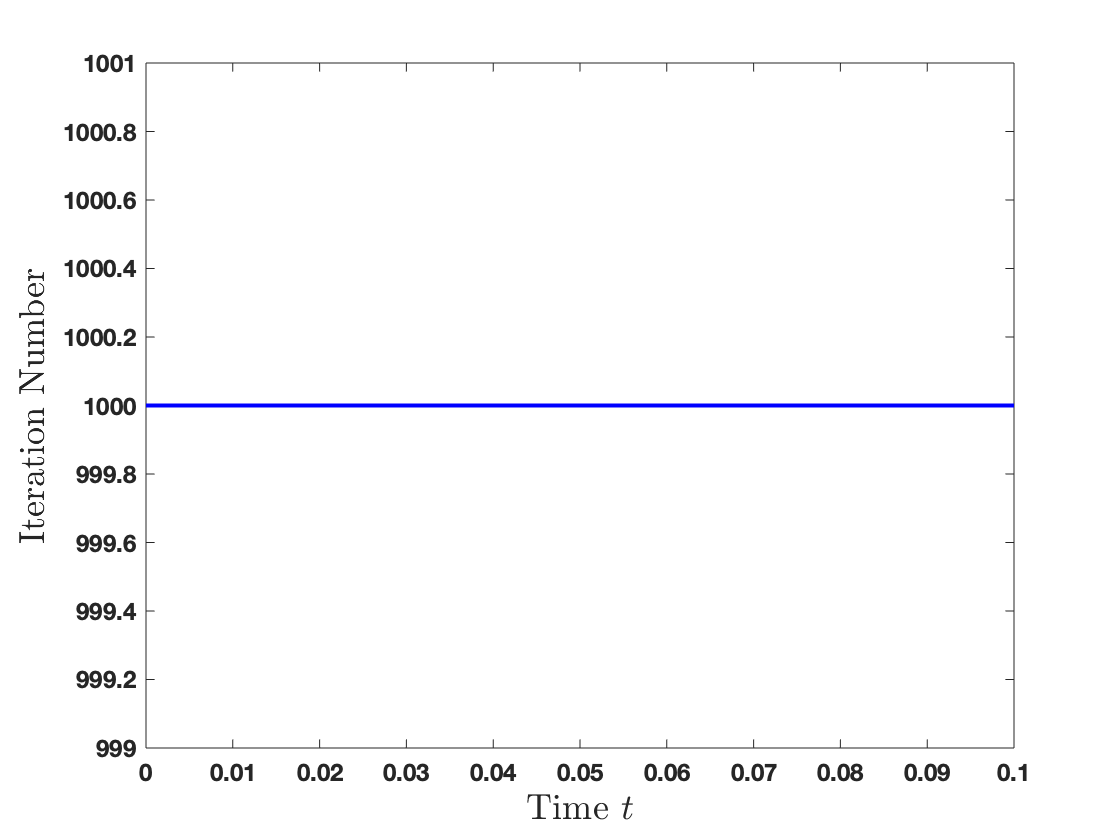}}
     \subfloat[BiCGstab iter]{\includegraphics[width=0.25\linewidth]{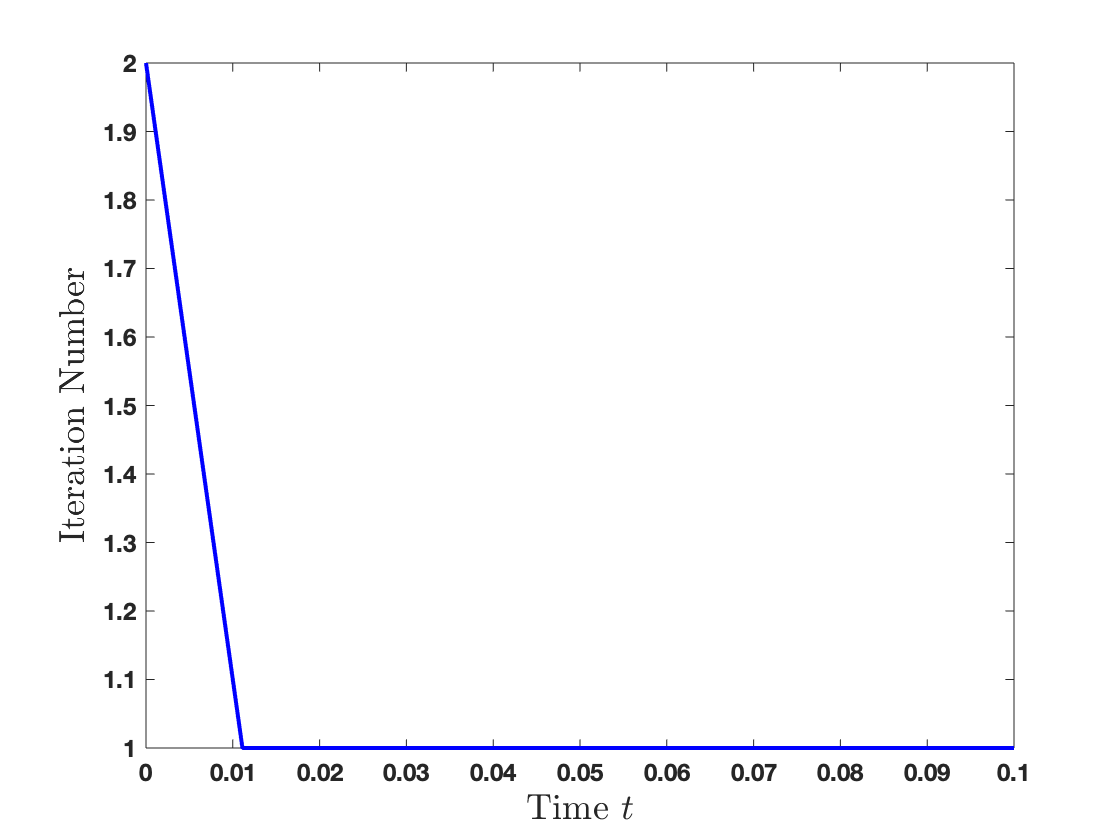}}
     \subfloat[TFQMR iter]{\includegraphics[width=0.25\linewidth]{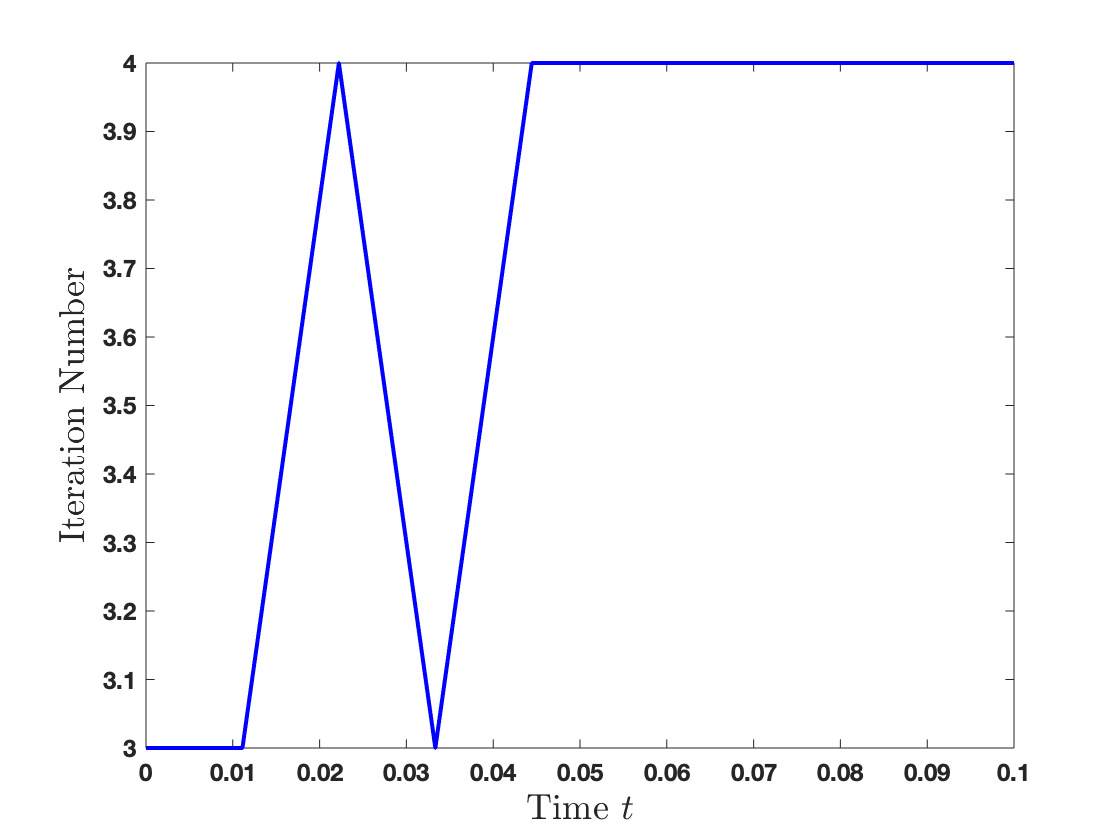}}
      \hspace{0.1in}
      \subfloat[GMRES relres]{\includegraphics[width=0.25\linewidth]{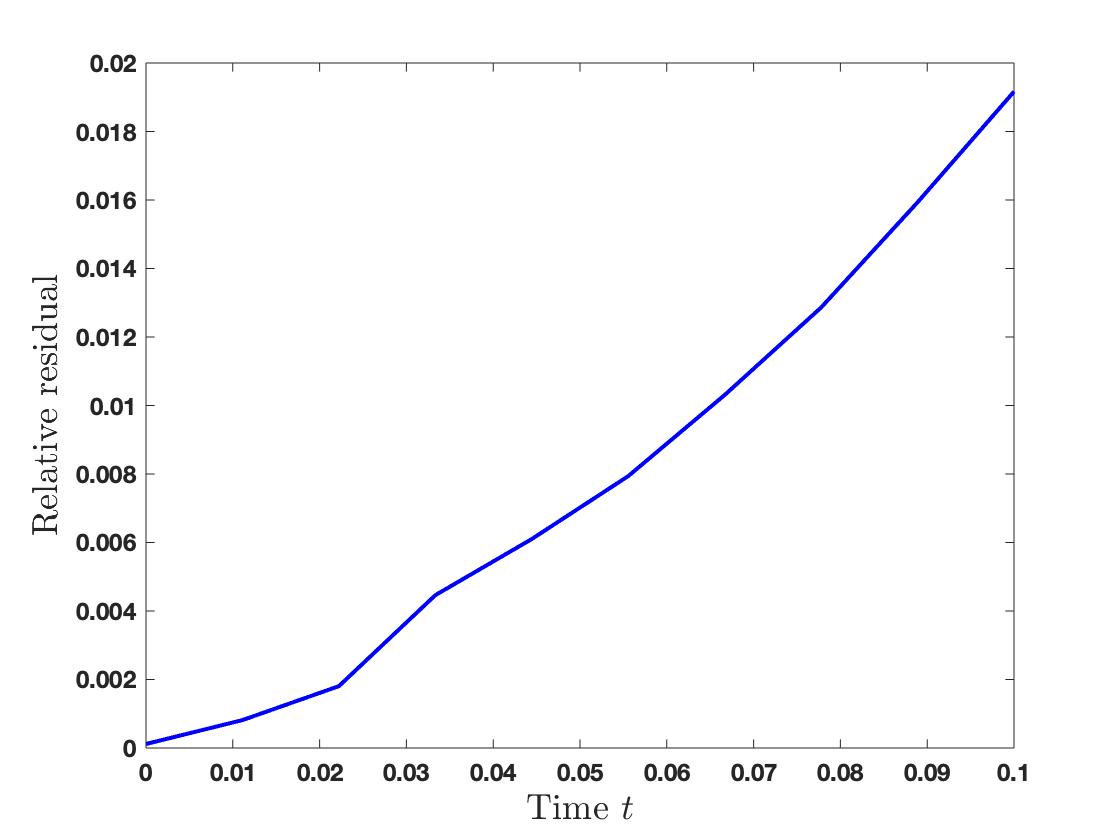}}
     \subfloat[BiCGstab relres]{\includegraphics[width=0.25\linewidth]{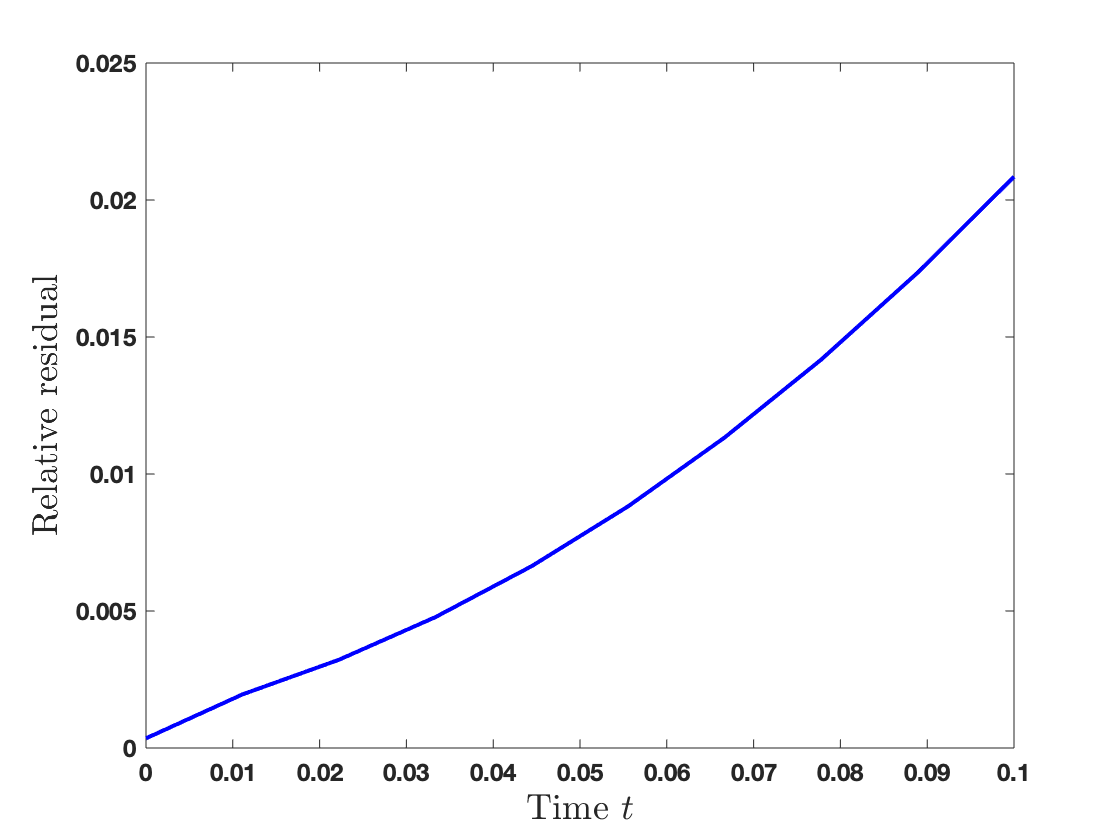}}
     \subfloat[TFQMR relres]{\includegraphics[width=0.25\linewidth]{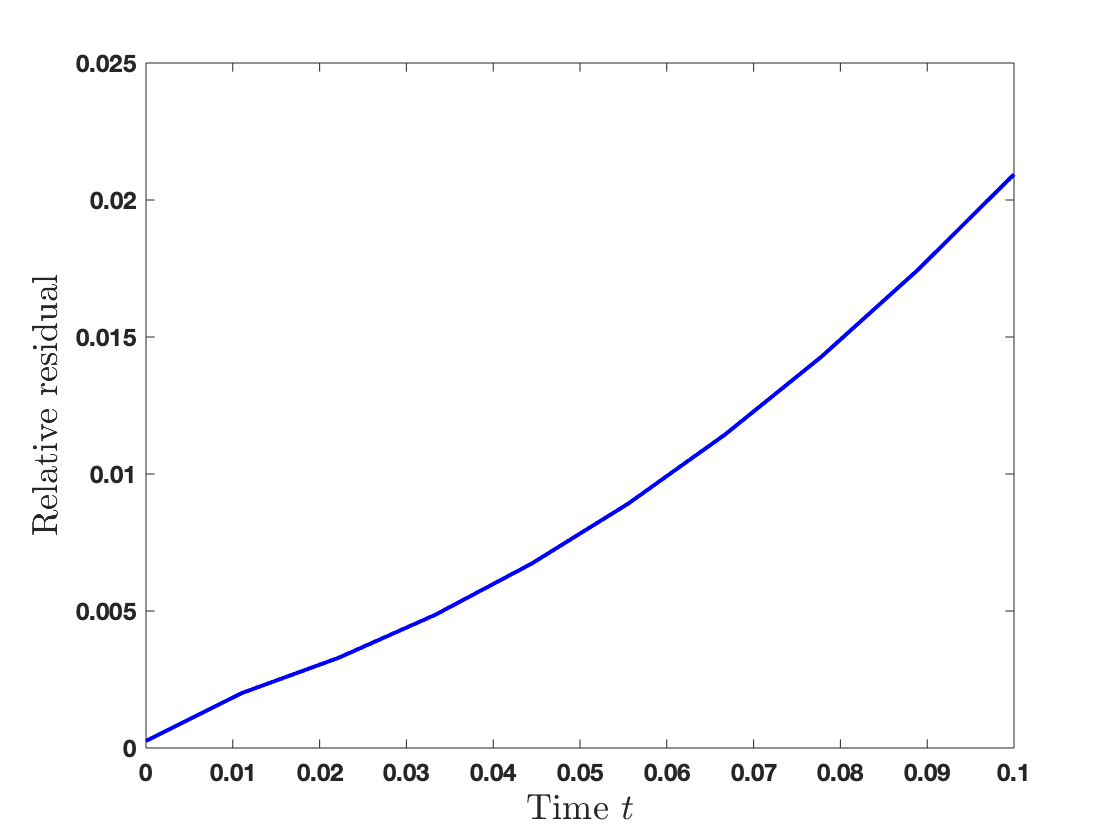}}
      \hspace{0.1in}
    \subfloat[GMRES with LU iter]{\includegraphics[width=0.25\linewidth]{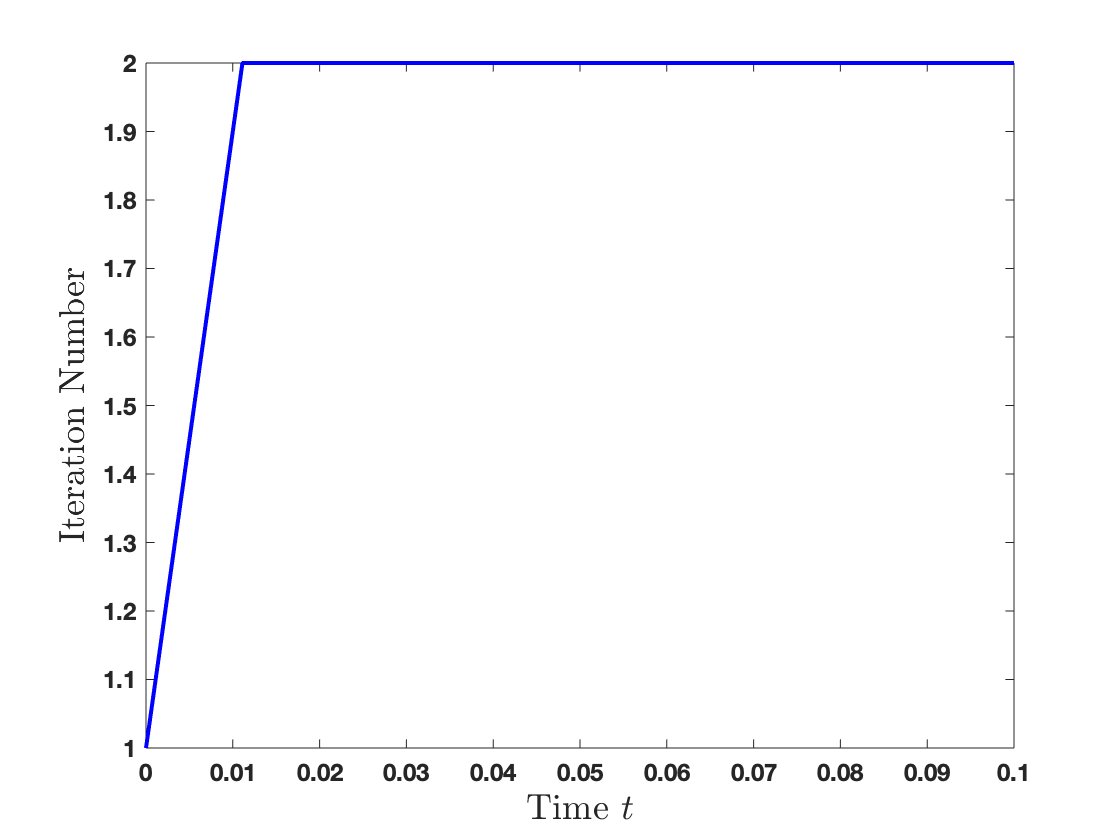}}
    \subfloat[BiCGstab with LU iter]{\includegraphics[width=0.25\linewidth]{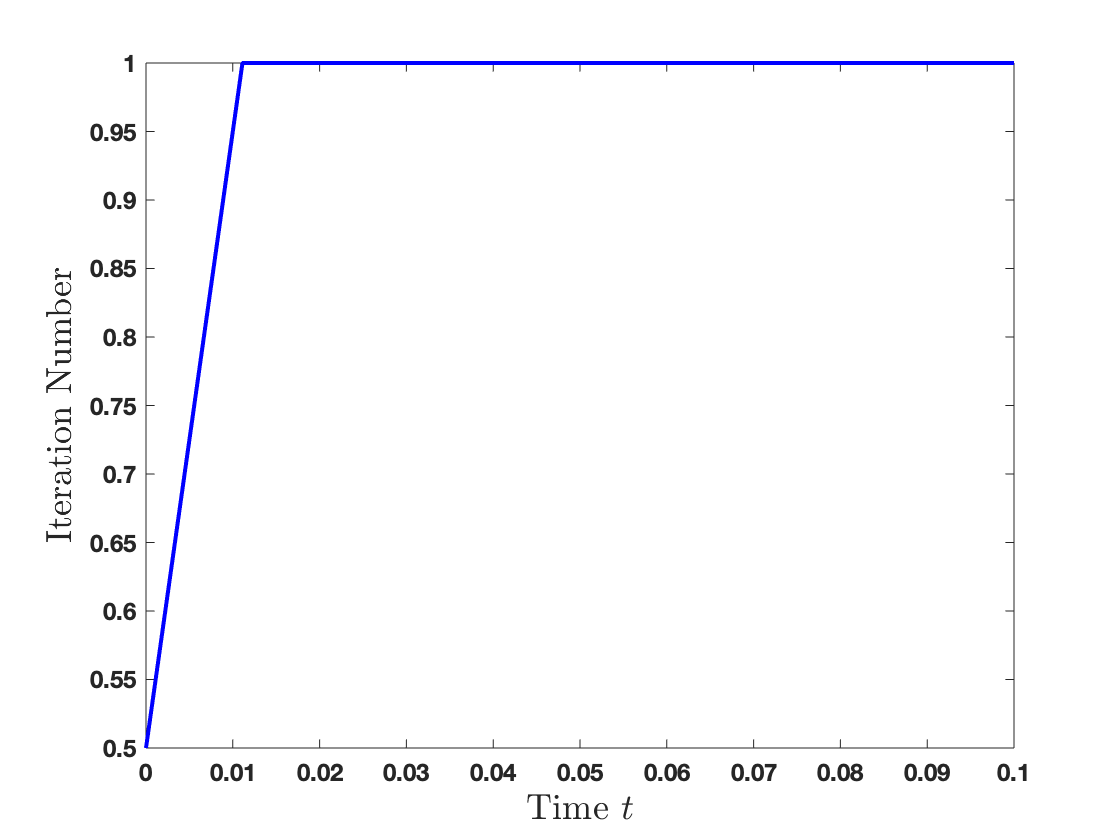}}
    % \hspace{0.1in}
    \subfloat[TFQMR with LU iter]{\includegraphics[width=0.25\linewidth]{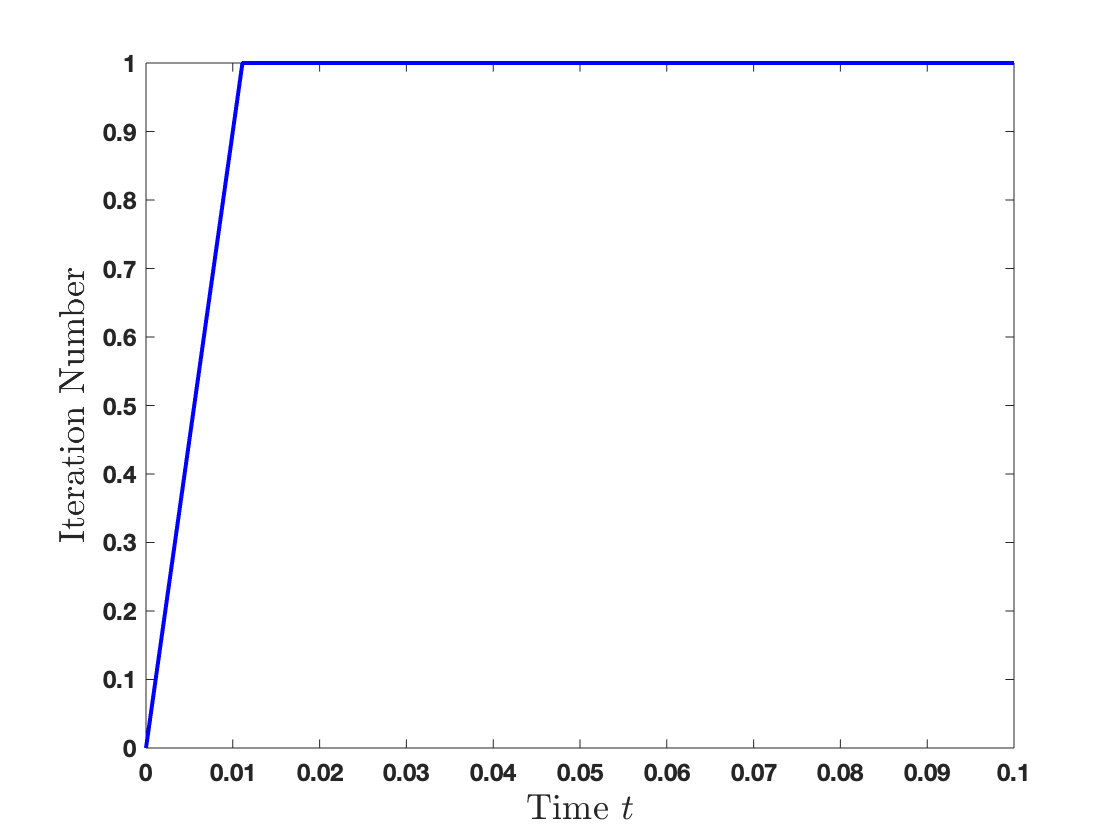}}
    \hspace{0.1in}
    \subfloat[GMRES with LU relres]{\includegraphics[width=0.25\linewidth]{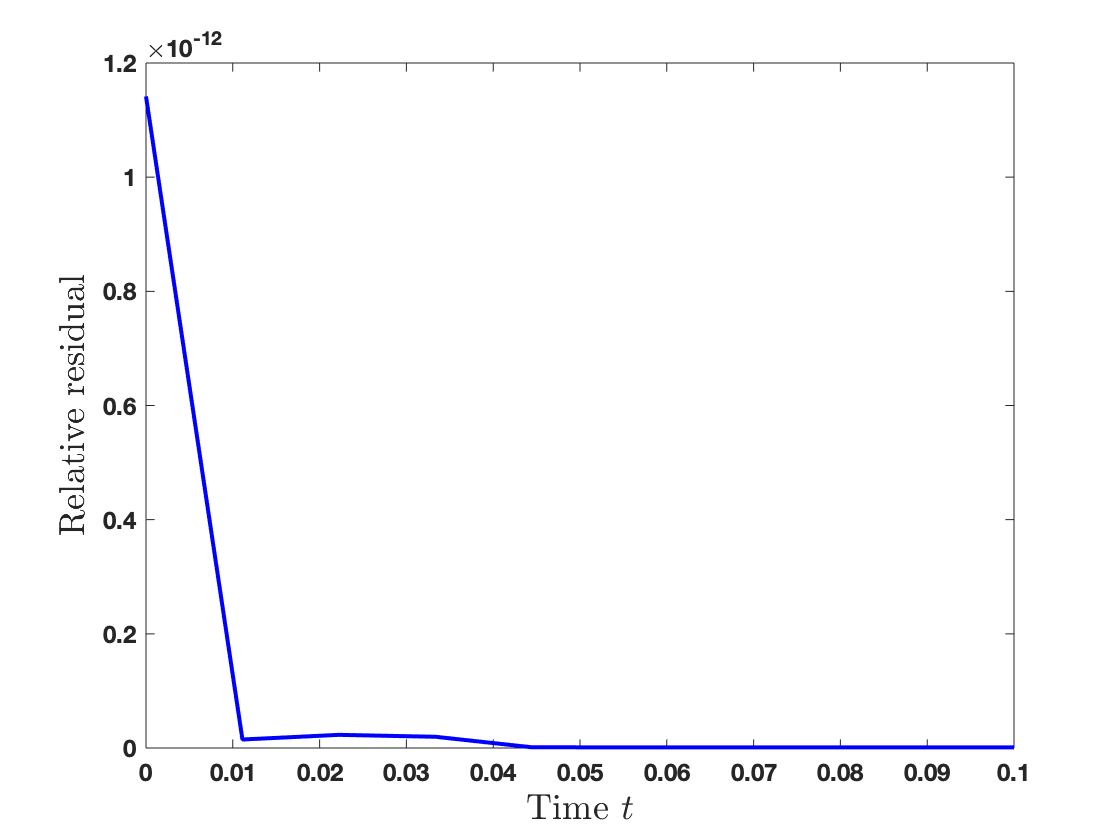}}
    \subfloat[BiCGstab with LU relres]{\includegraphics[width=0.25\linewidth]{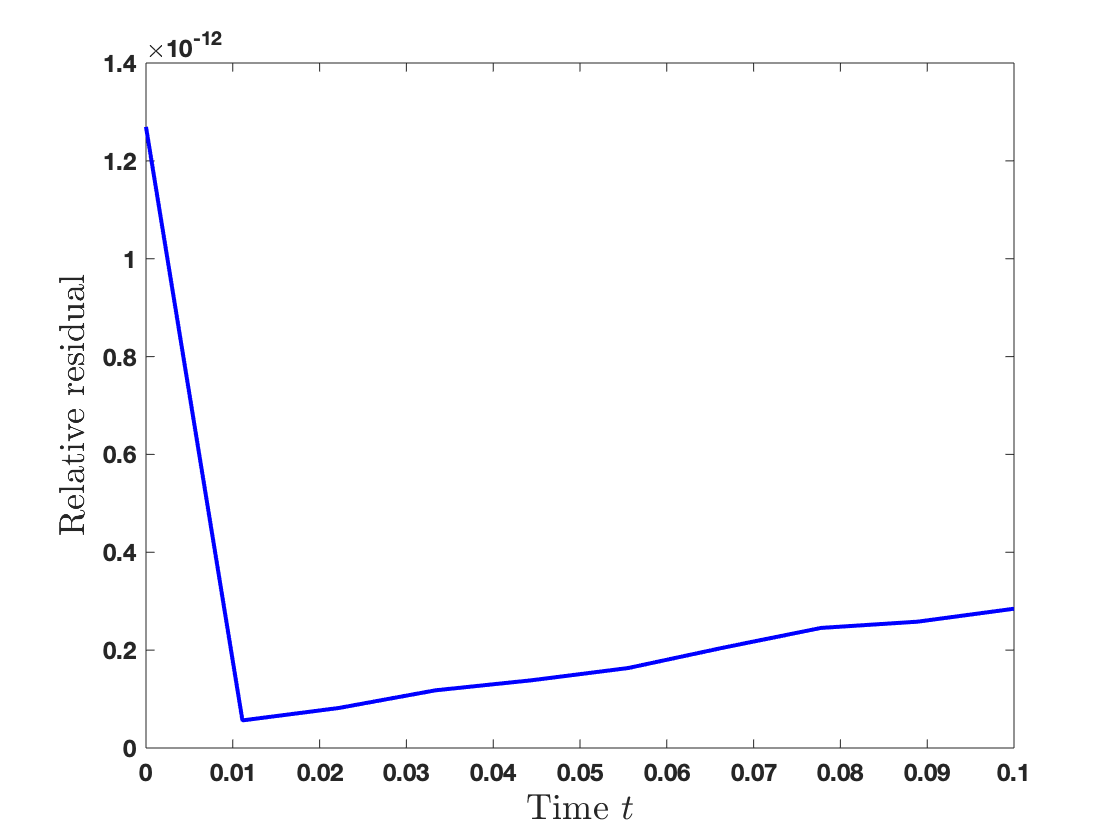}}
    % \hspace{0.1in}
    \subfloat[TFQMR with LU relres]{\includegraphics[width=0.25\linewidth]{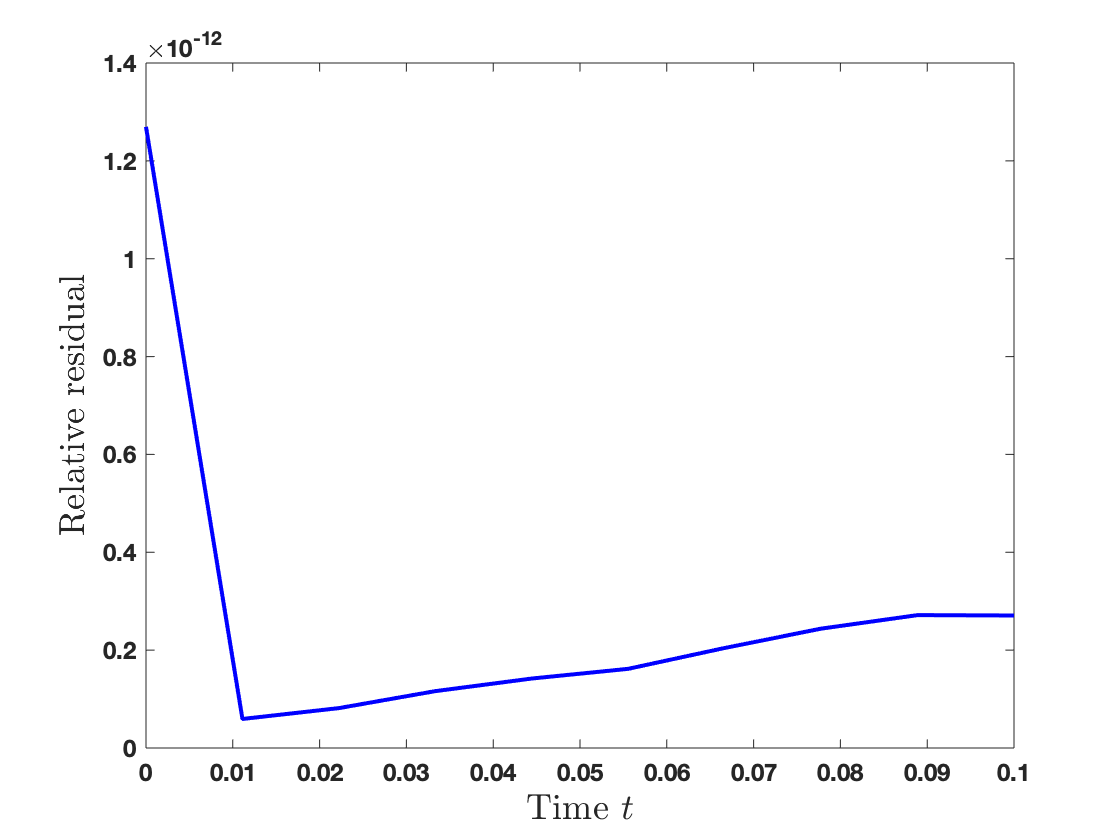}}
    \caption{Different linear solvers with $\alpha=0$, $N_x=1000$, $N_t=10$ up to the final time $T=0.1$ for 1D case. BDF1 without projection.}
    \label{fig:BDF1-3}
\end{figure}

\begin{figure}[htbp]
    \centering
    \subfloat[GMRES iter]{\includegraphics[width=0.25\linewidth]{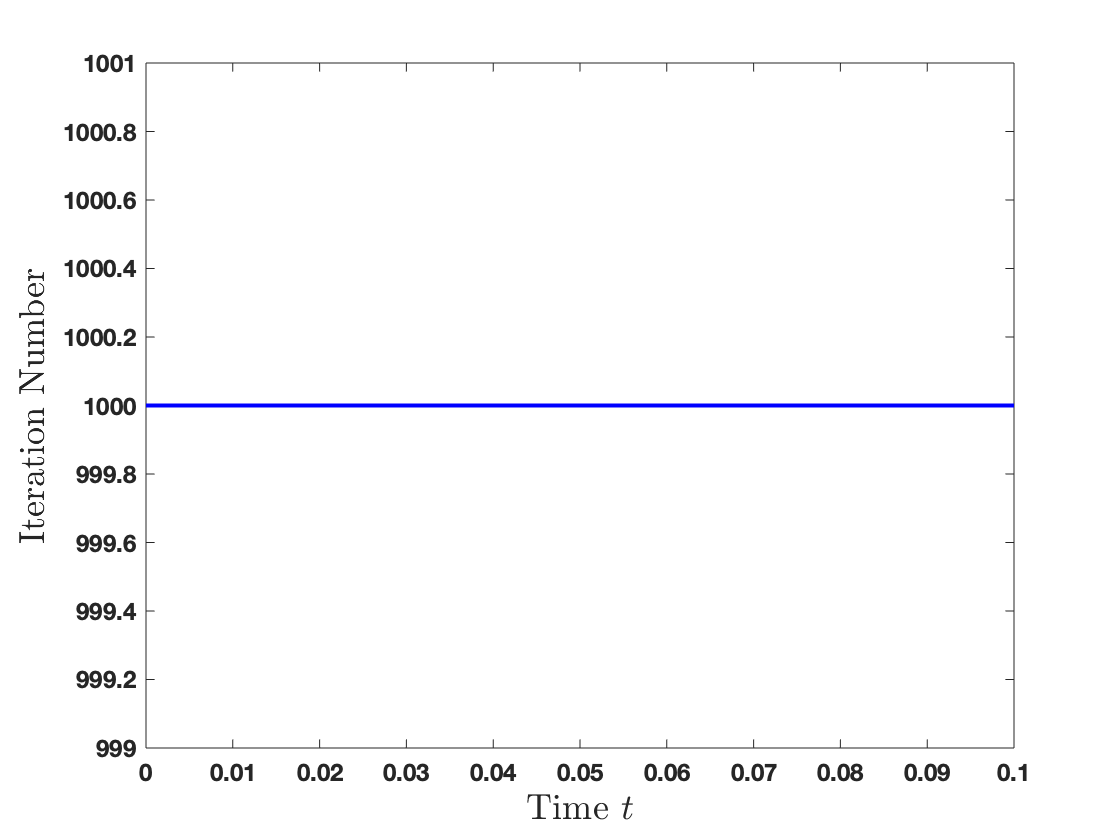}}
    \subfloat[BiCGstab iter]{\includegraphics[width=0.25\linewidth]{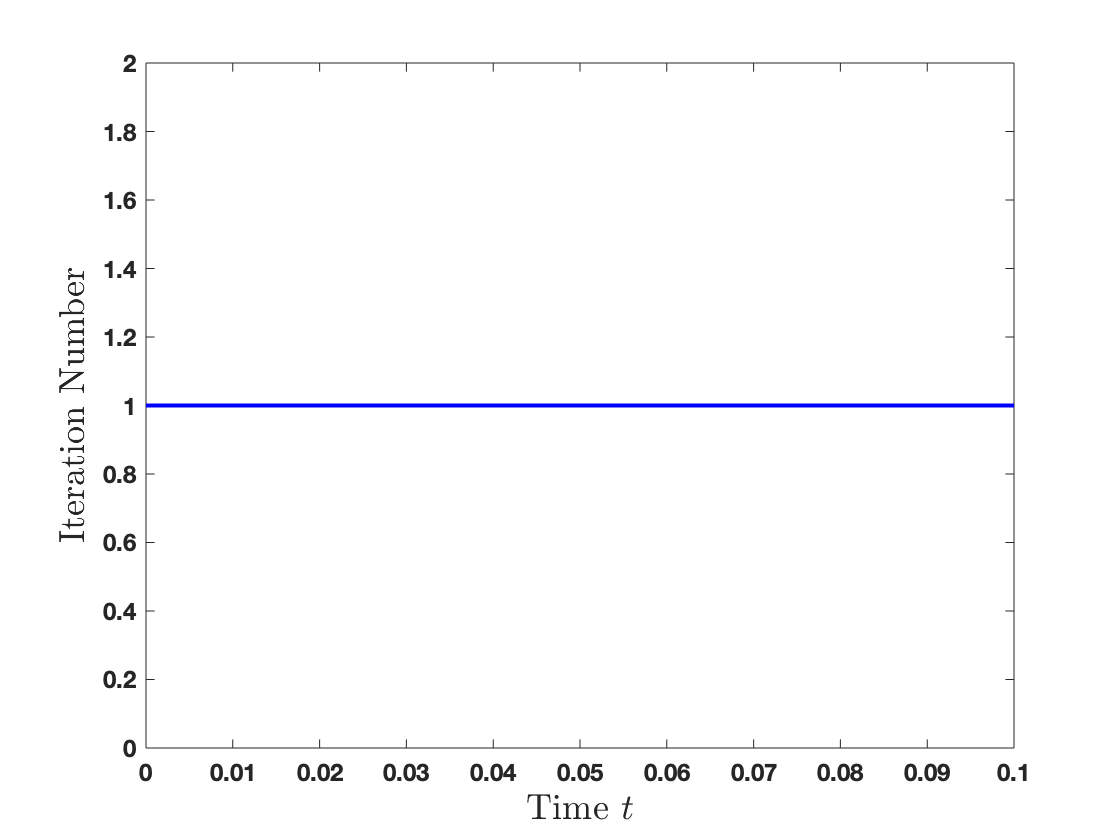}}
    % \hspace{0.1in}
    \subfloat[TFQMR iter]{\includegraphics[width=0.25\linewidth]{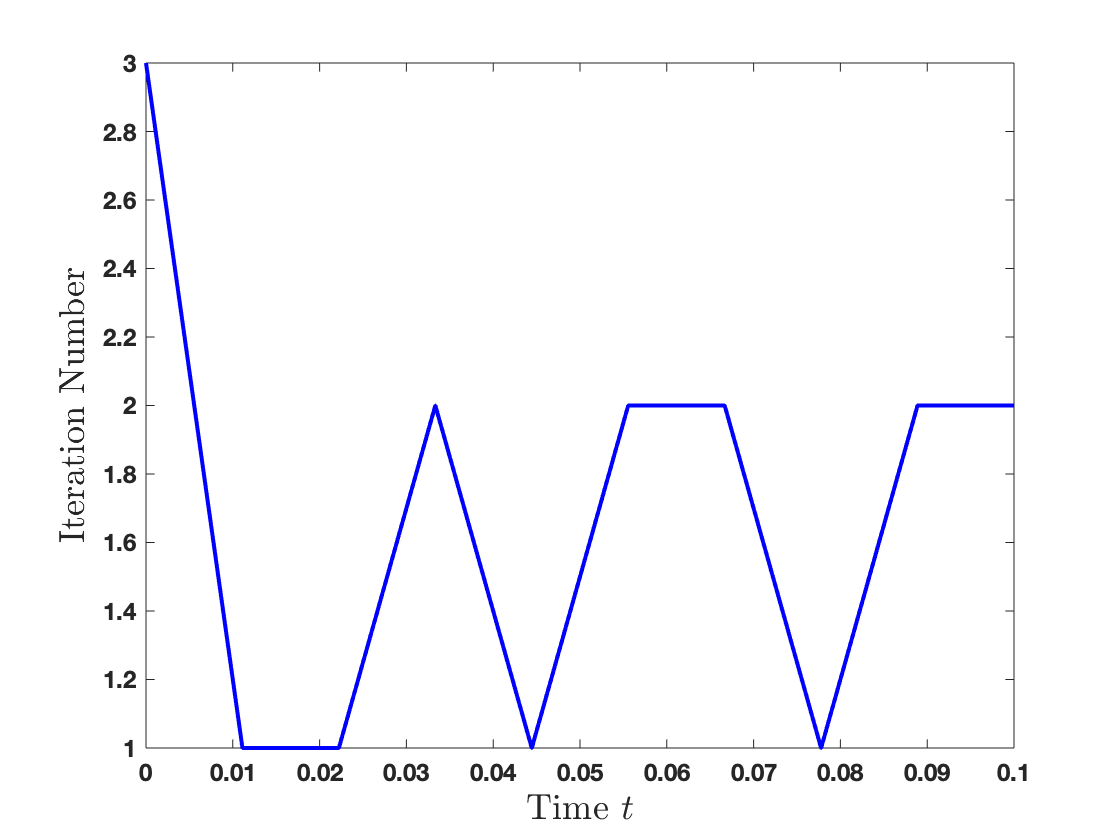}}
    \hspace{0.1in}
    \subfloat[GMRES relres]{\includegraphics[width=0.25\linewidth]{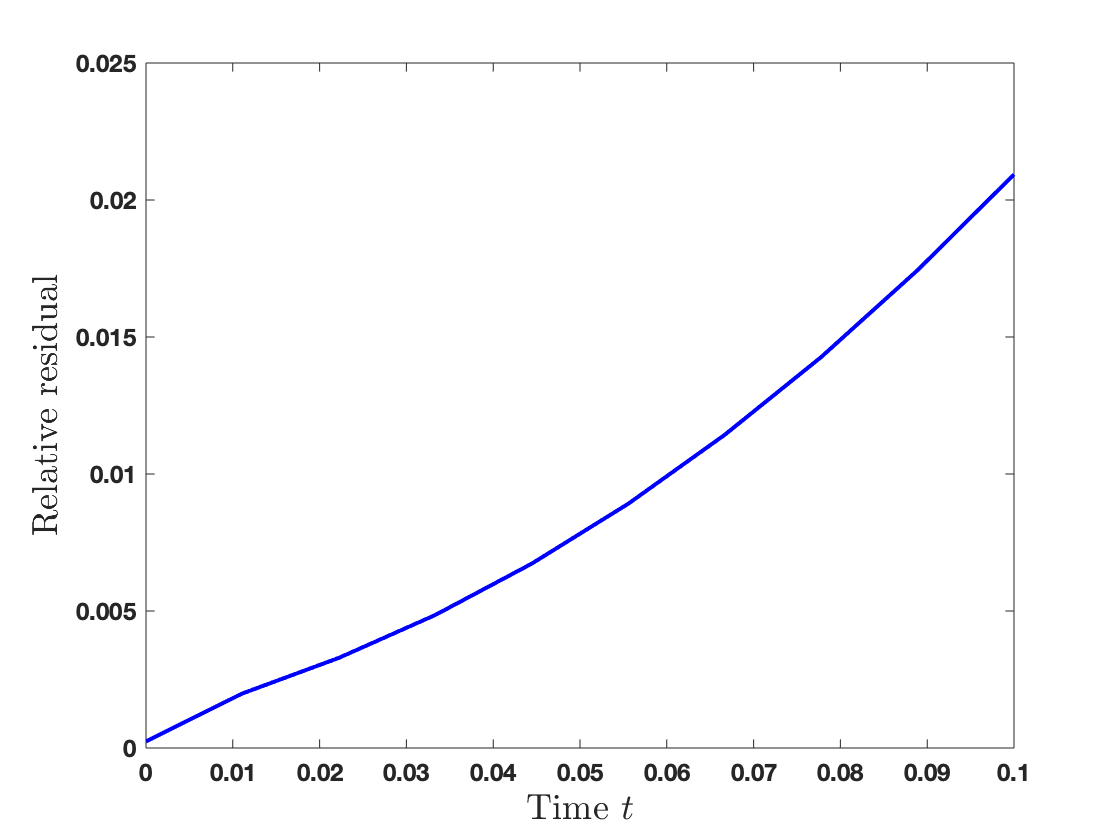}}
    \subfloat[BiCGstab relres]{\includegraphics[width=0.25\linewidth]{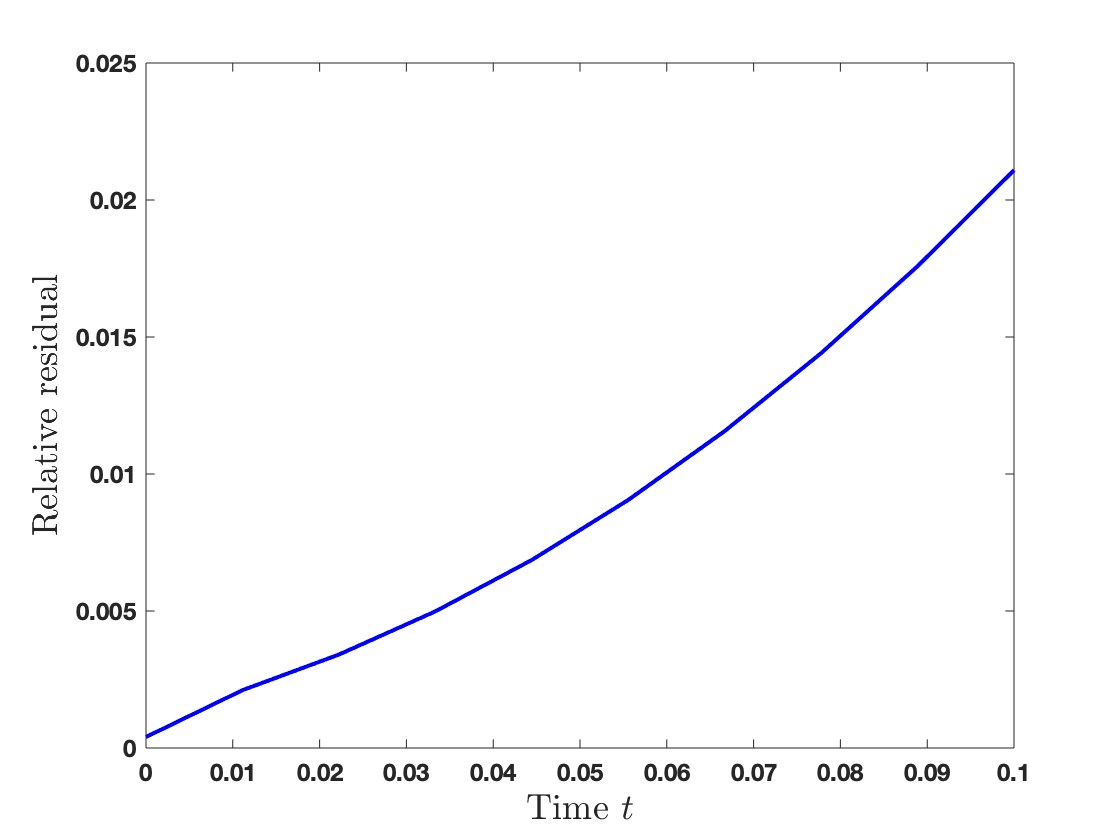}}
    % \hspace{0.1in}
    \subfloat[TFQMR relres]{\includegraphics[width=0.25\linewidth]{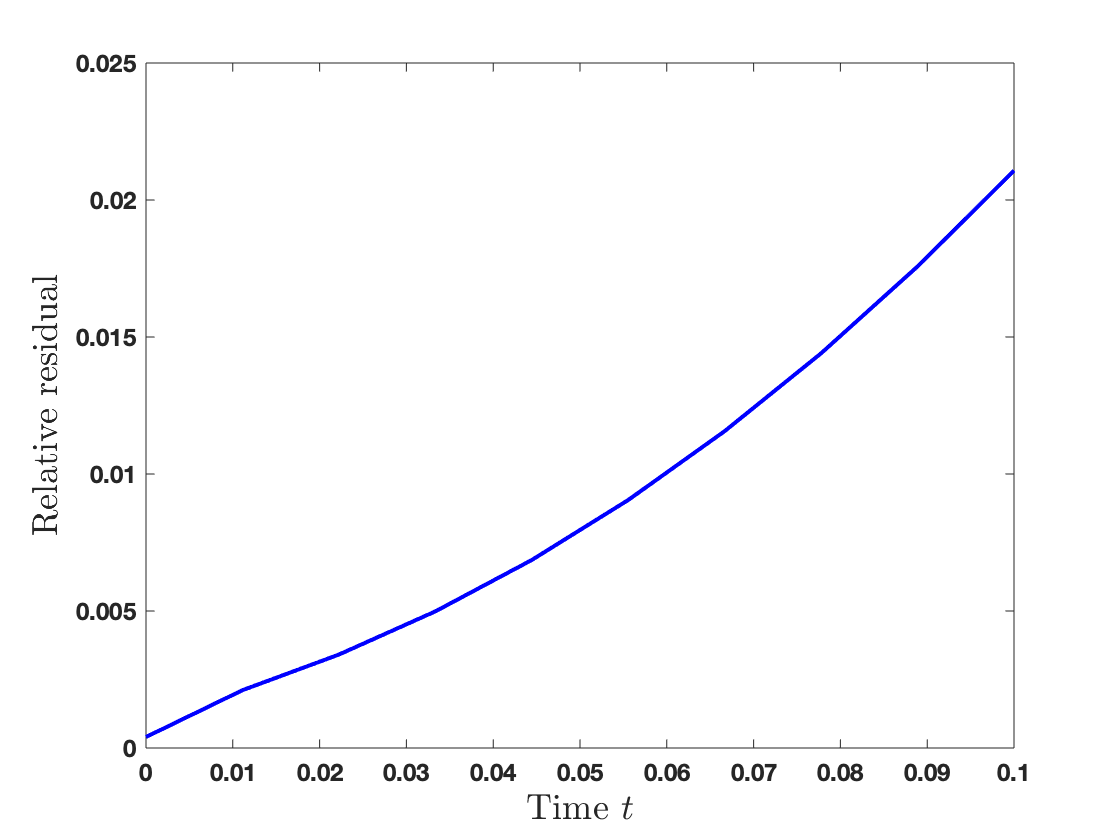}}
    \caption{Different linear solvers, namely GMRES, BiCGstab, TFQMR with $\alpha=0$, $N_x=10000$, $N_t=10$ up to the final time $T=0.1$ for the 1D case. BDF1 without projection.}
    \label{fig:BDF1-4}
\end{figure}

We test several iterative methods, like GMRES, BiCGstab and TFQMR and its preconditioning using LU decomposition. The results are presented in \Cref{fig:BDF1-1} for different linear solvers with $\alpha=0$, $N_x=64$, $N_t=100$ up to the final time $T=0.1$ in 1D case without projection.  Iteration counts and relative residuals are plotted for both original runs and those preconditioned by LU factorization. TFQMR yields nearly constant iteration numbers yet suffers from continuously growing residuals. GMRES and BiCGStab exhibit strongly oscillatory iteration counts and residuals. LU preconditioning drastically reduces iteration numbers and stabilizes residual behaviour across all three solvers.
\Cref{fig:BDF1-2} evaluates GMRES, BiCGStab and TFQMR for a 1D BDF1 discretized PDE with \(N_x=1000\), \(N_t=100\) and final time \(T=0.1\). Iteration counts and relative residuals over time are presented for both unpreconditioned and LU‑preconditioned variants. Without preconditioning, TFQMR exhibits severe oscillatory iteration numbers, while GMRES and BiCGStab keep steady iterations. LU preconditioning stabilizes iteration counts for all solvers yet fails to suppress residual growth, indicating residual degradation persists despite stabilized iterative behaviour. \Cref{fig:BDF1-3} compares GMRES, BiCGStab and TFQMR for a 1D BDF1‑discretized PDE, with \(N_x=1000\), \(N_t=10\), \(T=0.1\) and \(\alpha=0\). Iteration numbers and relative residuals over time are shown for unpreconditioned and LU‑preconditioned setups. Unpreconditioned TFQMR displays sharp iteration spikes. LU preconditioning yields nearly constant low iteration counts for all three solvers and drastically reduces initial residuals. Nevertheless, slow residual growth still occurs as time evolves, demonstrating that LU cannot fully eliminate time‑dependent residual accumulation. \Cref{fig:BDF1-4} investigates GMRES, BiCGStab and TFQMR for a 1D BDF1‑discretized PDE with \(N_x=10000\), \(N_t=10\), \(T=0.1\) and \(\alpha=0\). Iteration counts and relative residuals versus time are plotted for the unpreconditioned solvers. GMRES and BiCGStab maintain stable iteration numbers, whereas TFQMR produces pronounced oscillatory iteration behaviour. All three solvers exhibit continuous growth in relative residual over time. Even with refined spatial resolution, residual accumulation persists, and TFQMR remains susceptible to iterative oscillations for this problem.

\Cref{tab:num-1} reports wall‑clock computation time for GMRES, BiCGStab and TFQMR (with and without LU preconditioning) under various spatial \(N_x\) and temporal \(N_t\) resolutions for the 1D BDF1 scheme. Computational cost rises sharply with increasing \(N_x\). Unpreconditioned BiCGStab and TFQMR generally outperform GMRES. LU preconditioning accelerates small‑size problems, yet introduces substantial overhead for large \(N_x\). For fine spatial grids, pure iterative solvers are more time‑efficient than their LU‑preconditioned counterparts.

\begin{table}[htbp]
    \centering
    \begin{tabular}{|c|c|c|c|}
    \hline
        Method & $N_x$ & $N_t$& Total wall time (s) \\
        \hline
        GMRES & 64&100&  0.959815\\ 
        &1000&100&350.101823\\
        &1000&10&33.313915\\
        &10000&10&436.689131\\
        \hline
        GMRES with LU & 64&100&  1.034561\\ 
        & 1000&100&348.713269\\
        &1000&10&43.978948\\
        \hline
        BiCGstab & 64&100&  1.175190\\ 
        &1000&100&3.683081\\
        &1000&10&0.935558\\
        &10000&10&50.309583\\
        \hline
        BiCGstab with LU & 64&100&  0.630437\\ 
        &1000&100&355.375340\\
        &1000&10&40.506224\\
        \hline
        TFQMR & 64&100&  33.585289\\ 
        &1000&100&3.505883\\
        &1000&10&0.664085\\
        &10000&10&47.619316\\
        \hline
        TFQMR with LU & 64&100&  0.686889\\
        &1000&100&347.183437\\
        &1000&10&39.640307\\
        \hline
    \end{tabular}
    \caption{The computational cost using GMRES (with LU preconditioner), BiCGstab (with LU preconditioner) and TFQMR (with LU preconditioner) in the 1D case with different $N_x,N_t$, and $\alpha=0$. BDF1 without projection.}
    \label{tab:num-1}
\end{table}

Due to the preconditioning is tricky, we propose to use an aggregation AMG (AGMG) and compare it with GMRES. \Cref{tab:num-2} compares wall‑clock runtime of GMRES and AGMG for the 1D BDF1 discretization with \(N_t=100\), \(\alpha=0\) and two spatial resolutions \(N_x=100\) and 1000. Both solvers incur higher computational cost as spatial grid size increases. AGMG delivers drastically shorter execution time across both mesh settings. When \(N_x=1000\), GMRES suffers severe time overhead, while AGMG maintains efficiency. This demonstrates that the AGMG multigrid solver is far more competitive than standard GMRES for this specific PDE problem.

\begin{table}[htbp]
    \centering
    \begin{tabular}{|c|c|c|c|}
    \hline
        Method & $N_x$ & $N_t$& Total wall time (s) \\
        \hline
        GMRES &100 &100& 3.086975 \\ 
         &1000 &100& 353.388574 \\ 
        \hline
        AGMG &100 &100& 0.571916 \\
          &1000 &100& 3.899783 \\
        \hline
    \end{tabular}
    \caption{The computational cost using GMRES, AGMG in the 1D case with different $N_x=100,1000,N_t=100$, and $\alpha=0$. BDF1 without projection.}
    \label{tab:num-2}
\end{table}

\begin{figure}[htbp]
    \centering
    \subfloat[GMRES]{\includegraphics[width=0.5\linewidth]{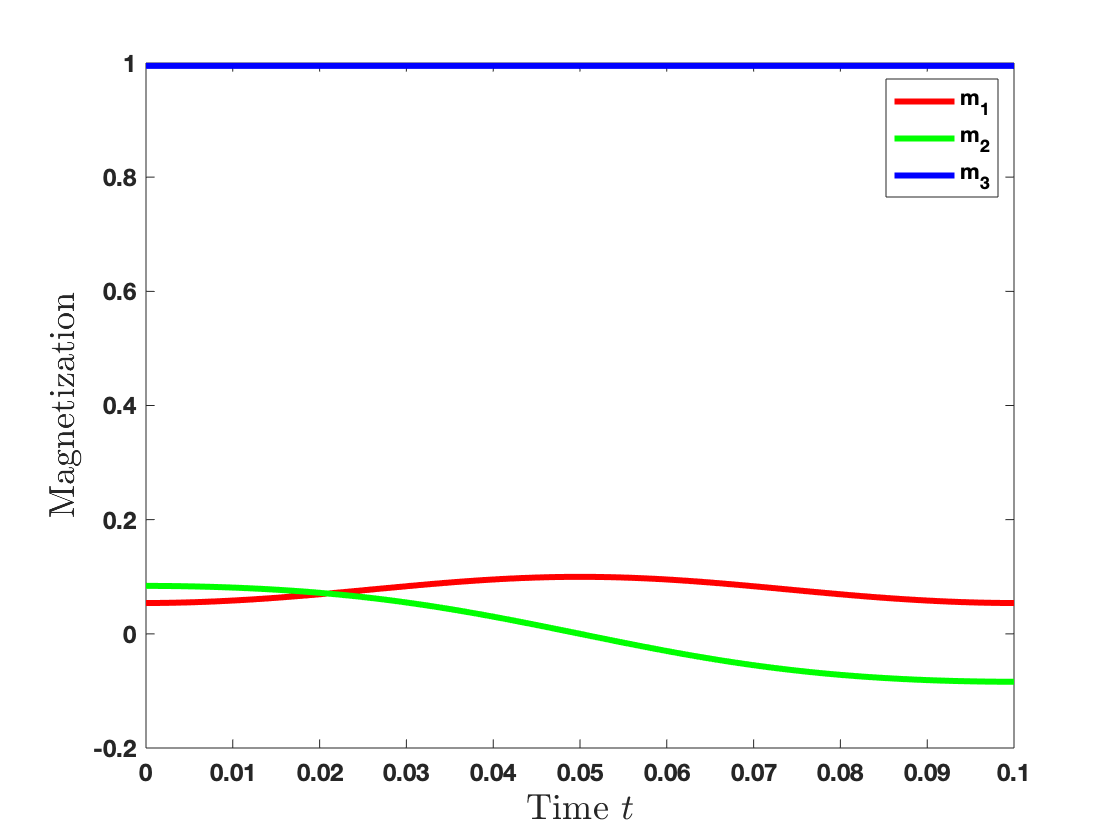}}
     \subfloat[AGMG]{\includegraphics[width=0.5\linewidth]{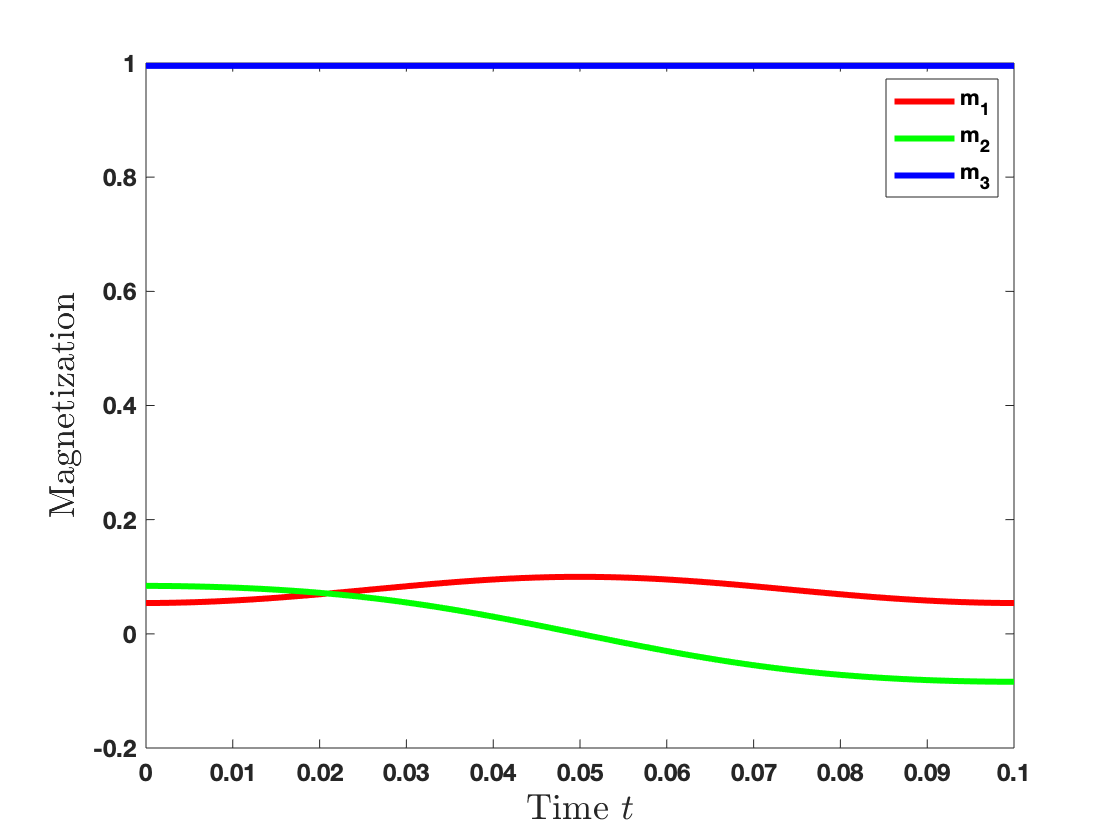}}
    \caption{Solution profiles in the 1D case with $N_x=100$ and $N_t=100$ up to the final time $T=0.1$ with $\alpha=0$. Using GMRES and AGMG. The initial condition $\m_0=[\cos(\cos(\pi x))\sin(0),\sin(\cos(\pi x))\sin(0),\cos(0)]^T$. BDF1 without projection.}
    \label{fig:num-1}
\end{figure}

Given the initial conditions as below,
\begin{align*}
    \m_0&=[\cos(\cos(\pi x))\sin(0),\sin(\cos(\pi x))\sin(0),\cos(0)]^T\quad \text{ in 1D,}\\
      \m_0&=[\cos(\cos(\pi x)\cos(\pi y))\sin(0),\sin(\cos(\pi x)\cos(\pi y))\sin(0),\cos(0)]^T \quad \text{ in 2D}\\
    \m_0&=[\cos(\cos(\pi x)\cos(\pi y)\cos(\pi z))\sin(0),\sin(\cos(\pi x)\cos(\pi y)\cos(\pi z))\sin(0),\cos(0)]^T \quad \text{ in 3D},
\end{align*}
we get the solution profiles at $T=0.1$ in \Cref{fig:num-1} which plots magnetization components \(m_1\), \(m_2\), \(m_3\) solved by GMRES and AGMG for the 1‑D BDF1 simulation with \(N_x=100\), \(N_t=100\), \(T=0.1\) and \(\alpha=0\). The two subplots show nearly identical time‑evolution curves for all three magnetization components. Despite their large differences in computational efficiency, GMRES and AGMG produce matching numerical solutions for this test problem. This confirms that AGMG can achieve the same solution quality as GMRES while greatly reducing runtime.

\Cref{tab:num-3} presents wall‑clock runtime of GMRES and AGMG for the 3‑D BDF1 simulation with projection, where \(N_t=40\), \(\alpha=0\), and uniform grids \(N_x=N_y=N_z = 20\) and 30. Both methods experience dramatic runtime growth as the 3D grid resolution increases. AGMG consumes slightly less computing time than GMRES for both mesh sizes. Nevertheless, the performance gap becomes narrow in three dimensions, indicating that AGMG’s multigrid advantage is less pronounced for this 3D test setting compared to the 1D counterpart.

\begin{table}[htbp]
    \centering
    \begin{tabular}{|c|c|c|c|}
    \hline
        Method & $N_x=N_y=N_z$ & $N_t$& Total wall time (s) \\
        \hline
        GMRES &20 &40& 46.698018 \\ 
         &30 &40& 419.867325 \\ 
        \hline
        AGMG &20 &40&39.314921  \\
          &30 &40& 415.967842 \\
        \hline
    \end{tabular}
    \caption{The computational cost using GMRES, AGMG in the 3D case with different $N_x=N_y=N_z=20,30,N_t=40$, and $\alpha=0$. BDF1 with projection.}
    \label{tab:num-3}
\end{table}

\begin{figure}[htbp]
    \centering
    \subfloat[GMRES Iter]{\includegraphics[width=0.25\linewidth]{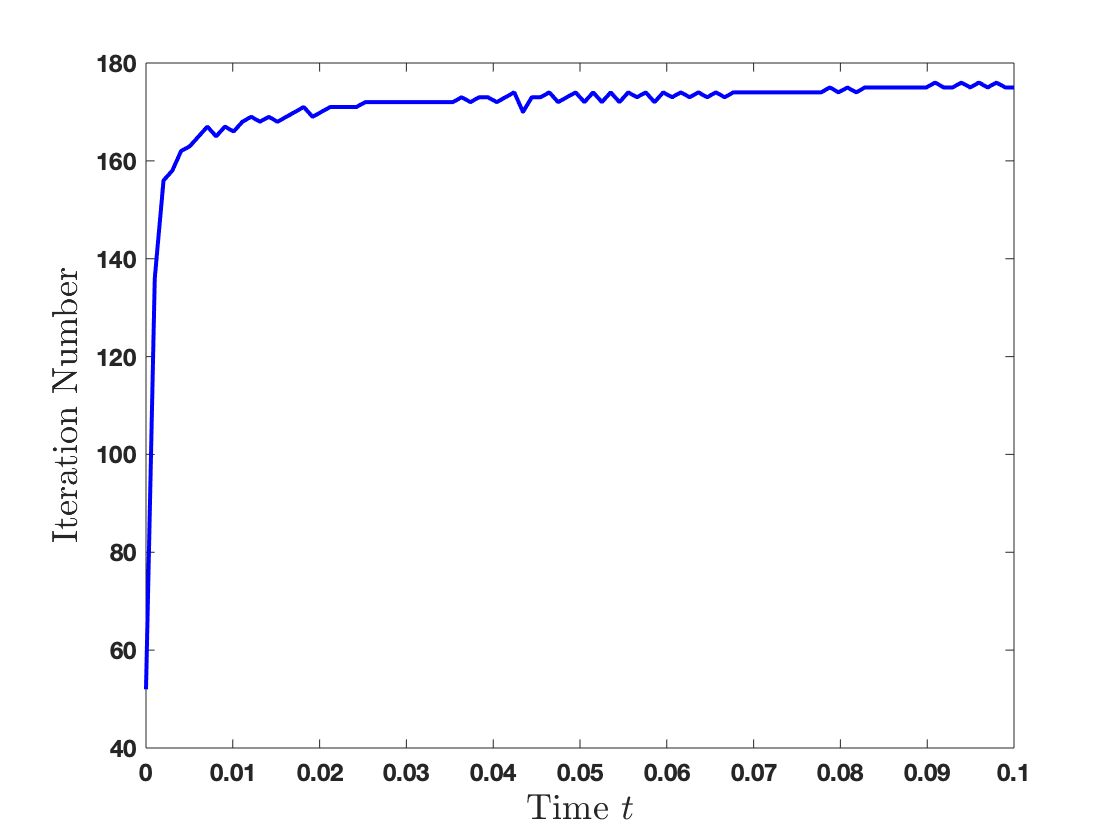}}
    \subfloat[GMRES relres]{\includegraphics[width=0.25\linewidth]{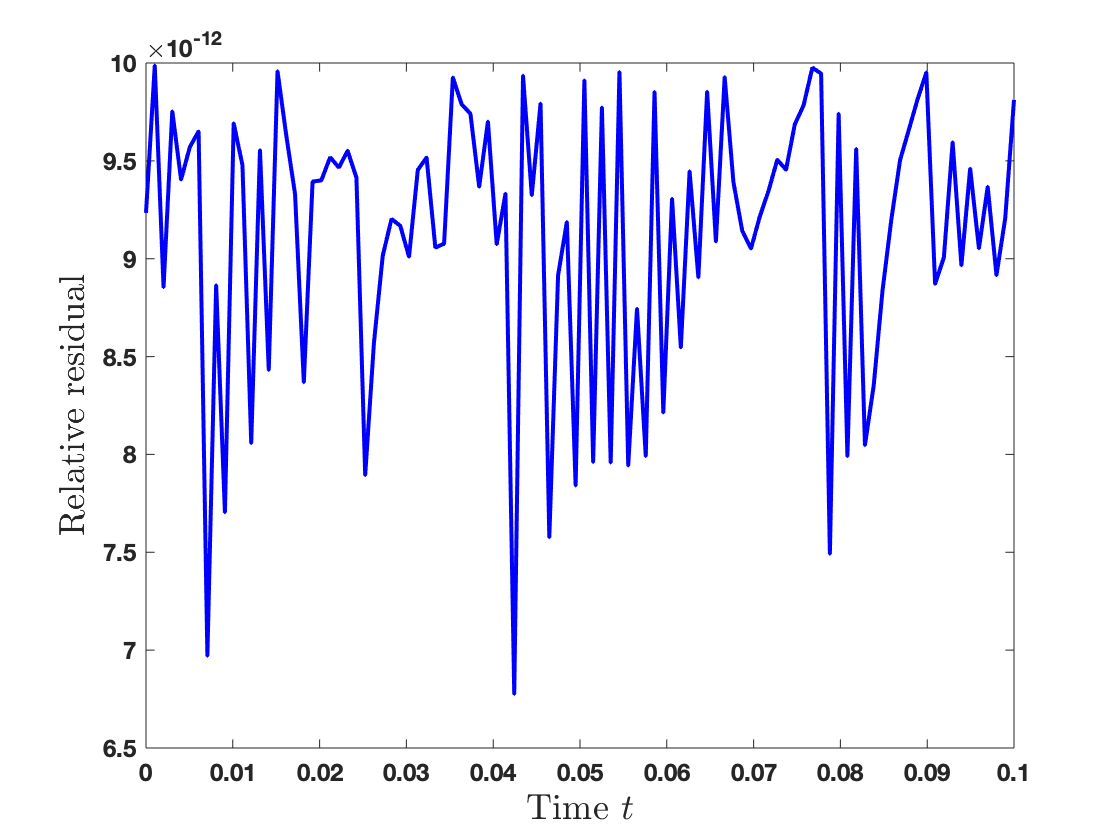}}
    % \hspace{0.1in}
    \subfloat[AGMG Iter]{\includegraphics[width=0.25\linewidth]{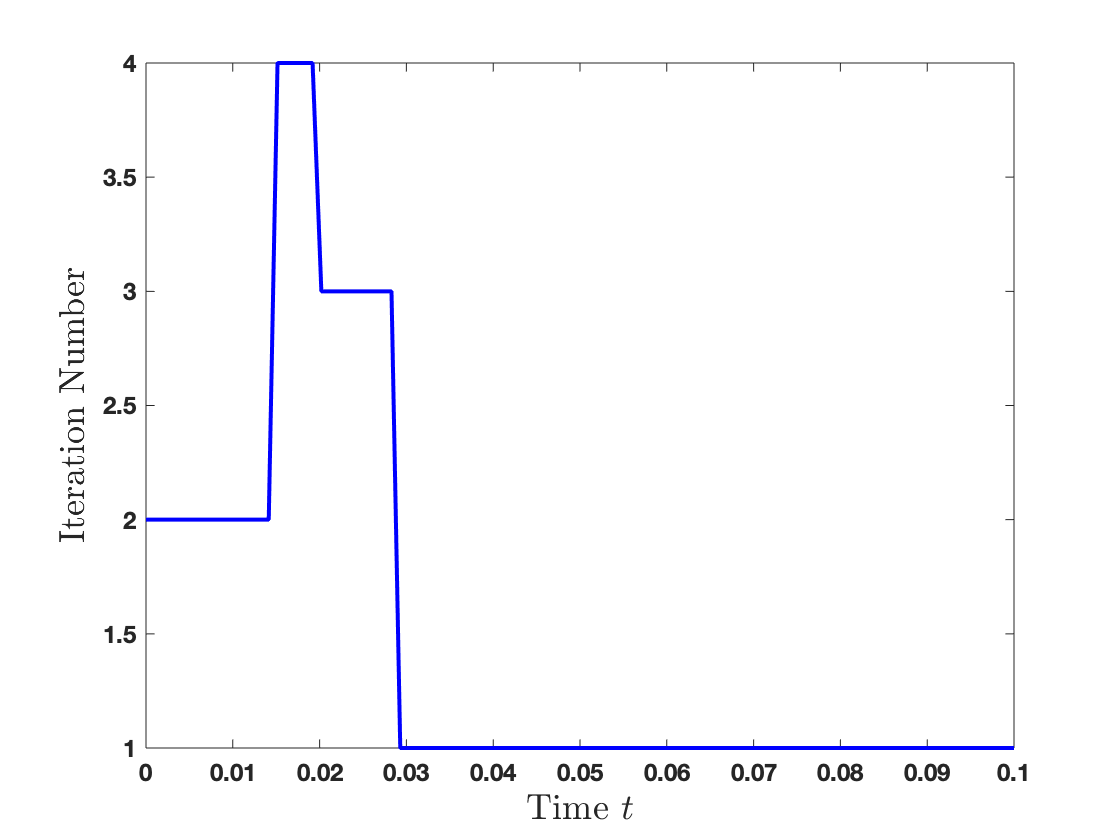}}
    \subfloat[AGMG relres]{\includegraphics[width=0.25\linewidth]{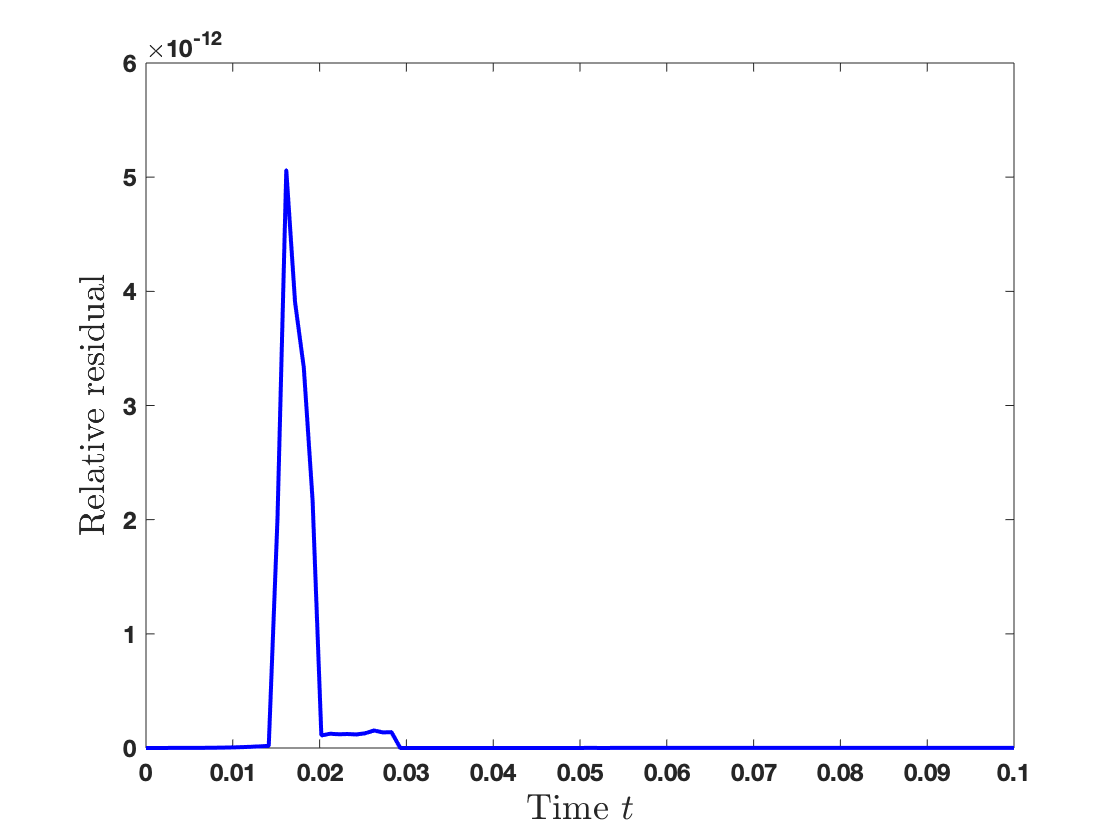}}
    \hspace{0.1in}
     \subfloat[GMRES Iter]{\includegraphics[width=0.25\linewidth]{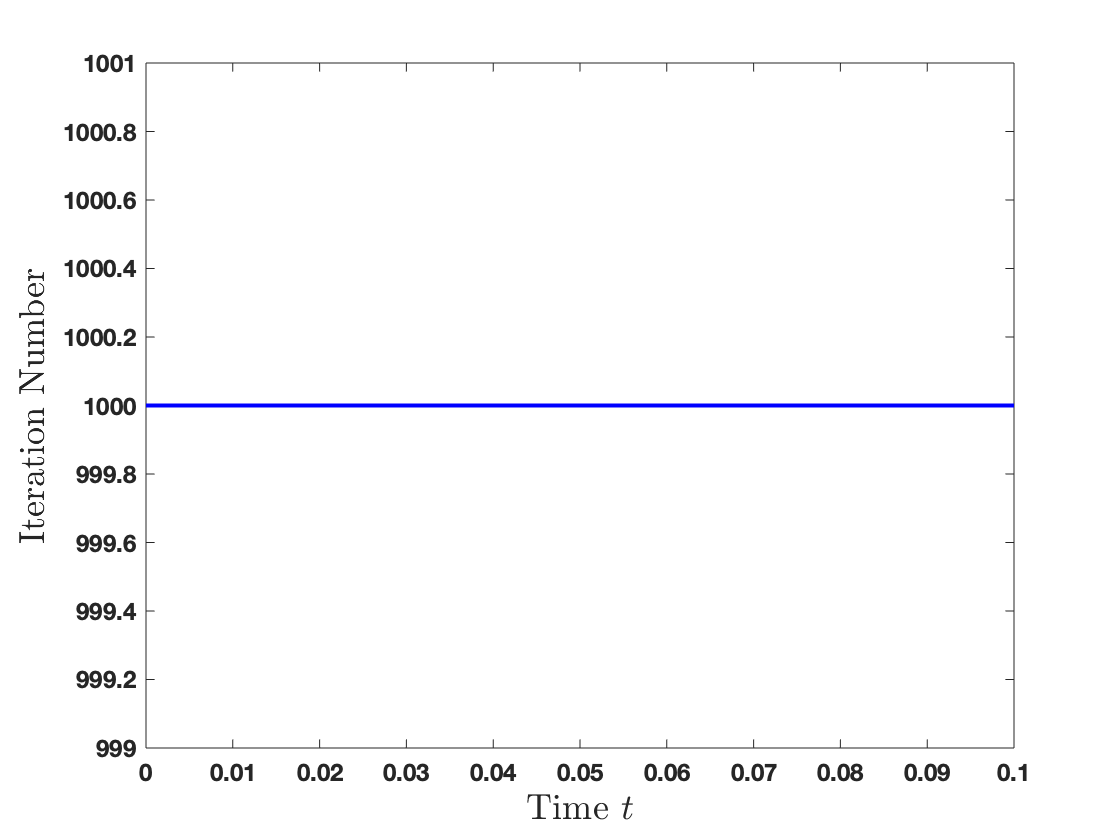}}
    \subfloat[GMRES relres]{\includegraphics[width=0.25\linewidth]{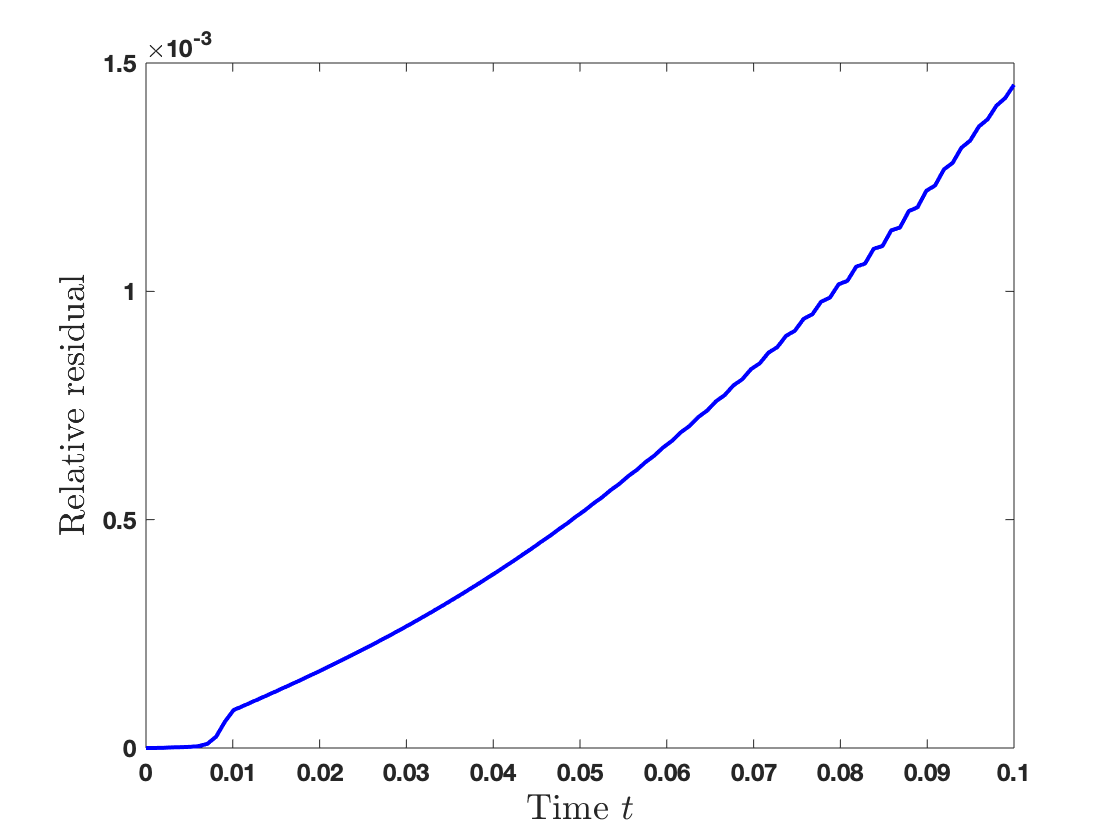}}
    % \hspace{0.1in}
    \subfloat[AGMG Iter]{\includegraphics[width=0.25\linewidth]{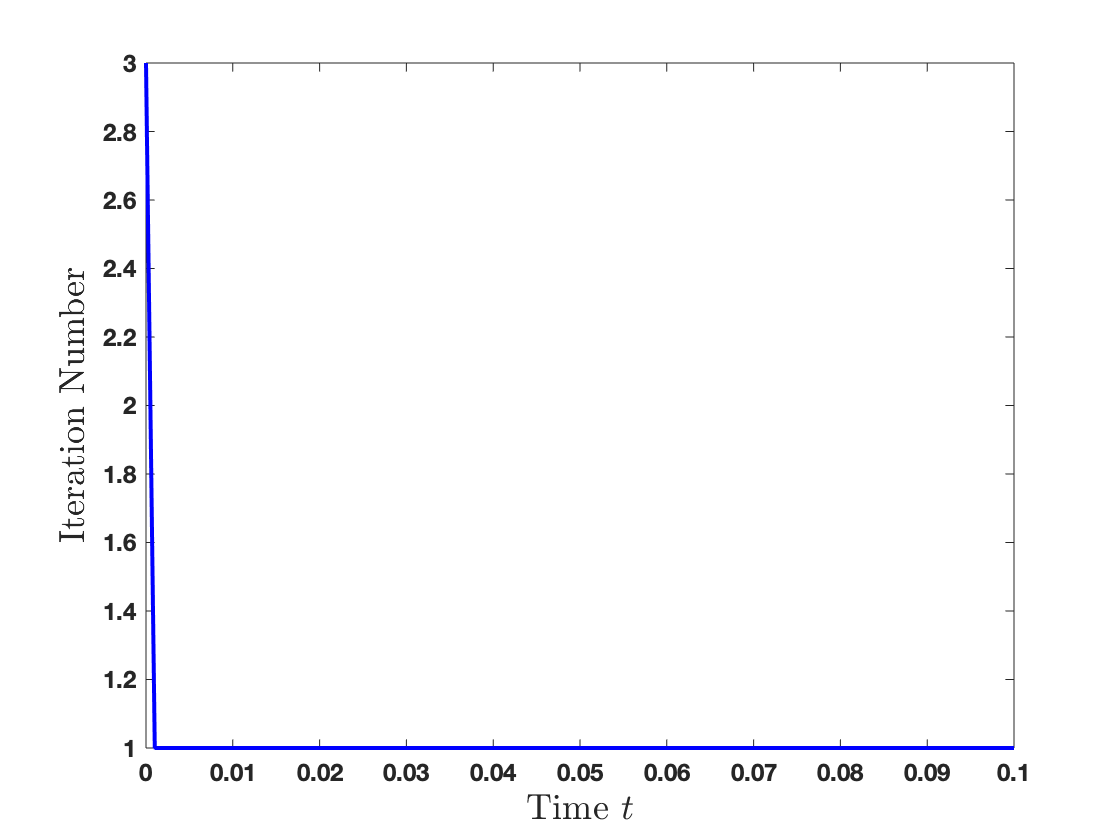}}
    \subfloat[AGMG relres]{\includegraphics[width=0.25\linewidth]{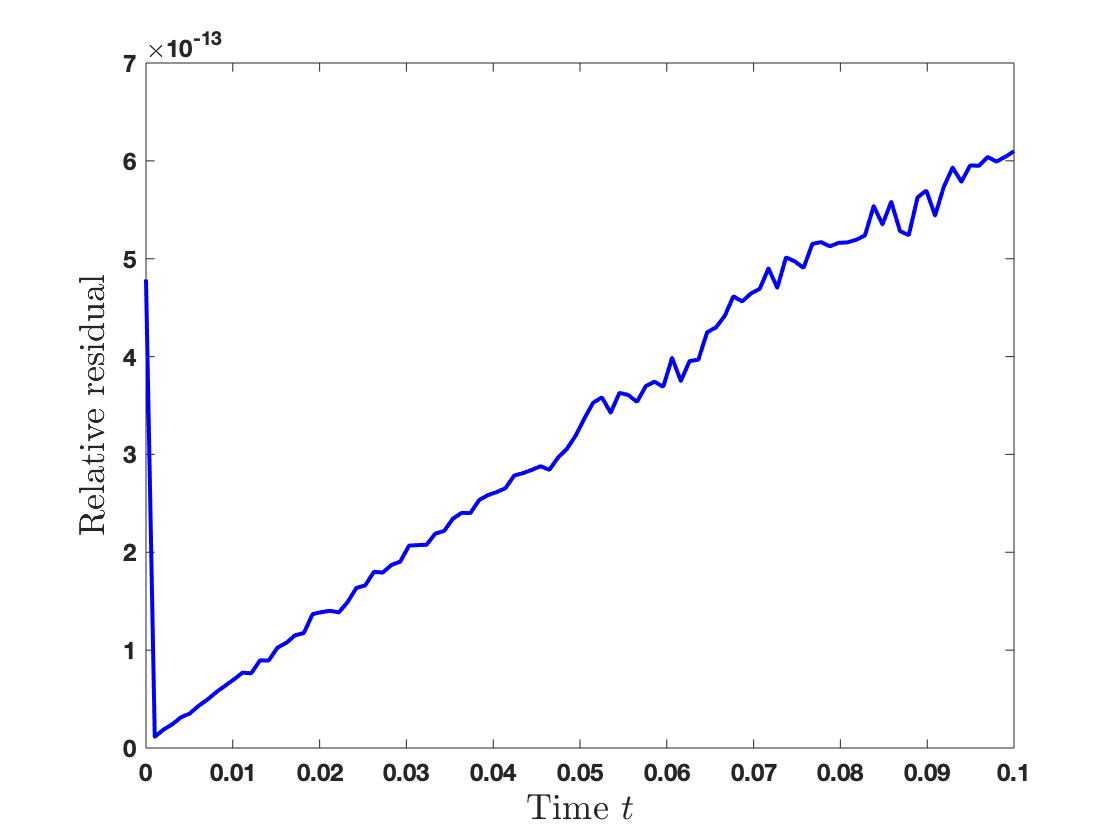}}
    \caption{Different linear solvers with $\alpha=0$, up to the final time $T=0.1$ for 1D case. Top rows: $N_x=100$, $N_t=100$ and Bottom rows: $N_x=1000$, $N_t=100$. BDF1 with projection.}
    \label{fig:BDF1-5}
\end{figure}

\Cref{fig:BDF1-5} illustrates iteration counts and relative residuals of GMRES and AGMG for the 1D BDF1 scheme with projection, under coarse (\(N_x=100\)) and fine (\(N_x=1000\)) spatial grids. GMRES requires far more iterations and exhibits strongly oscillatory or growing residuals over time, especially on the fine mesh. In contrast, AGMG maintains consistently low iteration numbers and well‑controlled relative residuals across both resolutions. This highlights AGMG’s superior robustness and convergence behaviour compared to standard GMRES for this problem.

% \begin{figure}[htbp]
%     \centering
%     \subfloat[gmres Iter]{\includegraphics[width=0.25\linewidth]{gmres_1D_Iter_v2_v1}}
%     \subfloat[gmres relres]{\includegraphics[width=0.25\linewidth]{gmres_1D_relres_v2_v1}}
%     \hspace{0.1in}
%     \subfloat[AGMG Iter]{\includegraphics[width=0.25\linewidth]{agmg_1D_iter_v2_v1}}
%     \subfloat[AGMG relres]{\includegraphics[width=0.25\linewidth]{agmg_1D_relres_v2_v1}}
%     \caption{Different linear solvers with $\alpha=0$, $N_x=1000$, $N_t=100$ up to the final time $T=0.1$ for 1D case (Here the flat 1000 iterations, means it exceeds the maximum of iteration number 1000).}
%     \label{fig:BDF1-6}
% \end{figure}

\begin{figure}[htbp]
    \centering
    \subfloat[GMRES Iter]{\includegraphics[width=0.25\linewidth]{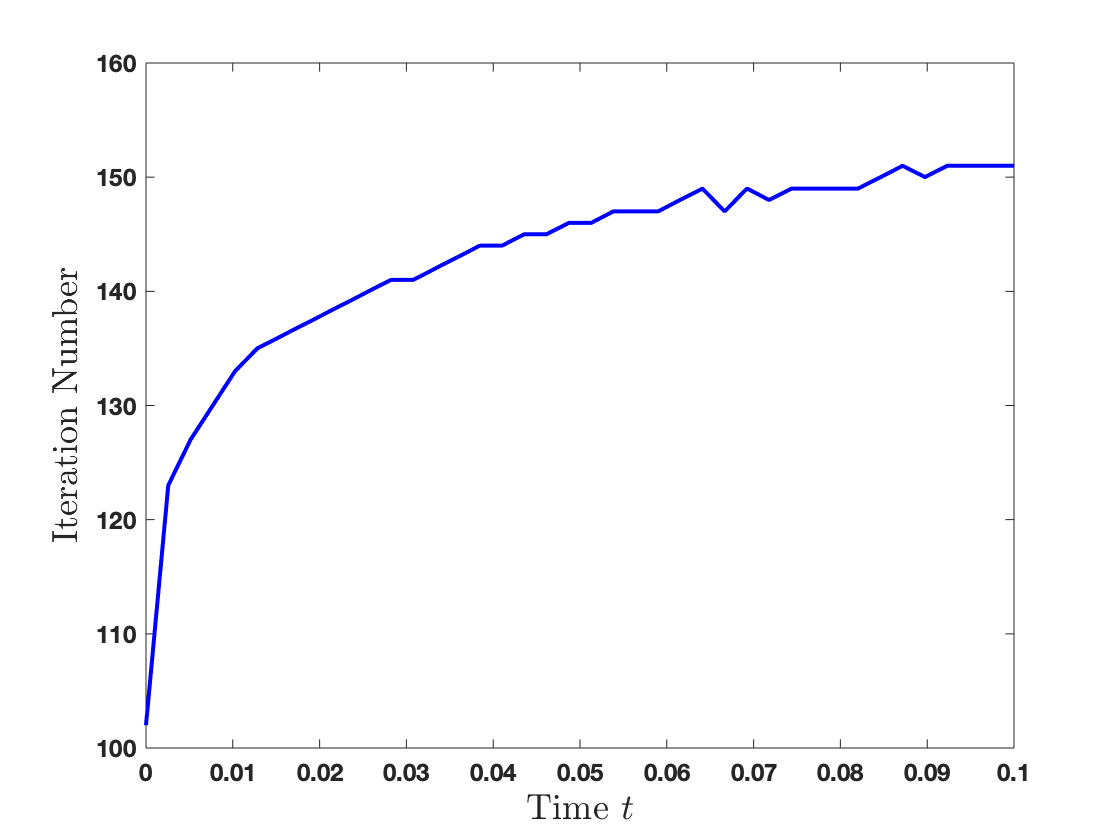}}
    \subfloat[GMRES relres]{\includegraphics[width=0.25\linewidth]{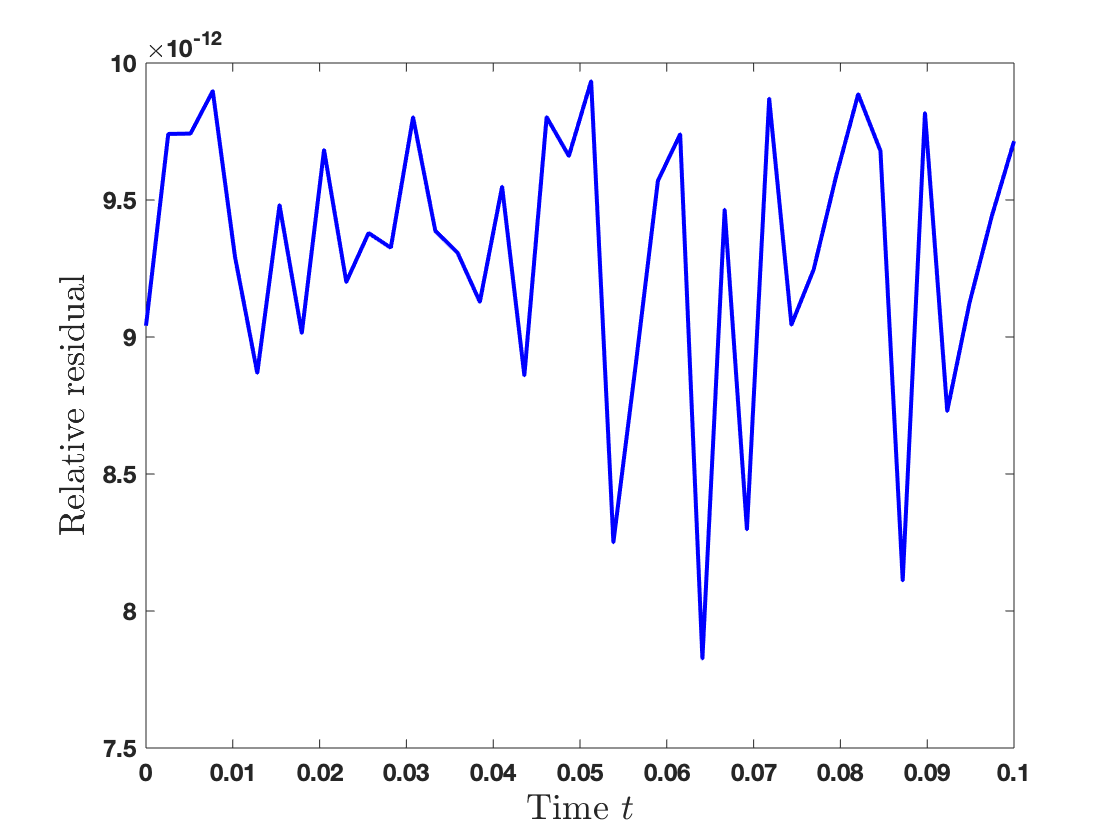}}
    % \hspace{0.1in}
    \subfloat[AGMG Iter]{\includegraphics[width=0.25\linewidth]{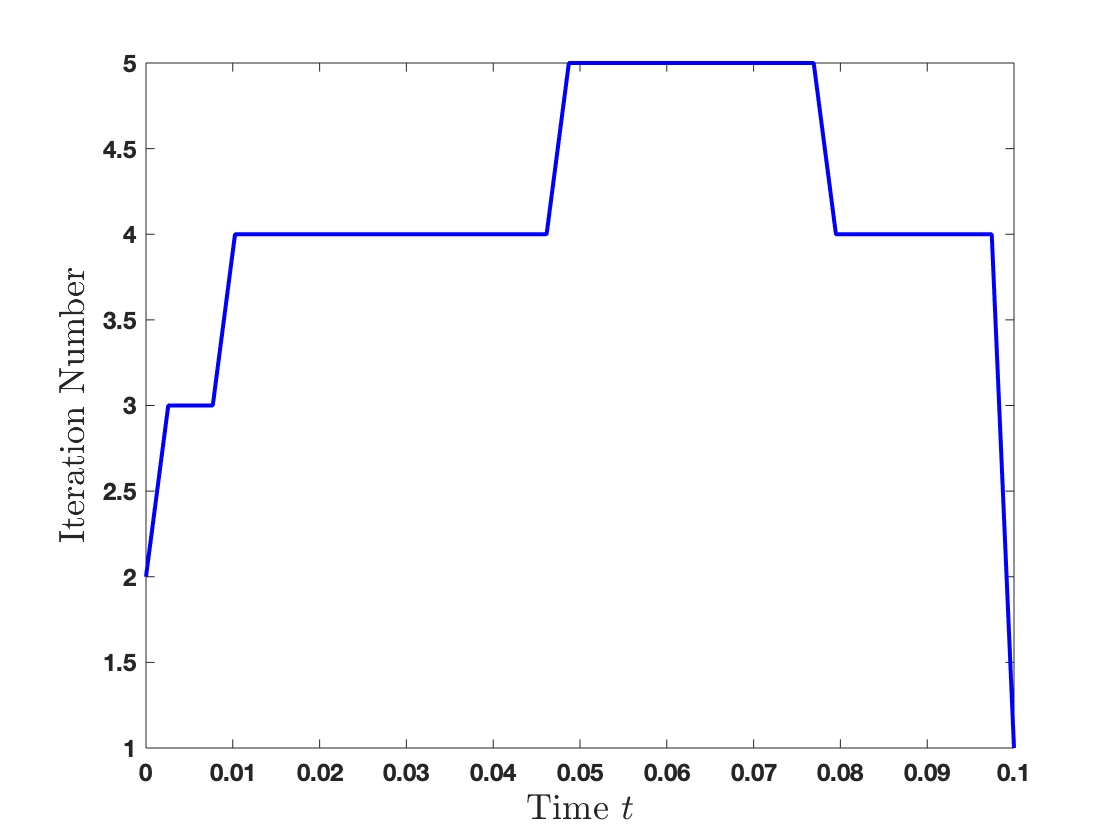}}
    \subfloat[AGMG relres]{\includegraphics[width=0.25\linewidth]{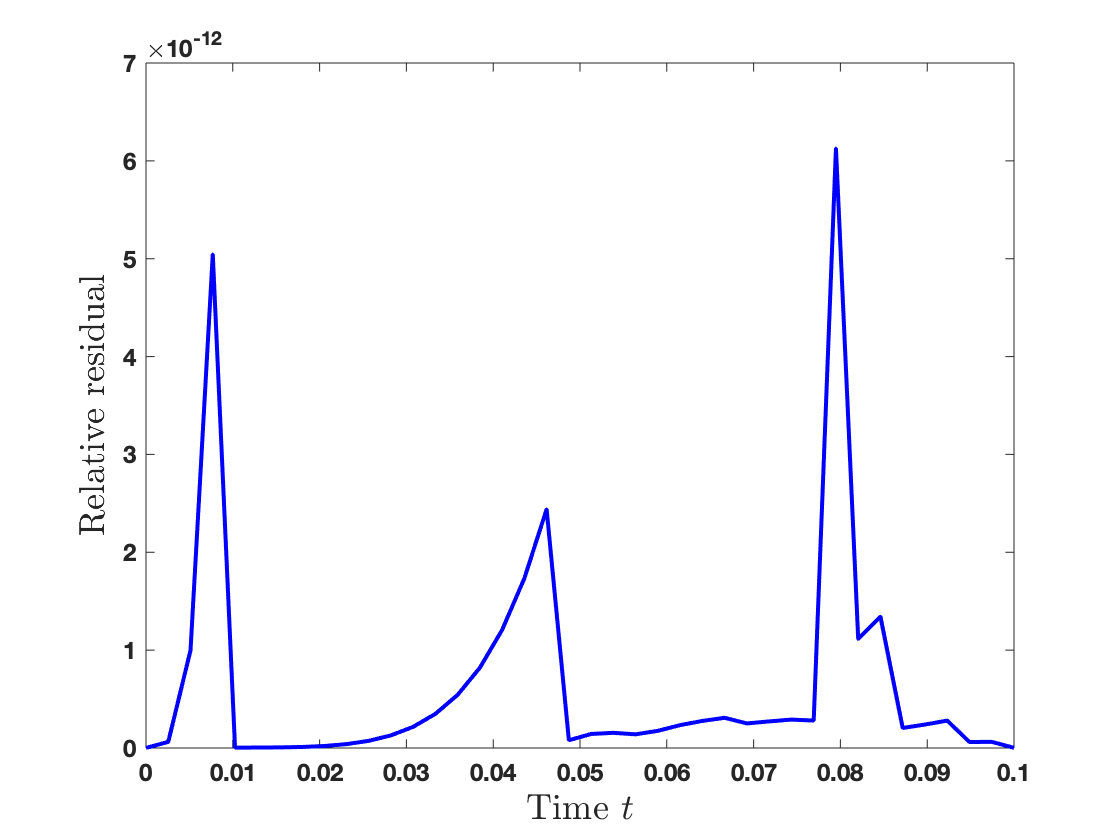}}
    \hspace{0.1in}
    \subfloat[GMRES Iter]{\includegraphics[width=0.25\linewidth]{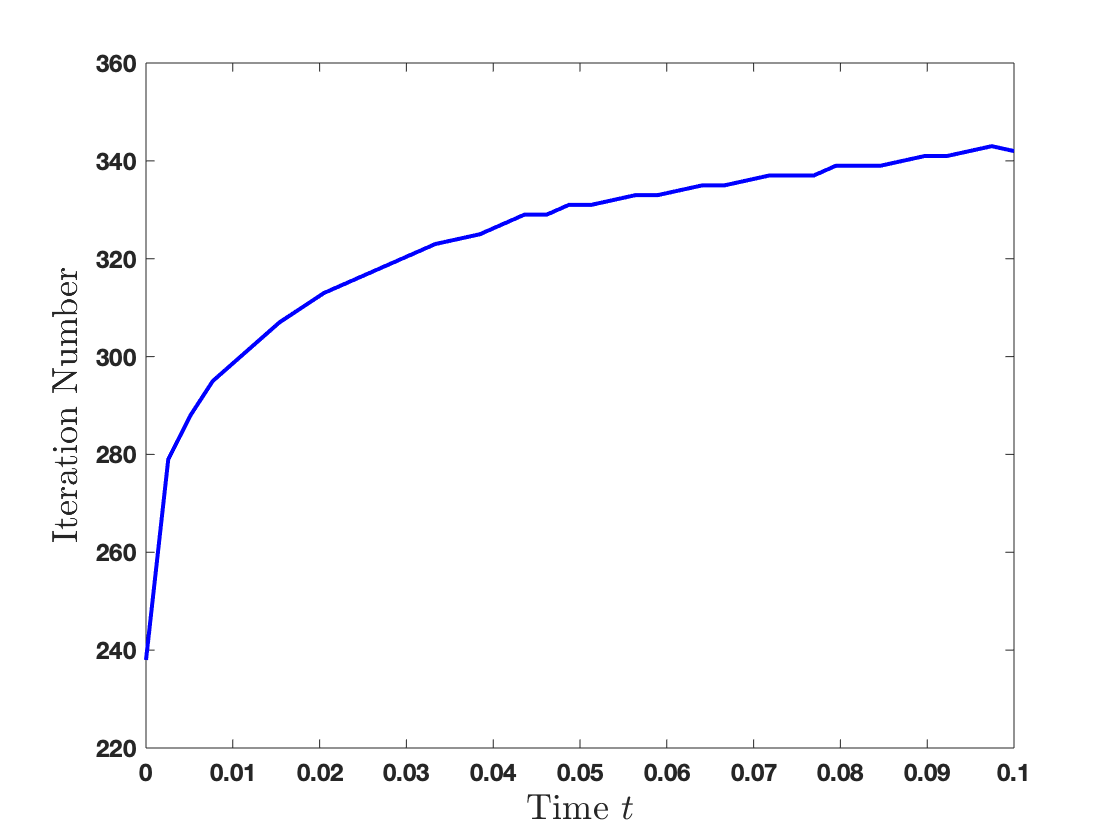}}
    \subfloat[GMRES relres]{\includegraphics[width=0.25\linewidth]{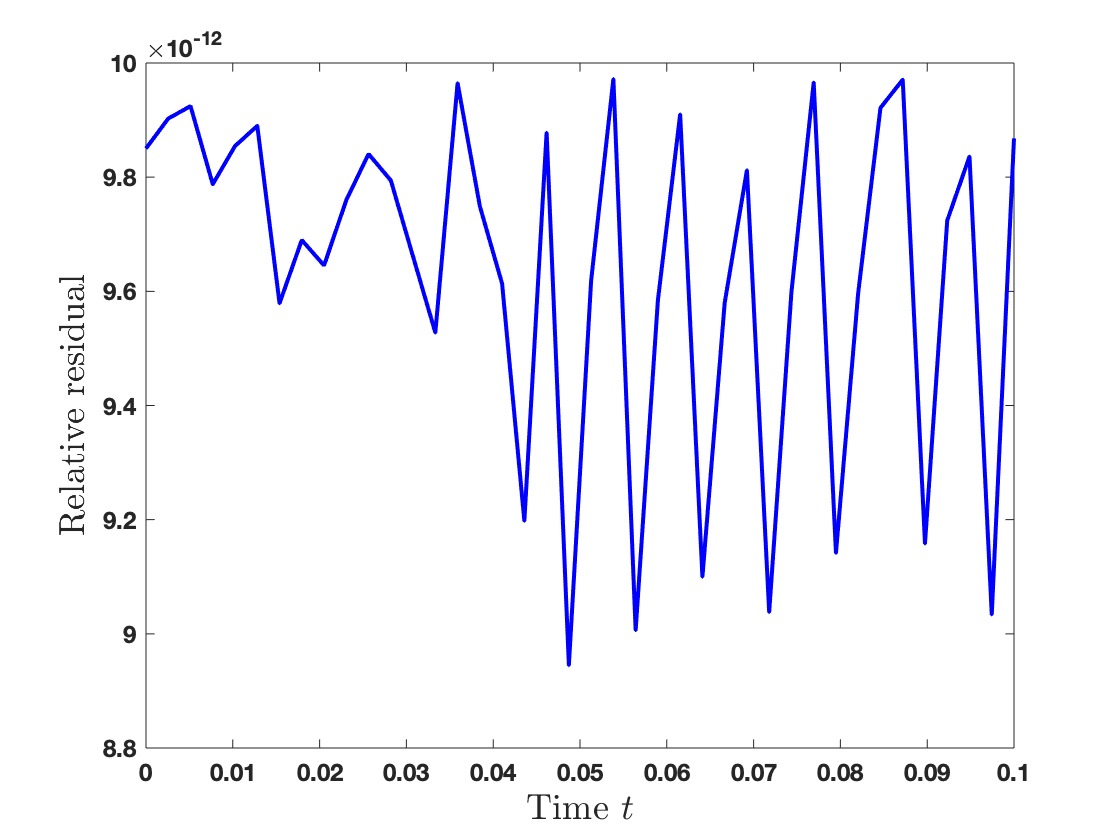}}
    % \hspace{0.1in}
    \subfloat[AGMG Iter]{\includegraphics[width=0.25\linewidth]{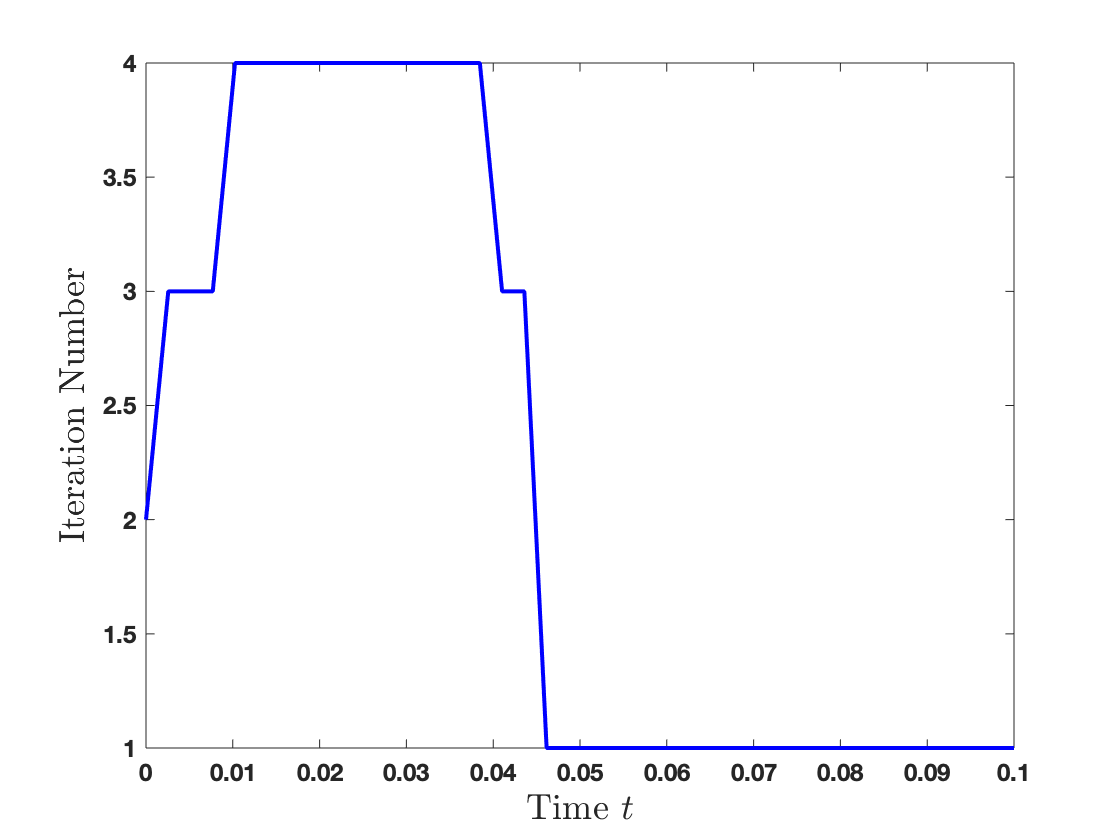}}
    \subfloat[AGMG relres]{\includegraphics[width=0.25\linewidth]{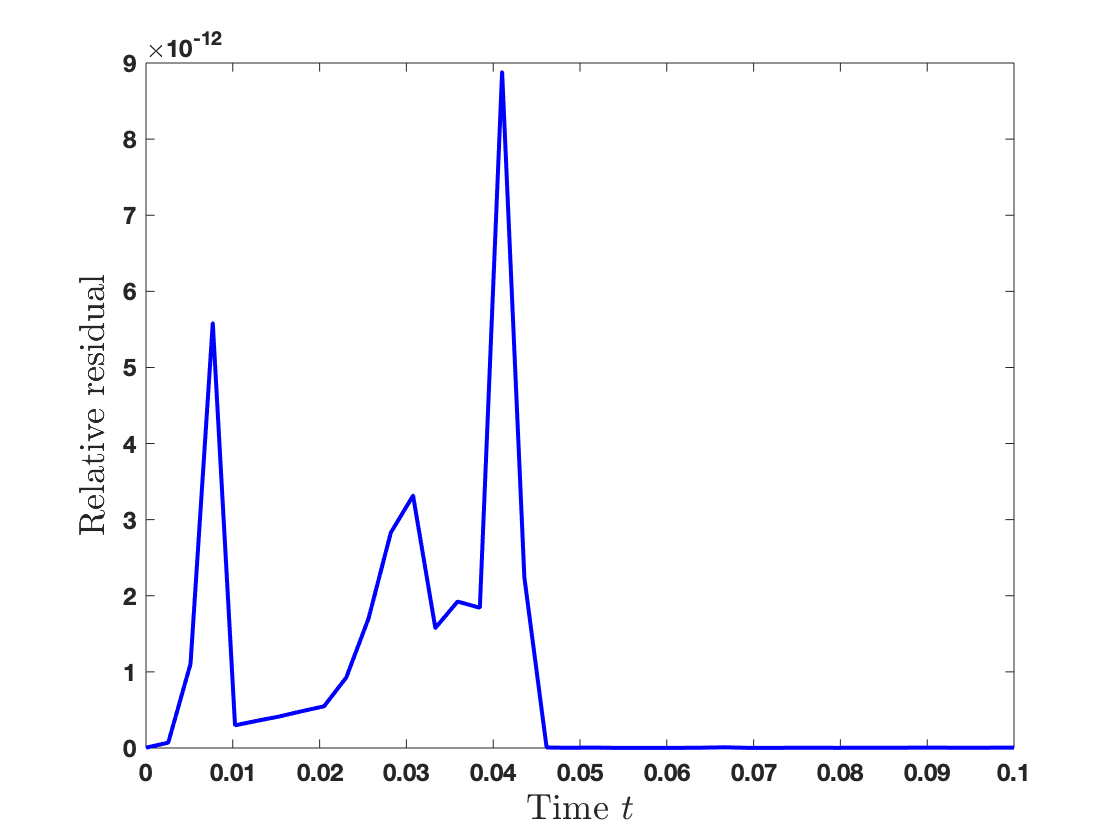}}
    \caption{Different linear solvers with $\alpha=0$, up to the final time $T=0.1$ for 3D case. Top rows: $N_x=N_y=N_z=20$, $N_t=40$ and Bottom rows: $N_x=N_y=N_z=30$, $N_t=40$. BDF1 with projection.}
    \label{fig:BDF1-7}
\end{figure}

\Cref{fig:BDF1-7} compares iteration counts and relative residuals of GMRES and AGMG for the 3‑D BDF1 scheme with projection on two uniform grids. GMRES requires substantially more iterations and shows persistent oscillatory relative residuals throughout temporal evolution. By contrast, AGMG maintains low iteration numbers and keeps residuals well‑controlled for both mesh resolutions. Even in three‑dimensional settings, AGMG exhibits more stable convergence behaviour than GMRES, though the performance gap becomes narrower than in the 1‑D counterpart.

\Cref{tab:num-4} summarizes wall‑clock time, average iteration counts and average relative residuals for GMRES and AGMG under various 3‑D mesh and time‑step configurations. GMRES often hits the maximum iteration limit of 1000, yielding extremely high runtime and degraded residual accuracy. In contrast, AGMG achieves far fewer average iterations, drastically shorter computation time and smaller residuals across most test cases. These metrics quantitatively verify AGMG’s prominent computational advantages over GMRES for this 3‑D BDF1‑projection framework.

% For other case,
\begin{table}[htbp]
    \centering
    \begin{tabular}{|c|c|c|c|c|c|c|c|}
    \hline
        Method & $N_x$ &$N_y$&$N_z$ & $N_t$& Total wall time (s)&average iter& average relres \\
        \hline
        GMRES &50&50&4 &40&700.675728 &5.636750000000000e+02&9.834080272744706e-12  \\ 
         & 100&50&4&10&704.968226  &1000& 4.484803064091582e-04 \\ 
         &100&4&4&10&123.035312&1000&1.555713012324725e-04\\
         &50&4&4&100&5.911617&56.189999999999998&9.166091267486675e-12\\
         &10&10&10&100&7.513409&12.320000000000000&6.676551578819269e-12\\
        \hline
        AGMG &50&50&4 &40&39.021992  &1.675000000000000&4.704135380268438e-13 \\
          & 100&50&4 &10& 31.563383&2.100000000000000&1.031574184679247e-14 \\
           &100&4&4&10&1.844813&1.100000000000000&1.483813562562323e-13\\
           &50&4&4&100&5.163407&3.640000000000000&3.195664477316403e-13\\
           &10&10&10&100&9.123604&8.789999999999999&3.442629545316088e-12\\
        \hline
    \end{tabular}
    \caption{The computational cost, average iteration over time, average relative residual over time using GMRES, AGMG in the 3D case with different $N_x,N_y,N_z,N_t$, and $\alpha=0$. In the runtime, we set the $maxit=1000$ and $tol=10^{-11}$. BDF1 with projection.}
    \label{tab:num-4}
\end{table}

\begin{figure}[htbp]
    \centering
    \subfloat[arrow, GMRES]{\includegraphics[width=0.3\linewidth]{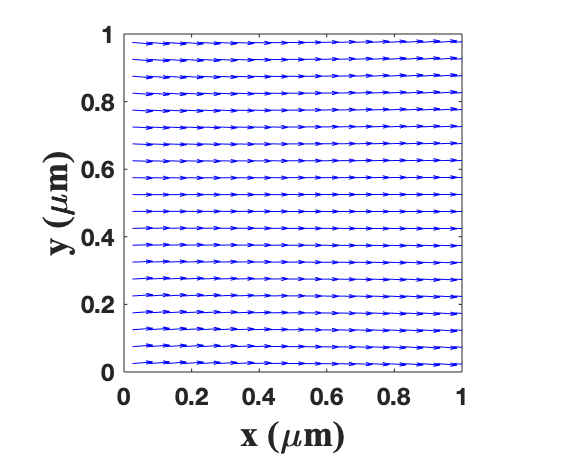}}
     \subfloat[color, GMRES]{\includegraphics[width=0.3\linewidth]{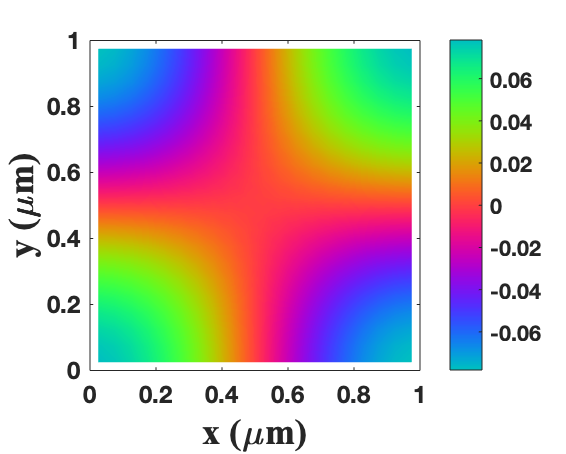}}
         \hspace{0.1in}
    \subfloat[arrow, AGMG]{\includegraphics[width=0.3\linewidth]{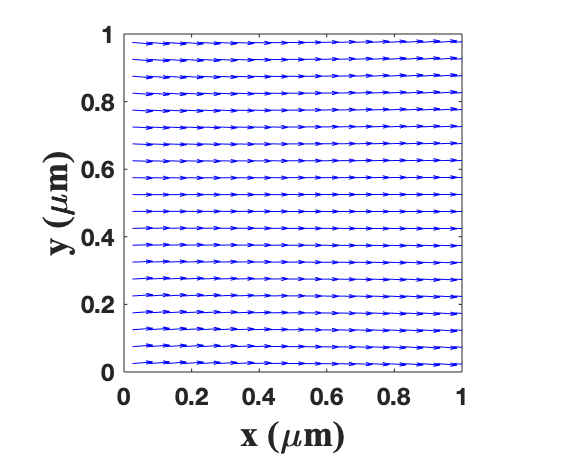}}
    \subfloat[color, AGMG]{\includegraphics[width=0.3\linewidth]{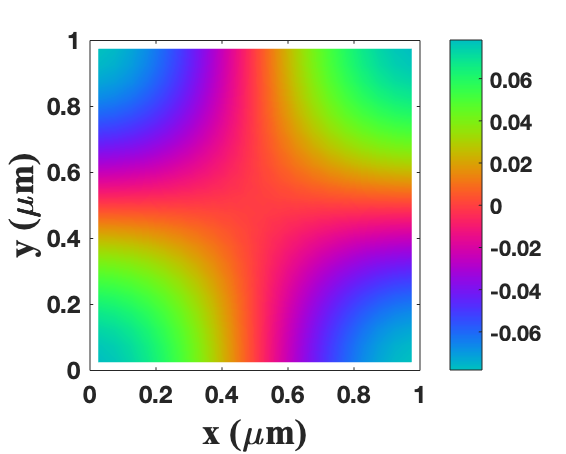}}
    \caption{The solution profiles in 3D case with $N_x=20,N_y=20,N_z=20$ and $N_t=40$ up to the final time $T=0.1$ with $\alpha=0$. Using GMRES and AGMG. The initial condition $\m_0=[\cos(\cos(\pi x)\cos(\pi y)\cos(\pi z))\sin(0),\sin(\cos(\pi x)\cos(\pi y)\cos(\pi z))\sin(0),\cos(0)]^T$. The profile take the slice at $z=1/2$. BDF1 with projection.}
    \label{fig:BDF-8}
\end{figure}

\Cref{fig:BDF-8} visualizes the magnetization solution at the slice \(z=1/2\) for the 3‑D BDF1‑projection simulation solved by GMRES and AGMG. Both arrow‑vector plots and colour‑contour plots demonstrate that the two solvers produce visually indistinguishable spatial magnetization profiles. Despite their large differences in iteration counts and computational overhead, AGMG yields numerical solutions fully consistent with GMRES. This confirms AGMG can retain solution accuracy while delivering superior computational efficiency for the micro‑magnetism model.

\Cref{tab:num-5} compares wall‑clock time, average iteration counts and average relative residuals between GMRES and AGMG for the two‑dimensional BDF1‑projection scheme across various spatial grids. GMRES consumes considerably longer runtime and requires many more iterations, especially for fine meshes. By contrast, AGMG yields drastically reduced iteration numbers and computational time while maintaining comparable small relative residuals. These quantitative results demonstrate AGMG’s prominent computational advantage over GMRES for the 2‑D micro‑magnetism problem.

\begin{table}[htbp]
    \centering
    \begin{tabular}{|c|c|c|c|c|c|c|}
    \hline
        Method & $N_x$ &$N_y$ & $N_t$& Total wall time (s)&average iter& average relres \\
        \hline
        GMRES &50&50&100&58.071110&1.367900000000000e+02&9.635084864234602e-12  \\ 
         &100&100&100&1434.450792&5.456700000000000e+02&9.910529473107728e-12 \\ 
         &50&20&100&7.935363&78.400000000000006&9.457445177659933e-12\\
         &20&50&100&7.057161&78.420000000000002&9.439176694610438e-12\\
         &80&80&100&432.312840&3.493300000000000e+02&9.868579608651962e-12\\
        \hline
        AGMG &50&50&100&15.665064&2.590000000000000&1.283841524630060e-12 \\
          &100&100&100&74.747998&1.330000000000000&4.785253903876194e-13 \\
           &50&20&100&5.735392&3.990000000000000&1.590648630286451e-12\\
           &20&50&100&4.546096&4.140000000000000&1.542554274184092e-12\\
           &80&80&100&29.024841&1.610000000000000&4.420535701386983e-13\\
        \hline
    \end{tabular}
    \caption{The computational cost, average iteration over time, average relative residual over time using GMRES, AGMG in the 2D case with different $N_x,N_y,N_t$, and $\alpha=0$. In the runtime, we set the $maxit=1000$ and $tol=10^{-11}$. BDF1 with projection.}
    \label{tab:num-5}
\end{table}

\begin{figure}[htbp]
    \centering
    \subfloat[arrow, GMRES]{\includegraphics[width=0.3\linewidth]{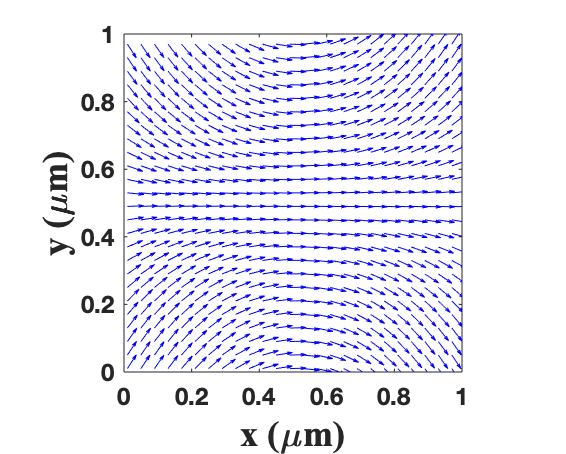}}
     \subfloat[color, GMRES]{\includegraphics[width=0.3\linewidth]{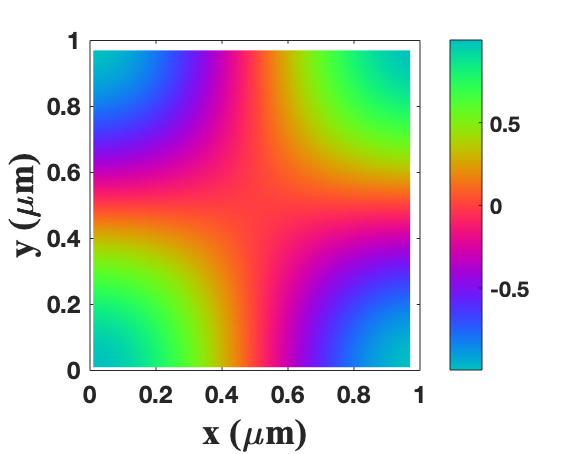}}
         \hspace{0.1in}
    \subfloat[arrow, AGMG]{\includegraphics[width=0.3\linewidth]{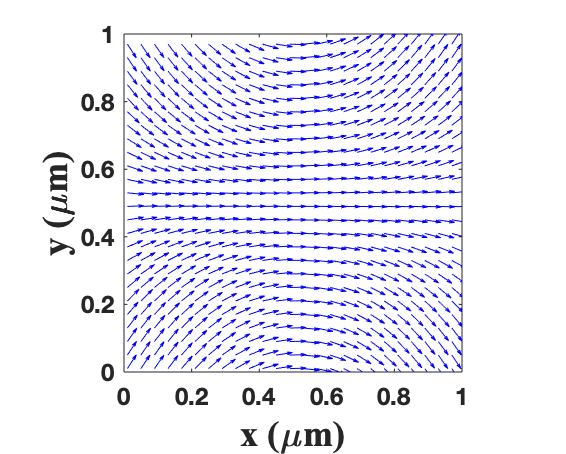}}
    \subfloat[color, AGMG]{\includegraphics[width=0.3\linewidth]{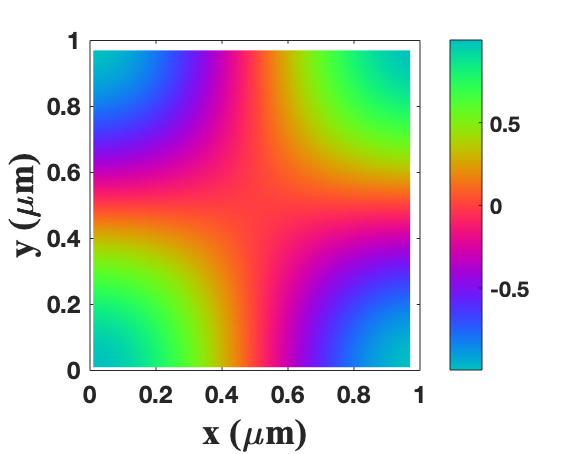}}
    \caption{The solution profiles in the 2D case with $N_x=50,N_y=50$ and $N_t=100$ up to the final time $T=0.1$ with $\alpha=0$. Using GMRES and AGMG. The initial condition $\m_0=[\cos(\cos(\pi x)\cos(\pi y))\sin(0),\sin(\cos(\pi x)\cos(\pi y))\sin(0),\cos(0)]^T$. BDF1 with projection.}
    \label{fig:BDF-9}
\end{figure}

\Cref{fig:BDF-9} presents arrow‑vector and colour‑contour plots of magnetization for the two‑dimensional BDF1‑projection simulation solved by GMRES and AGMG. Visually, the spatial magnetization profiles obtained from both solvers show no observable difference. Although AGMG achieves much lower iteration counts and computational cost than GMRES, it reproduces identical physical solutions. This observation validates that AGMG preserves numerical accuracy while delivering superior computational efficiency for the 2D micromagnetism model.

\begin{figure}[htbp]
    \centering
    \subfloat[arrow, GMRES]{\includegraphics[width=0.3\linewidth]{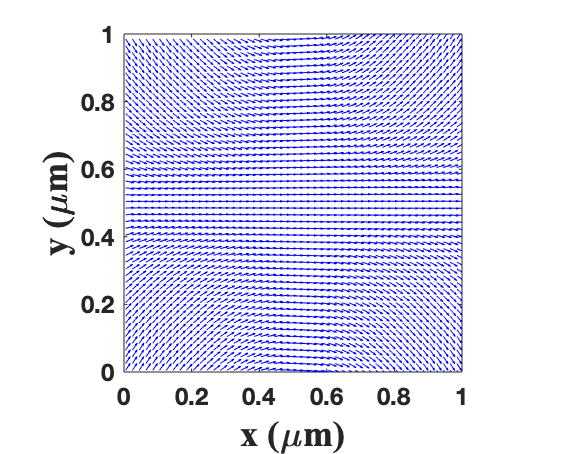}}
     \subfloat[color, GMRES]{\includegraphics[width=0.3\linewidth]{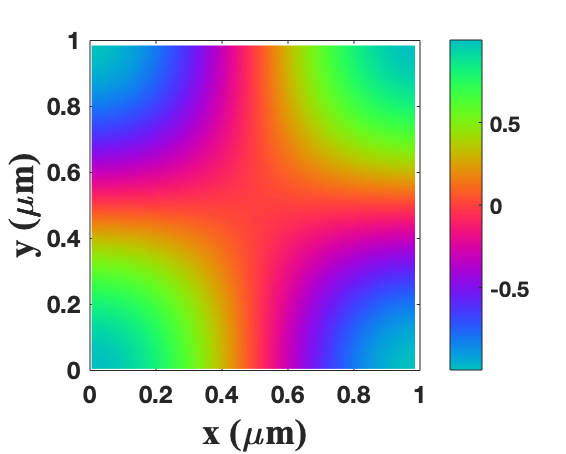}}
         \hspace{0.1in}
    \subfloat[arrow, AGMG]{\includegraphics[width=0.3\linewidth]{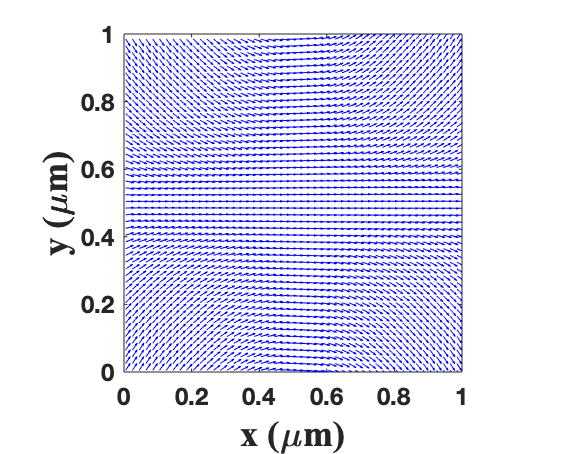}}
    \subfloat[color, AGMG]{\includegraphics[width=0.3\linewidth]{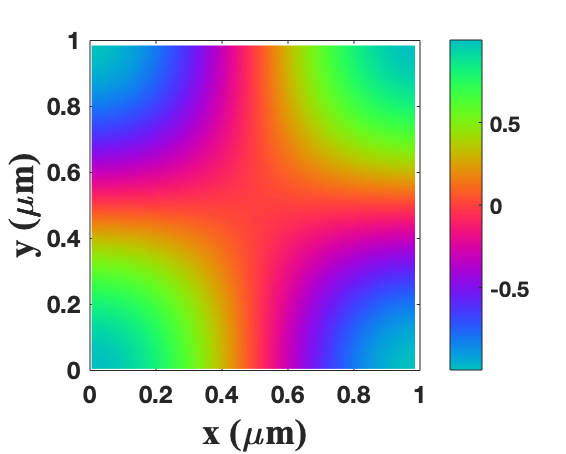}}
    \caption{The solution profiles in the 2D case with $N_x=100,N_y=100$ and $N_t=100$ up to the final time $T=0.1$ with $\alpha=0$. Using GMRES and AGMG. The initial condition $\m_0=[\cos(\cos(\pi x)\cos(\pi y))\sin(0),\sin(\cos(\pi x)\cos(\pi y))\sin(0),\cos(0)]^T$. BDF1 with projection.}
    \label{fig:BDF-10}
\end{figure}

\Cref{fig:BDF-10} visualizes the 2‑D magnetization profiles on a fine grid \(N_x=N_y=100\) for the BDF1‑projection scheme solved by GMRES and AGMG. Both arrow‑vector diagrams and colour‑contour maps demonstrate that the two solvers produce visually identical magnetization distributions. Even under dense spatial discretization, AGMG yields physically consistent solutions as GMRES. Combined with previous numerical metrics, this confirms AGMG maintains solution fidelity while offering greatly improved computational efficiency for fine‑grid micromagnetic simulations.

\subsection{Generalization to the damping case}
In this section, we test the GMRES and AGMG solver in 1D, 2D, 3D case with damping $\alpha=0.01$. \Cref{tab:num-6} presents computational metrics for GMRES and AGMG under the one‑dimensional BDF1 scheme without projection with \(\alpha=0.01\). For the fine grid \(N_x=1000\), GMRES hits the maximum iteration limit of 1000, yielding high wall‑clock time and deteriorated residual accuracy. In contrast, AGMG maintains very low average iteration counts, substantially shorter runtime and small residuals. These results demonstrate AGMG’s robust computational superiority even for non‑projected 1‑D micromagnetic problems.

\begin{table}[htbp]
    \centering
    \begin{tabular}{|c|c|c|c|c|c|}
    \hline
        Method & $N_x$ & $N_t$& Total wall time (s) &average iter &average relres\\
        \hline
        GMRES &100 &100& 2.600267 &1.695700000000000e+02&9.143869105167878e-12 \\ 
         &1000 &100& 322.504298 &1000&5.707928689918663e-04 \\ 
        \hline
        AGMG &100 &100& 0.466697 &1.810000000000000&5.009128833400892e-13 \\
          &1000 &100& 3.748856&1.010000000000000&2.409319292241509e-13 \\
        \hline
    \end{tabular}
    \caption{The computational cost, average iteration over time, average relative residual over time using GMRES, AGMG in the 1D case with different $N_x,N_t$, and  $\alpha=0.01$. BDF1 without projection.}
    \label{tab:num-6}
\end{table}

\Cref{tab:num-7} compares GMRES and AGMG for the 2‑D projected BDF1 scheme with damping parameter \(\alpha=0.01\) across various spatial grids. GMRES exhibits large average iteration counts and considerably longer wall‑clock time, whose computational overhead grows markedly for fine meshes. By contrast, AGMG retains nearly constant low iteration numbers, achieves much shorter runtime, and maintains small relative residuals. These quantitative results verify that AGMG delivers superior computational efficiency without sacrificing accuracy for damped two‑dimensional micromagnetic simulations.

\begin{table}[htbp]
    \centering
    \begin{tabular}{|c|c|c|c|c|c|c|}
    \hline
        Method & $N_x$ &$N_y$ & $N_t$& Total wall time (s)&average iter& average relres \\
        \hline
        GMRES &50&50&100&45.079057&1.339400000000000e+02&9.657318474538825e-12  \\ 
         % &100&100&100&&& \\ 
         &50&20&100&6.163240&77.170000000000002&9.370022527455704e-12\\
         &20&50&100&6.324561&77.159999999999997&9.387602394501081e-12\\
         &80&80&100&395.501187&3.397800000000000e+02&9.854547478942513e-12\\
        \hline
        AGMG &50&50&100&10.444269&3.130000000000000&1.156627660977494e-12 \\
          % &100&100&100&&& \\
           &50&20&100&5.772690&4.420000000000000&5.473021691057962e-13\\
           &20&50&100&4.378038&4.430000000000000&6.367888642123739e-13\\
           &80&80&100&31.877171&2.500000000000000&4.890880713350495e-13\\
        \hline
    \end{tabular}
    \caption{The computational cost, average iteration over time, average relative residual over time using GMRES, AGMG in the 2D case with different $N_x,N_y,N_t$, and $\alpha=0.01$ up to the final time $T=0.1$. In the runtime, we set the $maxit=1000$ and $tol=10^{-11}$. BDF1 with projection.}
    \label{tab:num-7}
\end{table}

\Cref{tab:num-8} compares GMRES and AGMG for 3‑D projected BDF1 micromagnetic simulations with \(\alpha=0.01\) under diverse mesh configurations. GMRES frequently reaches the maximum iteration limit of 1000, leading to excessive wall‑clock time and degraded residual accuracy for large‑size problems. In contrast, AGMG maintains consistently small iteration counts, drastically reduces computational runtime, and preserves tight relative residuals. This demonstrates AGMG’s outstanding robustness and efficiency for challenging three‑dimensional micromagnetic computations.

\begin{table}[htbp]
    \centering
    \begin{tabular}{|c|c|c|c|c|c|c|c|}
    \hline
        Method & $N_x$ &$N_y$&$N_z$ & $N_t$& Total wall time (s)&average iter& average relres \\
        \hline
        GMRES &50&50&4 &40&620.723074&5.477000000000000e+02&9.858476506896158e-12  \\ 
         & 100&50&4&10&893.643937&1000&3.965404773605235e-04\\ 
         &100&4&4&10&157.278015&1000&1.360703271204140e-04\\
         &50&4&4&100&10.453303&55.450000000000003&9.194578163656896e-12\\
         &10&10&10&100&8.039507&12.310000000000000&6.570329286522540e-12\\
        \hline
        AGMG &50&50&4 &40&48.436074&2.825000000000000&4.924872746306352e-13 \\
          & 100&50&4 &10& 27.753632&1.100000000000000&1.996216469484103e-13 \\
           &100&4&4&10&2.543446&1.100000000000000&1.195257611992792e-13\\
           &50&4&4&100&5.468379&4.430000000000000&1.094427125751366e-12\\
           &10&10&10&100&9.438602&8.760000000000000&3.246262981192344e-12\\
        \hline
    \end{tabular}
    \caption{The computational cost, average iteration over time, average relative residual over time using GMRES, AGMG in the 3D case with different $N_x,N_y,N_z,N_t$, and $\alpha=0.01$ up to the final time $T=0.1$. In the runtime, we set the $maxit=1000$ and $tol=10^{-11}$. BDF1 with projection.}
    \label{tab:num-8}
\end{table}

\Cref{tab:num-9} reports numerical comparisons between GMRES and AGMG for three‑dimensional BDF1 simulations without projection, with \(\alpha=0.01\) on various spatial‑temporal grids. GMRES often saturates at the maximum iteration limit of 1000, yielding heavy computational overhead and poor residual accuracy for large‑scale 3D configurations. By contrast, AGMG retains small iteration numbers, achieves drastically reduced wall‑clock time and maintains satisfactory residuals. These outcomes highlight AGMG’s compelling efficiency and robustness for non‑projected 3D micromagnetic problems.

\begin{table}[htbp]
    \centering
    \begin{tabular}{|c|c|c|c|c|c|c|c|}
    \hline
        Method & $N_x$ &$N_y$&$N_z$ & $N_t$& Total wall time (s)&average iter& average relres \\
        \hline
        GMRES &50&50&4 &40&522.956746&5.477250000000000e+02&9.855887219964113e-12  \\ 
         & 100&50&4&10&687.980963&1000&3.961776568537530e-04\\ 
         &100&4&4&10&128.572499&1000&1.359310738073514e-04\\
         &50&4&4&100&5.752891&55.479999999999997&9.136203668205477e-12\\
         &10&10&10&100&8.326046&12.310000000000000&6.567601241117217e-12\\
        \hline
        AGMG &50&50&4 &40&32.315438&2.825000000000000&4.914594767554422e-13 \\
          & 100&50&4 &10&20.012595 &1.100000000000000&1.996990180282928e-13 \\
           &100&4&4&10&1.648144&1.100000000000000&1.194501100351557e-13\\
           &50&4&4&100&5.097735&4.430000000000000&1.094391404823389e-12\\
           &10&10&10&100&9.354798&8.760000000000000&3.244218176745128e-12\\
        \hline
    \end{tabular}
    \caption{The computational cost, average iteration over time, average relative residual over time using GMRES, AGMG in the 3D case with different $N_x,N_y,N_z, N_t$, and $\alpha=0.01$ up to the final time $T=0.1$. In the runtime, we set the $maxit=1000$ and $tol=10^{-11}$. BDF1 without projection.}
    \label{tab:num-9}
\end{table}

It is interesting that the projection step $\m_h^{n+1}=\m_h^{n+1}/|\m_h^{n+1}|$ does not affect the iteration performance of GMRES and AGMG, which is shown in \Cref{tab:num-8} and \Cref{tab:num-9}.

% \begin{figure}[htbp]
%     \centering
%     \subfloat[gmres Iter]{\includegraphics[width=0.25\linewidth]{gmres_3D_Iter_v2_v1}}
%     \subfloat[gmres relres]{\includegraphics[width=0.25\linewidth]{gmres_3D_relres_v2_v1}}
%     \hspace{0.1in}
%     \subfloat[AGMG Iter]{\includegraphics[width=0.25\linewidth]{agmg_3D_iter_v2_v1}}
%     \subfloat[AGMG relres]{\includegraphics[width=0.25\linewidth]{agmg_3D_relres_v2_v1}}
%     \caption{Different linear solvers with $\alpha=0$, $N_x=N_y=N_z=30$, $N_t=40$ up to the final time $T=0.1$ for 3D case.}
%     \label{fig:BDF1-8}
% \end{figure}

% \subsection{The simulations with only exchange field}

\subsection{Generalization to different initial conditions}

In this section, we test different initial conditions for 3D as below, 
\begin{itemize}
    \item Previous test with the forcing term. In this test, we set the forcing term to be zero, and take the initial condition not like $\m_0=[\cos(\cos(\pi x)\cos(\pi y)\cos(\pi z))\sin(0),\sin(\cos(\pi x)\cos(\pi y)\cos(\pi z))\sin(0),\cos(0)]^T$, since that the trivial solution will be obtained. So we choose the initial condition as 
    \begin{align*}
        \m_0=[\cos(\cos(\pi x)\cos(\pi y)\cos(\pi z))\sin(0.01),\sin(\cos(\pi x)\cos(\pi y)\cos(\pi z))\sin(0.01),\cos(0.01)]^T,
    \end{align*}
    which we refer to the initial C1 state. The solution profile is presented in \Cref{fig:init-1}. The result visualizes magnetization profiles under C1‑state initialization for the projected 3D micromagnetic test with damping \(\alpha=0.01\). The initial magnetic configuration is displayed in subplots (a)‑(b). At final time \(t=0.1\), both AGMG ((c)‑(d)) and GMRES ((e)‑(f)) produce nearly identical vector‑arrow and color‑mapped magnetization patterns. Visually indistinguishable solution profiles confirm that AGMG yields accurate physical results consistent with standard GMRES, while delivering superior computational efficiency observed in prior numerical tables.
    \begin{figure}[htbp]
    \centering
    \subfloat[arrow, Initial C1 State]{\includegraphics[width=0.35\linewidth]{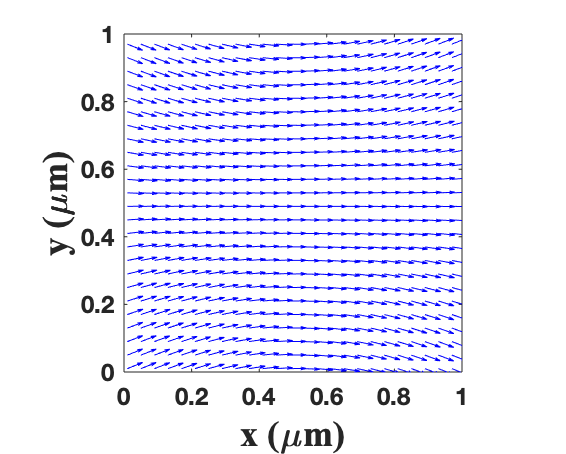}}
    \subfloat[color, Initial C1 State]{\includegraphics[width=0.35\linewidth]{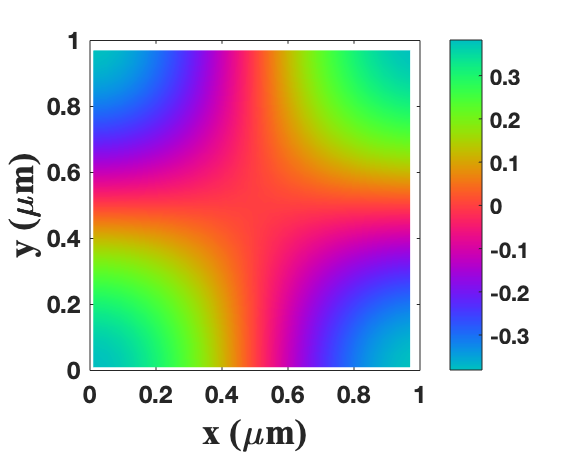}}
    \hspace{0.1in}
     \subfloat[arrow, agmg, $t=0.1$]{\includegraphics[width=0.35\linewidth]{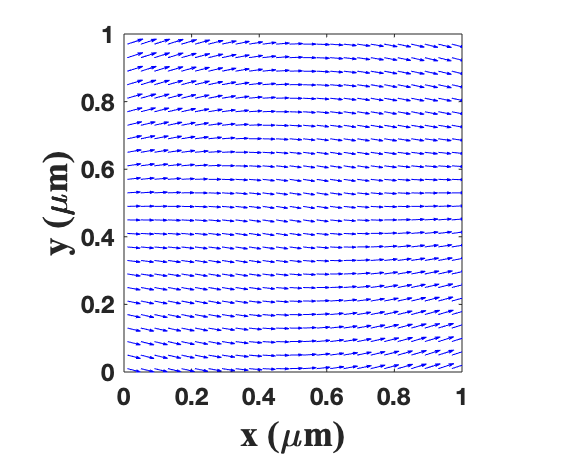}}
     \subfloat[color, agmg, $t=0.1$]{\includegraphics[width=0.35\linewidth]{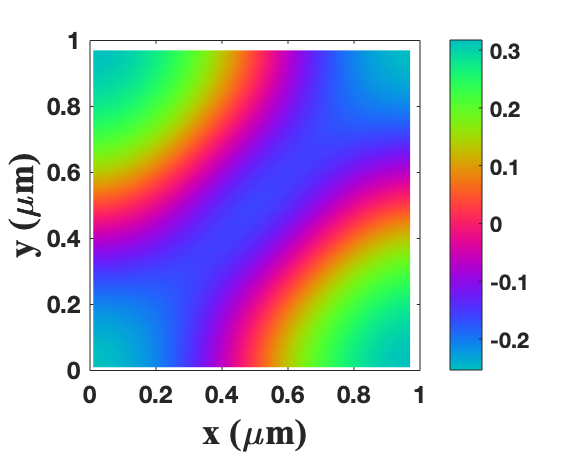}}
     \hspace{0.1in}
     \subfloat[arrow, gmres, $t=0.1$]{\includegraphics[width=0.35\linewidth]{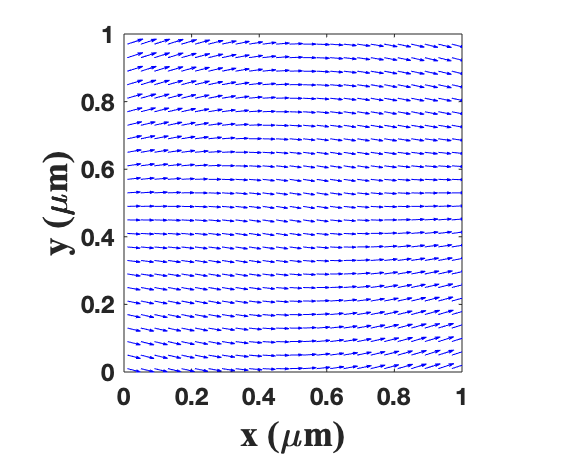}}
     \subfloat[color, gmres, $t=0.1$]{\includegraphics[width=0.35\linewidth]{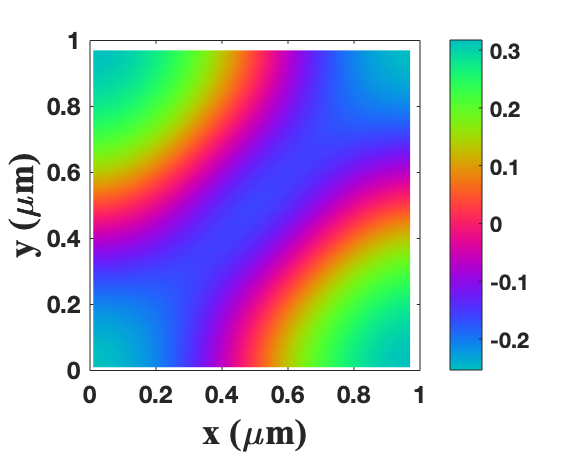}}
    \caption{The solution profiles given C1 state initialization. Using Matlab for the 3D case, with damping $\alpha=0.01$, with projection. $N_x=N_y=50,Nz=4$, $N_t=40$ up to the final time $T=0.1$. Top row with agmg, bottom row with gmres.}
    \label{fig:init-1}
\end{figure}
The comparison is presented in \Cref{tab:init-1}. It presents comparative numerical assessments of GMRES and AGMG for projected BDF1‑based 3D micromagnetic computations initialized from the C1 state with damping \(\alpha=0.01\). Under several mesh‑time discretization combinations, GMRES suffers from excessive iterative cycles and prolonged computing time, even hitting the predefined iteration ceiling and yielding deteriorated residual precision. In sharp contrast, AGMG maintains bounded iteration counts, markedly cuts down computational overhead, and sustains favorable residual convergence across test configurations.
\begin{table}[htbp]
    \centering
    \begin{tabular}{|c|c|c|c|c|c|c|c|}
    \hline
        Method & $N_x$ &$N_y$&$N_z$ & $N_t$& Total wall time (s)&average iter& average relres \\
        \hline
        GMRES &50&50&4 &40&227.780569&3.530250000000000e+02&9.873557363370828e-12  \\ 
         &100&4&4&10&118.253155&1000&4.849567504332029e-06\\
         &50&4&4&100&4.154236&42.850000000000001&9.087021667227475e-12\\
         &10&10&10&100&4.974056&9.449999999999999&7.510565789774815e-12\\
        \hline
        AGMG &50&50&4 &40&20.772403&7.800000000000000&5.879876263337615e-13\\
           &100&4&4&10&1.221880&1.300000000000000&7.325215014967852e-15\\
           &50&4&4&100&3.712899&4&3.663610246651624e-15\\
           &10&10&10&100&6.241831 &8&3.461028498051210e-12\\
        \hline
    \end{tabular}
    \caption{The computational cost, average iteration over time, average relative residual over time using GMRES, AGMG in the 3D case with different $N_x,N_y,N_z, N_t$, and $\alpha=0.01$ given initial C1 state up to the final time $T=0.1$. In the runtime, we set the $maxit=1000$ and $tol=10^{-11}$. BDF1 with projection.}
    \label{tab:init-1}
\end{table}
    \item S state initialization: The initialization is given by
    \[
 \m_0(x)=
\begin{cases}
\bigl(0,\,1,\,0\bigr), & x<\dfrac{x_{\text{max}}}{5}\;\text{or}\; x>\dfrac{4x_{\text{max}}}{5},\\[6pt]
\bigl(1,\,0,\,0\bigr), & \dfrac{x_{\text{max}}}{5}\le x \le \dfrac{4x_{\text{max}}}{5},
\end{cases}
\]
with $x_{max}=1$.
The solution profile is presented in \Cref{fig:init-2}. The solution profile using AGMG and GMRES are comparable.
    \begin{figure}[htbp]
    \centering
    \subfloat[arrow, Initial S State]{\includegraphics[width=0.35\linewidth]{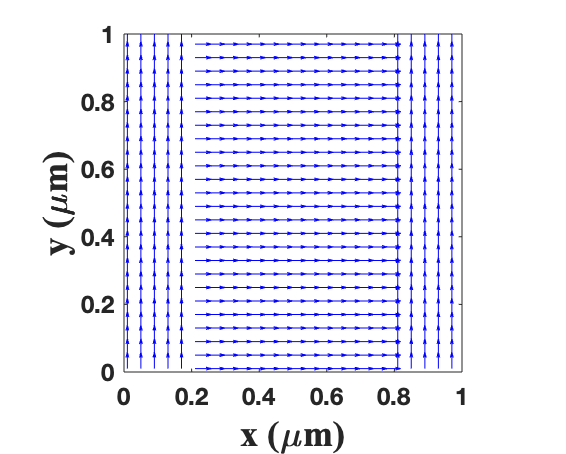}}
    \subfloat[color, Initial S State]{\includegraphics[width=0.35\linewidth]{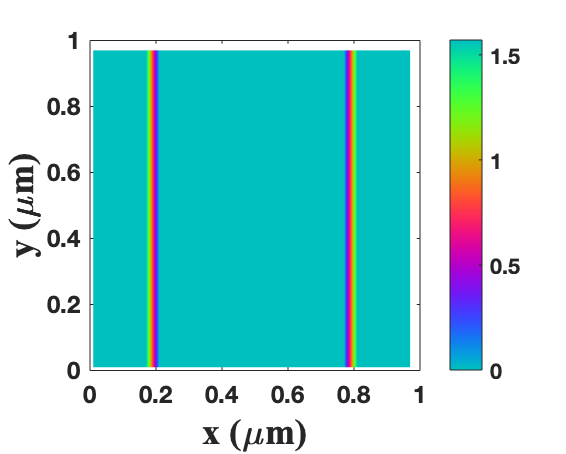}}
    \hspace{0.1in}
     \subfloat[arrow, agmg, $t=0.1$]{\includegraphics[width=0.35\linewidth]{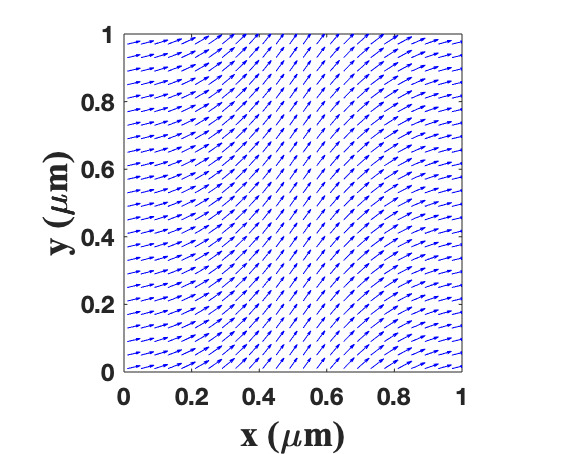}}
     \subfloat[color, agmg, $t=0.1$]{\includegraphics[width=0.35\linewidth]{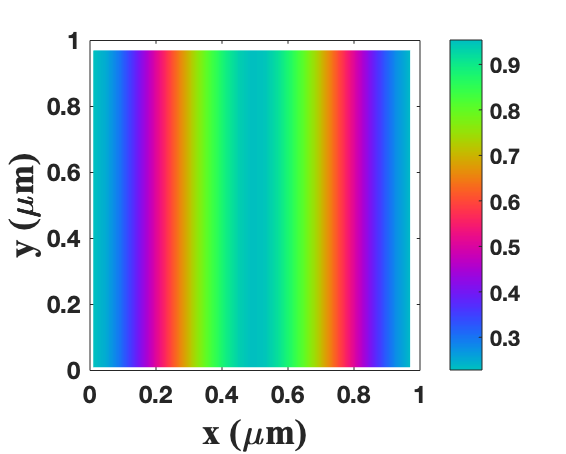}}
     \hspace{0.1in}
     \subfloat[arrow, gmres, $t=0.1$]{\includegraphics[width=0.35\linewidth]{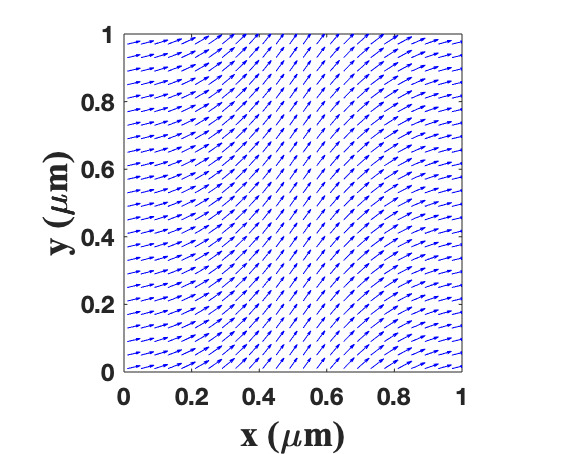}}
     \subfloat[color, gmres, $t=0.1$]{\includegraphics[width=0.35\linewidth]{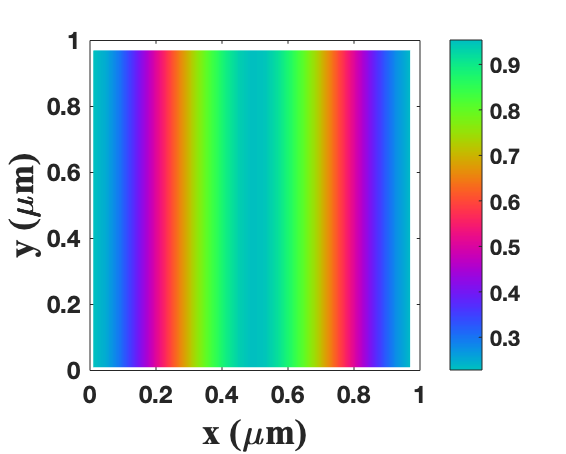}}
    \caption{The solution profiles given S state initialization. Using Matlab for 3D case, with damping $\alpha=0.01$, with projection. $N_x=N_y=50,Nz=4$, $N_t=40$ up to the final time $T=0.1$. Top row with agmg, bottom row with gmres.}
    \label{fig:init-2}
\end{figure}
The comparison is presented in \Cref{tab:init-2}. It quantifies the numerical performance of GMRES and AGMG for projected BDF1 3D micromagnetic simulations starting from the S‑state initial configuration with damping \(\alpha=0.01\). Across diverse spatial‑temporal discretization setups, GMRES exhibits scaling‑dependent performance degradation, consuming substantial runtime and accumulating large iteration counts. By contrast, the AGMG preconditioned solver delivers nearly mesh‑independent iteration numbers, achieves prominent speed‑up factors, and maintains satisfactory residual convergence for all tested discretization scenarios.
\begin{table}[htbp]
    \centering
    \begin{tabular}{|c|c|c|c|c|c|c|c|}
    \hline
        Method & $N_x$ &$N_y$&$N_z$ & $N_t$& Total wall time (s)&average iter& average relres \\
        \hline
        GMRES &50&50&4 &40&147.534499&2.823750000000000e+02&9.169961413394674e-12 \\ 
         &100&4&4&10&111.118406&9.449000000000000e+02&4.906533407786123e-08\\
         &50&4&4&100&5.630164&82.150000000000006&7.491625195093052e-12\\
         &10&10&10&100&4.813901&11.279999999999999&4.067289446664573e-12\\
        \hline
        AGMG &50&50&4 &40&25.854048&1.175000000000000&2.011205351726279e-13\\
           &100&4&4&10&1.521969&1.200000000000000&4.456186620678823e-13\\
           &50&4&4&100&3.957046&1.400000000000000&1.633089675435382e-13\\
           &10&10&10&100&5.800985&11.810000000000000&4.153792750468288e-12\\
        \hline
    \end{tabular}
    \caption{The computational cost, average iteration over time, average relative residual over time using GMRES, AGMG in the 3D case with different $N_x,N_y,N_z, N_t$, and $\alpha=0.01$ given initial S state up to the final time $T=0.1$. In the runtime, we set the $maxit=1000$ and $tol=10^{-11}$. BDF1 with projection.}
    \label{tab:init-2}
\end{table}
    \item C state initialization: The initailization is given by
    \[
\m_0(x)=
\begin{cases}
\bigl(0,\; 1,\; 0\bigr), & x<\dfrac{x_{\text{max}}}{5},\\[6pt]
\bigl(0,-1,\; 0\bigr), & x>\dfrac{4x_{\text{max}}}{5},\\[6pt]
\bigl(1,\; 0,\; 0\bigr), & \dfrac{x_{\text{max}}}{5}\le x \le \dfrac{4x_{\text{max}}}{5}.
\end{cases}
\]
with $x_{max}=1$.
The solution profile is presented in \Cref{fig:init-3}. The solution profiles are comparable.
\begin{figure}[htbp]
    \centering
    \subfloat[arrow, Initial C State]{\includegraphics[width=0.35\linewidth]{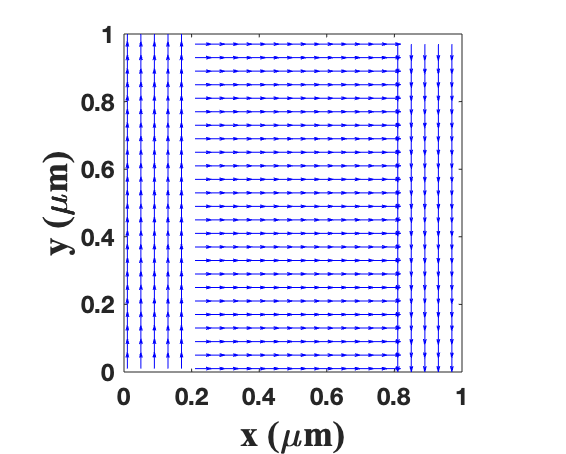}}
    \subfloat[color, Initial C State]{\includegraphics[width=0.35\linewidth]{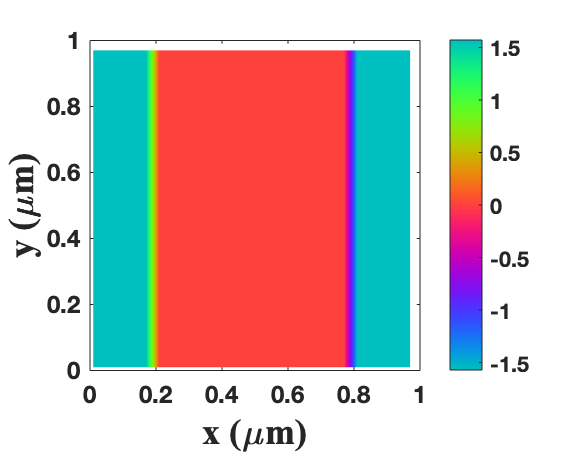}}
    \hspace{0.1in}
     \subfloat[arrow, AGMG, $t=0.1$]{\includegraphics[width=0.35\linewidth]{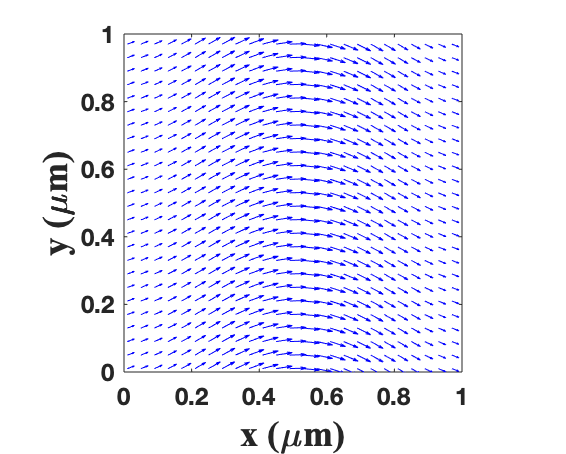}}
     \subfloat[color, AGMG, $t=0.1$]{\includegraphics[width=0.35\linewidth]{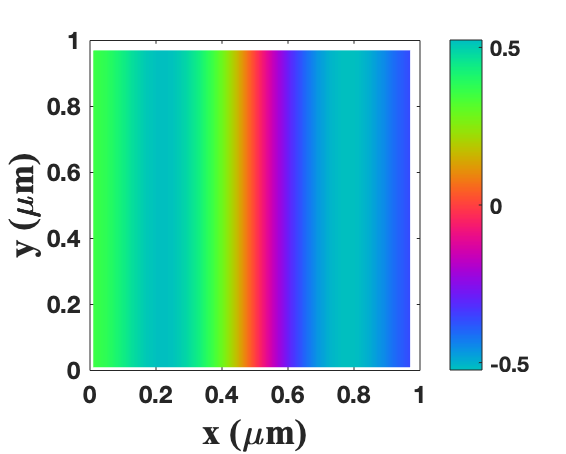}}
     \hspace{0.1in}
     \subfloat[arrow, GMRES, $t=0.1$]{\includegraphics[width=0.35\linewidth]{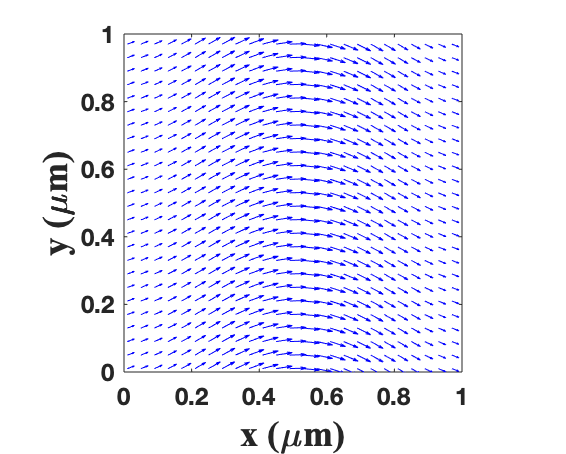}}
     \subfloat[color, GMRES, $t=0.1$]{\includegraphics[width=0.35\linewidth]{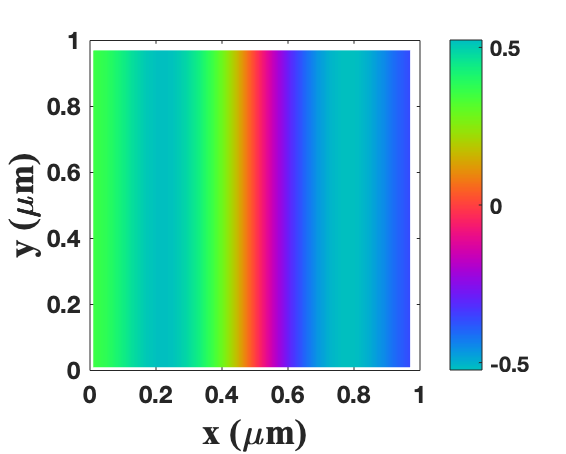}}
    \caption{The solution profiles given C state initialization. Using Matlab for 3D case, with damping $\alpha=0.01$, with projection. $N_x=N_y=50,Nz=4$, $N_t=40$ up to the final time $T=0.1$. Top row with agmg, bottom row with gmres.}
    \label{fig:init-3}
\end{figure}
The comparison is presented in \Cref{tab:init-3}. It assesses GMRES and AGMG for projected BDF1 3D micromagnetic simulations with C‑state initialization. GMRES suffers from high iteration counts and heavy computational burden, whereas AGMG yields prominent acceleration while maintaining favorable residual accuracy across discretizations.
\begin{table}[htbp]
    \centering
    \begin{tabular}{|c|c|c|c|c|c|c|c|}
    \hline
        Method & $N_x$ &$N_y$&$N_z$ & $N_t$& Total wall time (s)&average iter& average relres \\
        \hline
        GMRES &50&50&4 &40&123.555034 &2.531500000000000e+02&9.085076356441198e-12 \\ 
         &100&4&4&10&99.262196&9.023000000000000e+02&9.908446750786463e-12\\
         &50&4&4&100&5.584694&82.469999999999999&7.217054436455491e-12\\
         &10&10&10&100&4.890206&11.869999999999999&4.281400781314147e-12\\
        \hline
        AGMG &50&50&4&40&26.221334&2.300000000000000&8.853025431798482e-13\\
           &100&4&4&10&1.423810 &1.200000000000000&4.265561633465099e-13\\
           &50&4&4&100&4.231381&2.910000000000000&1.754236799563387e-12\\
           &10&10&10&100&5.708775&12.170000000000000&4.510504751616761e-12\\
        \hline
    \end{tabular}
    \caption{The computational cost, average iteration over time, average relative residual over time using GMRES, AGMG in the 3D case with different $N_x,N_y,N_z, N_t$, and $\alpha=0.01$ given initial C state up to the final time $T=0.1$. In the runtime, we set the $maxit=1000$ and $tol=10^{-11}$. BDF1 with projection.}
    \label{tab:init-3}
\end{table}
    \item Flower state initialization: The initialization is given by
    \[
\m_0(x,y)=
\begin{cases}
\bigl(0,\; 1,\; 0\bigr), & x<\dfrac{x_{\text{max}}}{10},\ y<\dfrac{y_{\text{max}}}{2},\\[4pt]
\bigl(0,-1,\; 0\bigr), & x<\dfrac{x_{\text{max}}}{10},\ y>\dfrac{y_{\text{max}}}{2},\\[4pt]
\bigl(0,-1,\; 0\bigr), & x>\dfrac{9x_{\text{max}}}{10},\ y<\dfrac{y_{\text{max}}}{2},\\[4pt]
\bigl(0,\; 1,\; 0\bigr), & x>\dfrac{9x_{\text{max}}}{10},\ y>\dfrac{y_{\text{max}}}{2},\\[4pt]
\bigl(1,\; 0,\; 0\bigr), & \text{otherwise},
\end{cases}
\]
with $x_{max}=1$, and $y_{max}=1$. The solution profile is presented in \Cref{fig:init-3}.
\begin{figure}[htbp]
    \centering
    \subfloat[arrow, Initial C State]{\includegraphics[width=0.35\linewidth]{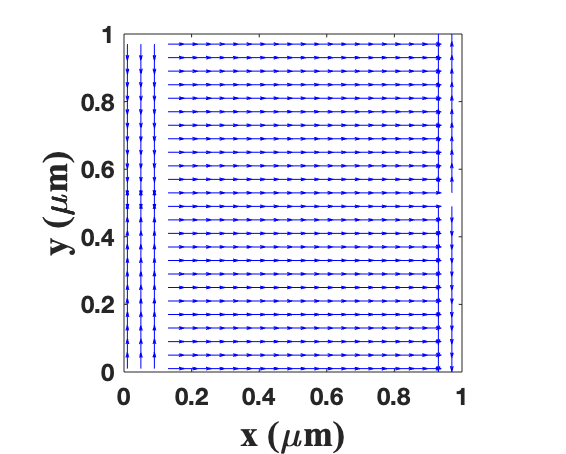}}
    \subfloat[color, Initial C State]{\includegraphics[width=0.35\linewidth]{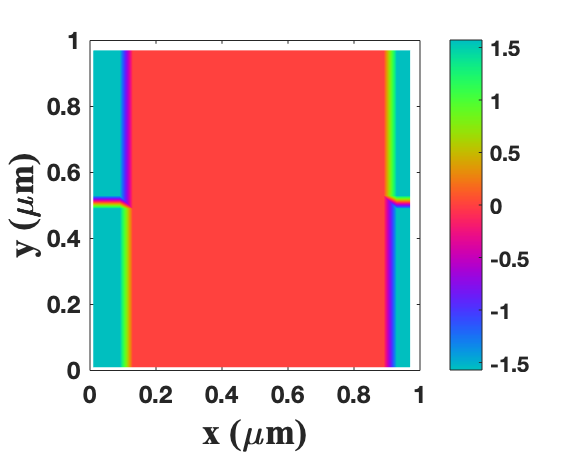}}
    \hspace{0.1in}
     \subfloat[arrow, AGMG, $t=0.1$]{\includegraphics[width=0.35\linewidth]{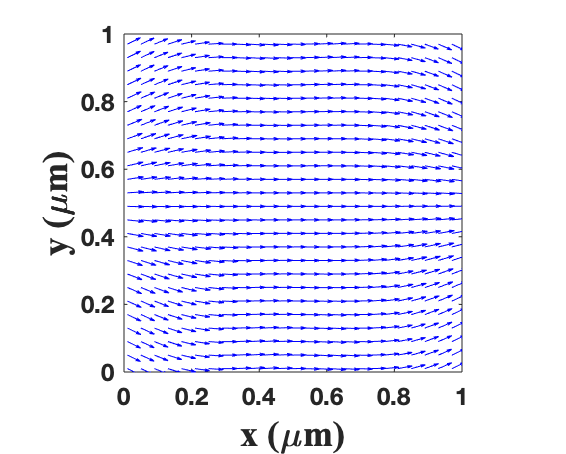}}
     \subfloat[color, AGMG, $t=0.1$]{\includegraphics[width=0.35\linewidth]{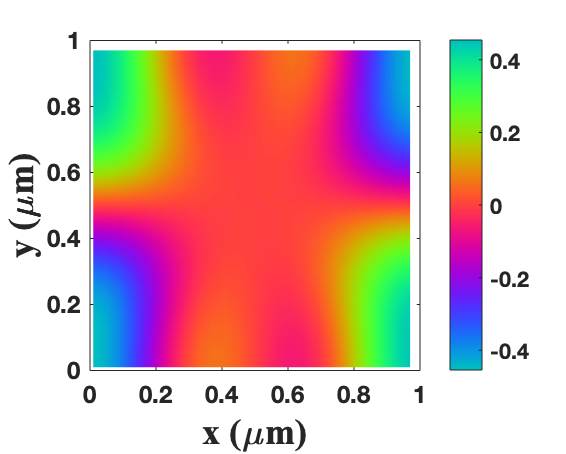}}
     \hspace{0.1in}
     \subfloat[arrow, GMRES, $t=0.1$]{\includegraphics[width=0.35\linewidth]{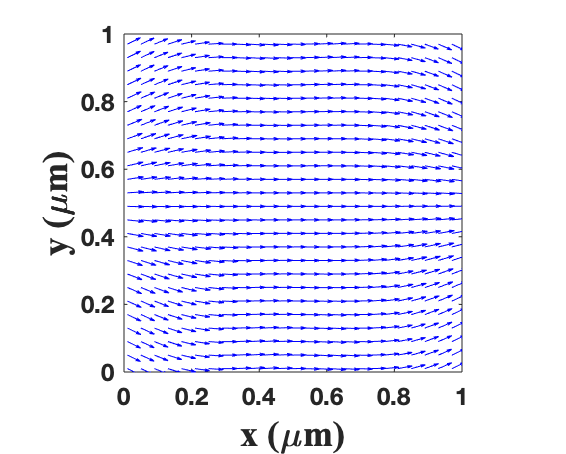}}
     \subfloat[color, GMRES, $t=0.1$]{\includegraphics[width=0.35\linewidth]{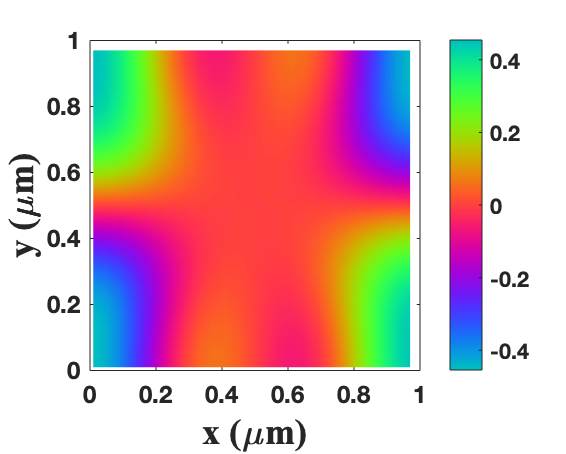}}
    \caption{The solution profiles given Flower state initialization. Using Matlab for 3D case, with damping $\alpha=0.01$, with projection. $N_x=N_y=50,Nz=4$, $N_t=40$ up to the final time $T=0.1$. Top row with agmg, bottom row with gmres.}
    \label{fig:init-3}
\end{figure}
The comparison is presented in \Cref{tab:init-3}.
\begin{table}[htbp]
    \centering
    \begin{tabular}{|c|c|c|c|c|c|c|c|}
    \hline
        Method & $N_x$ &$N_y$&$N_z$ & $N_t$& Total wall time (s)&average iter& average relres \\
        \hline
        GMRES &50&50&4 &40&1297.981502&7.896000000000000e+02&1.146390151299056e-11 \\ 
         &100&4&4&10&158.063117&1000&0.001846593003202\\
         &50&4&4&100&12.013690&1.300300000000000e+02&8.341205024161826e-12\\
         &10&10&10&100&8.734229&19.390000000000001&5.441520081854459e-12\\
        \hline
        AGMG &50&50&4&40&29.671035&7.650000000000000&9.620581330060751e-13\\
           &100&4&4&10&2.997722&1.200000000000000&3.047166910560058e-13\\
           &50&4&4&100&6.055821&4.220000000000000&1.972125479010863e-12\\
           &10&10&10&100&9.418841&11.560000000000000&4.154244734214298e-12\\
        \hline
    \end{tabular}
    \caption{The computational cost, average iteration over time, average relative residual over time using GMRES, AGMG in the 3D case with different $N_x,N_y,N_z, N_t$, and $\alpha=0.01$ given initial Flower state up to the final time $T=0.1$. In the runtime, we set the $maxit=1000$ and $tol=10^{-11}$. BDF1 with projection.}
    \label{tab:init-3}
\end{table}
\end{itemize}

% \subsection{The simulations with only exchange field, anisotropy field, stray field}

% \subsection{The simulations with only exchange field, anisotropy field, stray field, external field}

% \subsection{Norm preserving tests}

% \subsection{Initial conditions and stability tests}

% \subsection{Stability tests}

% \subsection{Hypre lib}

% In our previous simulations, we use PCR-GMRES the solver in Hypre. We specify the mpirun -n with $n=1$. In this section, We still take the simple case with $\alpha=0$ with projection.

% \begin{table}[htbp]
%     \centering
%     \begin{tabular}{c|c}
%     \hline
%         Solver & solver id  \\
%         \hline
%         BoomerAMG & 0 \\
%         \hline
%         ParCSRPCG & 50 \\
%         \hline
%          ParCSRPCG with preconditioner BoomerAMG & 1 \\
%         \hline
%         ParCSRPCG with preconditioner ParaSails & 8 \\
%         \hline
%         \red ParCSRGMRES & 11 \\
%         \hline
%         ParCSRFlexGMRES & 12 \\
%         \hline
%         ParCSRFlexGMRES with preconditioner BoomerAMG & 121 \\
%         \hline
%         ParCSRLGMRES & 13 \\
%         \hline
%         ParCSRBiCGSTAB & 14 \\
%         \hline
%         ParCSRHybrid & 15 \\
%         \hline
%     \end{tabular}
%     \caption{Hypre test for 3D case.}
%     \label{tab:placeholder}
% \end{table}

\subsection{Micromagnetics Simulations}

We apply ParCSRGMRES and AGMG solver to simulate the dynamics of the real ferromagnetic model.

\begin{figure}[htbp]
    \centering
    \subfloat[arrow, Initial S State]{\includegraphics[width=0.35\linewidth]{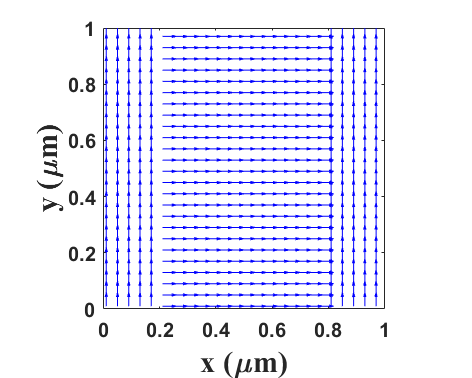}}
    \subfloat[color, Initial S State]{\includegraphics[width=0.35\linewidth]{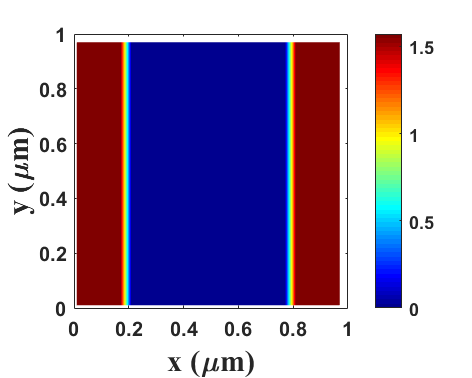}}
    \hspace{0.1in}
     \subfloat[arrow, $t=0.1$]{\includegraphics[width=0.35\linewidth]{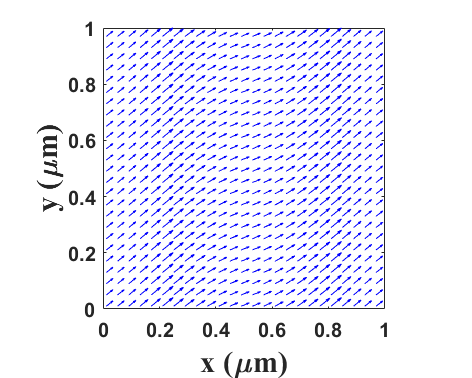}}
     \subfloat[arrow, $t=0.1$]{\includegraphics[width=0.35\linewidth]{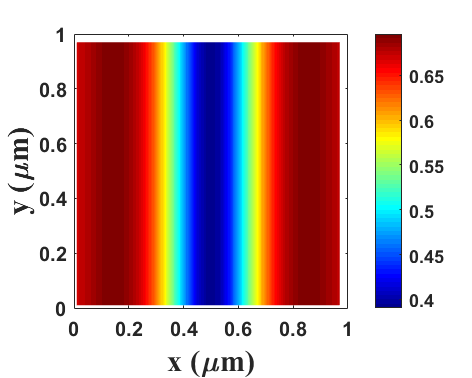}}
    \caption{The solution profile at the final time $T=0.1$ by usage of GMRES (matrix in ParCSR form) for the simulations with the S state initialization.}
    \label{fig:micro-1}
\end{figure}

In our previous simulations, see references \cite{xie2020second,xie2025enhancing}, we use the GMRES (the matrix in ParCSR form) for the dynamics of simulation. We find out it cost a lot of time and slower. If we use the S state initialization. We take the domain $\Omega=[0,1]^3$, we take the number of subdomains $N_x=N_y=50, N_z=4$, and the damping $\alpha=0$ up to the final time $T=0.1$ and $N_t=100$. The results are presented in \Cref{fig:micro-1}.

\section{Conclusions and discussions}
\label{sec:conclusions}

In this work, we develop a fast aggregation-based multigrid method for semi-implicit schemes of the Landau–Lifshitz equation. As the number of grid points increases, we obtain large-scale sparse nonsymmetric linear systems. Conventional iterative solvers such as GMRES suffer from growing iteration counts as the matrix dimension expands. Simulations of magnetic materials require resolving structures of fine characteristic sizes, which necessitates extremely small mesh sizes and results in linear systems of very large dimension. In such scenarios, standard iterative methods perform poorly. Preconditioned GMRES and related algorithms impose stringent requirements on the design and selection of preconditioners. The aggregation-based multigrid method proposed in this paper achieves higher efficiency compared with conventional iterative approaches. Classical multigrid methods are mainly devised for linear systems with symmetric structures, whereas the aggregation-based multigrid presented herein is specially tailored for the target problem. We demonstrate the computational efficiency of the proposed method through numerical tests, and apply it to simulations of magnetic materials. In our tests, the projection step does not affect the performance of AGMG and GMRES iterations. We will investigate fast multigrid methods for semi-implicit schemes of the Maxwell–Landau–Lifshitz equation \cite{bavnas2010efficient,bavnas2008convergent} in our future work. We will develop the parallel version of the current multigrid method for the micromagnetics.

\section*{Data availability}
The data will be made available on reasonable request.

\section*{Conflict of Interest Statement}
The authors have no conflicts of interest to declare. 

\section*{Acknowledgments}
This work is partially supported by the Basic Research Program of Jiangsu Province under Grant BK20250468, and the Research and Development Fund of XJTLU under Grant RDF-24-01-015.

\vspace{1cm}

\bibliographystyle{elsarticle-num-names}
\bibliography{references.bib}

\end{document}